\documentstyle{article}
\begin{document}

\def\R{{\bf R}}
\def\C{{\bf C}}
\def\Z{{\bf Z}}
\def\Q{{\bf Q}}
\def\A{{\cal A}}
\def\O{{\cal O}}
\def\L{{\cal L}}

\vskip2.5mm

\begin{center}
{\bf Secants of Abelian Varieties, Theta Functions, and
Soliton Equations}
\footnote[1]{This work is supported by RFFI, grant 96-01-01889.}
\end{center}

\begin{center}
{\bf I. A. Taimanov}
\end{center}

\vskip2.5mm

\begin{center}
{\bf Introduction}
\end{center}

The main aim of the present article is to overview the
applications of soliton equations to the geometry of Abelian
varieties.

Novikov was the first to indicate a possibility of this
by conjecturing that Jacobi varieties are exactly
principally-polarised Abelian varieties such that the
Kadomtsev--Petviashvili equation integrates in their theta functions.

This conjecture followed a stormy starting period of the development of
finite-zone integration of soliton equations (the Korteweg--de Vries
equation, the Toda lattice, the sine--Gordon equation, etc.) (see
\cite{DMN,TS}) which led to the method of Baker--Akhieser functions.
This method was proposed by Krichever for constructing
finite-zone solutions to the Kadomtsev--Petviashvili equation
and describing commutative rings of ordinary differential
operators of rank $1$ (\cite{Kr1,Kr2}).  The latter results inspired
the Novikov conjecture.  Its proof  by Shiota (\cite{Shi1}) led to
solution of one of the oldest and prominent problems of algebraic
geometry, the Riemann--Schottky problem.

At the same time, Arbarello and De Concini (\cite{AC1}) showed that
the Novikov conjecture interweaves tightly with the ideas of Gunning
who proposed to describe Jacobi varieties as
admitting sufficiently many trisecants. The last condition amounts
to validity of special theta function identities
(\cite{G}). In their research Arbarello and De Concini used
the observation of Mumford that many soliton equations (the
Korteweg--de Vries, Kadomtsev--Petviashvili, and sine--Gordon
equations) are ``hidden'' in the Fay trisecant formula (\cite{Mm2}).

The present article gives a survey of these papers on the
Riemann--Schottky problem as well as of the papers on the analogue of
the Novikov conjecture for Prym varieties.
The analogue is abstracted further because the
Veselov--Novikov, Landau--Lifschitz and BKP equations integrated in
Prym theta functions are also ``hidden'' in the Fay and
Beauville--Debarre quadrisecant formulae.

Chapter 1 contains a brief introduction to the analytic theory of
theta functions (see also \cite{GH,Mm1}) which can be considered as
self-contained if supplemented with information on
cohomologies with coefficients in sheaves, for instance, in the amount
of Chern's book \cite{C}.

Chapter 2 provides necessary information on the theta functions of
Jacobi and Prym (the detailed proofs are exposed, for instance, in
\cite{Fay1,GH}) and the inference of the trisecant and quadrisecant
formulae.

Chapters 3 and 4 are devoted to application of soliton equations to the
theory of Jacobi and Prym varieties.
Therewith necessary information on finite-zone solutions to non-linear
equations and on Baker--Akhieser functions is exposed in brief in \S 8.
This information is given in detail in the survey articles
on finite-zone theory \cite{DMN,Kr2,KN2,Dubr1} (see also
\cite{BBME,Kr89}).  We hope that the present survey completes them
in the part related to application of finite-zone theory to the
geometry of Abelian varieties.

\newpage

\begin{center}
{\bf Chapter 1. Abelian varieties and theta functions}
\end{center}

\vskip2.5mm

\begin{center}
{\bf \S 1. A condition for a complex torus to be
algebraic}
\end{center}

\vskip2.5mm

{\bf 1.1. Algebraic varieties.}

A complex manifold is called a
{\it projective algebraic variety}
if it is embeddable into the complex projective space $\C P^n$
as the set of zeros of a system of homogeneous polynomials.

By definition, $\C P^n$ is the quotient space of
$\C ^{n+1} \setminus \{0\}$ for the following action of
$\C^*= \C \setminus \{0\}$:
$$
(z^0,z^1,\dots,z^n) \rightarrow
(\lambda z^0, \lambda z^1, \dots, \lambda z^n), \lambda \in \C^*.
$$
This space is endowed with the Fubini--Study metric
written in terms of the homogeneous coordinates
$(z^0 : \dots : z^n)$ as
$$
\frac{(\sum_j z^j {\bar z}^j) \cdot (\sum_k dz^k d{\bar z}^k)
- (\sum_j z^j d{\bar z}^j)\cdot (\sum_k {\bar z}^k dz^k)}{(\sum_k z^k
{\bar z}^k)^2}.
$$
Assign to this Hermitian metric the uniquely-determined form
$$
\omega =
\frac{\sqrt{-1}}{2\pi} \cdot \frac{(\sum_j z^j {\bar z}^j) \cdot
\sum_k dz^k \wedge d{\bar z}^k - (\sum_j z^j d{\bar z}^j)
\wedge (\sum_k {\bar z}^k dz^k)}{(\sum_k z^k {\bar z}^k)^2}.
$$
Straightforward computations show the validity of
the following:

1) $\omega$ is a closed form;

2) $[\omega]$ is an integer cohomology class:
$[\omega] \in H^2(\C P^n; \Z)$.

Complex manifolds with the above properties
fall into special classes:

1) a compact manifold $M$ is
called a {\it K\"ahler manifold}
if there exists a Hermitian
metric $h_{ij} dz^i d{\bar z}^j$ on $M$ such that the form
\begin{equation}
\omega=
\frac{\sqrt{-1}}{2\pi} h_{ij} dz^i \wedge d{\bar z}^j
\label{1}
\end{equation}
associated with this metric is closed.
The metric $h_{ij}$ is said to be {\it K\"ahler};

2) a K\"ahler manifold $M$ is called
a {\it Hodge manifold} if the form $\omega$ associated with the
K\"ahler metric on $M$ represents an integer cohomology class,
$[\omega] \in H^2(M ; \Z)$.
This form $\omega$ is said to be {\it Hodge}.

Unless stated otherwise, repeated upper and lower indices in the same
expression imply summation in (\ref{1}) and in the
sequel.

Give an invariant definition of $\omega$.
It suffices to do that locally, in the tangent space at a point.
Decompose the Hermitian form
$$
(v,w) = h_{ij}v^i{\bar w}^j
$$
into the sum of its real and imaginary parts:
\begin{equation}
H(v,w) = H_R(v,w) + \sqrt{-1}H_I(v,w)
\label{2}
\end{equation}
where the real forms $H_R$ and $H_I$ are $\R$-bi-linear.
Since $H$ is Hermitian, $H_R$ is symmetric and
positive definite whereas $H_I$ is skew-symmetric.
Since $H(\sqrt{-1}v,w) = \sqrt{-1}H(v,w)$, we have
\begin{equation}
H_R(v,w) = H_I(\sqrt{-1}v,w).
\label{3}
\end{equation}
Hence the Hermitian metric is determined by the skew-symmetric form
$H_I$ alone.  Given $x \in M$,
put
\begin{equation}
\omega(x)(v,w) = -\frac{1}{\pi}
H_I(x)(v,w)
\label{4}
\end{equation}
where $v,w \in T_x M$.

The equivalence of (\ref{1}) and (\ref{4})
is clear for diagonal metrics. The general case reduces to this one by
change of a linear basis.

{\bf 1.2. Theorems of Kodaira and Chow.}

Clearly the properties that a manifold is a
K\"ahler or Hodge manifold are inherited by submanifolds.  For
instance, let $Y$ be a submanifold of a K\"ahler manifold $X$. Then the
K\"ahler metric on $X$ induces a K\"ahler metric on $Y$. This also
holds in the case of Hodge manifolds.  Hence every projective algebraic
variety is a Hodge manifold. Kodaira proved that the converse statement
is also valid.

{\bf 1.2.1. Kodaira Theorem.}
{\sl A complex manifold is projective algebraic if and
only if it is a Hodge manifold.}

Proof of this theorem consists in
embedding a Hodge manifold into a projective space of
sufficiently high dimension and applying the following

{\bf 1.2.2. Chow Theorem.}
{\sl An analytic subset $X \subset \C P^n$ is distinguished in
$\C P^n$ as the set of zeros of a system of homogeneous polynomials.}

For complex tori the equivalence of the property of being a Hodge
manifold and that of being an algebraic variety was established by
Lefschetz (see \S 4).

{\bf 1.3. A criterion for a complex torus to be a Hodge manifold.}

A complex torus $\C^n/\Lambda$, with $\Lambda$ a lattice of
rank $2n$ in $\C^n$, is called an {\it Abelian variety} (or an
{\it Abelian torus}) if it is a projective algebraic variety. By
the Lefschetz theorem, this is equivalent to the
property to be a Hodge torus.

{\bf 1.3.1. Riemann Criterion.}
{\sl A complex torus $M=\C^n/\Lambda$ is a Hodge manifold if
and only if there exists a complex linear basis
for $\C^n$ such that in this basis the
lattice $\Lambda$ is written as
\begin{equation}
\Lambda = \Delta_\delta N_1 + \Omega N_2, \ \
N_1, N_2 \in {\Z}^n,
\label{5}
\end{equation}
where

1) $\Delta_{\delta}$ stands for the diagonal $(n \times n)$-matrix
$diag(\delta_1,\dots,\delta_n)$ with integer entries,
$\delta_j > 0$, and $\delta_k$ divides $\delta_{k+1}$ for every
$k=1,\dots,(n-1)$;

2) the matrix $\Omega$ is symmetric;

3) the matrix $Im\, \Omega$ is positive definite.}

Proof of the Riemann criterion.

Suppose that the complex torus $M =
\C^n/\Lambda$ is a Hodge manifold.
Take a Hermitian metric
$h_{ij}dz^i d{\bar z}^j$ on $\C^n$
invariant under translations by the vectors of $\Lambda$
and inducing a metric that defines a Hodge structure on $M$.

{\bf 1.3.2.} {\sl Let $x^1,\dots,x^{2n}$ be $\R$-linear
coordinates on $\R^{2n}=\C^n$ and let $\omega =
\omega_{ij}dx^i \wedge dx^j$ be a closed $2$-form on
$\C^n/\Lambda$. Define an averaging operator  as follows
$$
\omega_{ij} \rightarrow {\tilde \omega}_{ij}= \frac {\int_{\Pi}
\omega_{ij}(z^1,\dots,z^n) dz^1 \wedge \dots \wedge dz^n} {\int_{\Pi}
dz^1 \wedge \dots \wedge dz^n}
$$
where $\Pi$ is the fundamental domain of
$\Lambda$. Then the averaging operator transforms
$\omega$ into a $2$-form with constant coefficients which is
cohomologous to $\omega$. Therewith Hodge forms are transformed into
Hodge forms.}

The proof is clear.

{\bf 1.3.3.} {\sl Every closed $2$-form $\omega$ on
$M = \C^n/\Lambda$ defines a skew-symmetric form $Q$ on
$\Lambda$. The value of $Q$ at a pair $(u,v)$ is equal to
the integral of $\omega$ over the $2$-torus spanned by
$u$ and $v$. If $[\omega] \in H^2(M;\Z)$ then $Q$
is integer-valued.  If $\omega$ is a Hodge form then  $Q$ is
non-degenerate.}

Let $v,w \in \Lambda$ and let $\Pi_{v,w}$ be the parallelogram in
$\C^n$ spanned by $v$ and $w$.  This parallelogram projects onto
the $2$-cycle $T_{v,w}$ in $M$. Define $Q$ as follows
$$
Q(v,w) = \int_{T_{v,w}} \omega.
$$
By (\ref{1}), if $\omega$ is a Hodge form
then $\omega^n$ is proportional to the volume form, and we infer that
$Q$ is non-degenerate. This proves Proposition 1.3.3.

{\bf 1.3.4. Frobenius Criterion.}
{\sl Let $M = \C^n/\Lambda$ be a complex torus with $\Lambda$
a lattice of rank $2n$. The torus $M$ is a Hodge manifold if and
only if there exists a skew-symmetric real $\R$-bi-linear form $Q$
on $\C^n \times \C^n$ such that

1) the form $\langle v,w \rangle = Q(\sqrt{-1}v,w)$ is symmetric and
positive definite;

2) $Q$ is integer-valued on $\Lambda \times \Lambda$.}

Proof of the Frobenius criterion.

Let $M$ be a Hodge torus. Take a Hermitian metric with constant
coefficients on $M$ generating a Hodge form.
By (\ref{4}) and 1.3.3, the form $\omega$ associated with
this metric satisfies conditions 1 and 2.

Let $Q$ satisfies the hypotheses of 1.3.4.
Then the Hermitian metric $\langle v,w \rangle = \pi (Q(\sqrt{-1}v,w)
+ \sqrt{-1}Q(v,w))$ generates the associated Hodge form.

The Frobenius criterion is established.

{\bf 1.3.5.}
{\sl Let $Q$ be a non-degenerate integer skew-symmetric
form on $\Lambda=\Z^{2n}$.
Then in a suitable basis
$\lambda_1,\dots,\lambda_{2n}$ the form $Q$ is written as
\begin{equation}
Q = \left( \begin{array}{cc} 0 & \Delta_{\delta} \\ -\Delta_{\delta}
& 0 \end{array} \right)
\label{6}
\end{equation}
where $\Delta_{\delta}$ stands for the diagonal integer matrix
$diag(\delta_1,\dots,\delta_n)$, $\delta_k>0$, and $\delta_j$ divides
$\delta_{j+1}$ for every $j=1,\dots,n-1$.  An $n$-tuple
$\delta_k$ satisfying these conditions is an invariant of $Q$.}

We prove this proposition.

Take the minimal value
$\delta_1 =Q(\lambda_1,\lambda_{n+1})$ among all possible positive
values of the form
$Q(\lambda,\lambda')$. Denote by $\Lambda'$ the
orthogonal complement to the sublattice
$\Z\{\lambda_1,\lambda_{n+1}\}$ in
$\Lambda$. We have
$$
\lambda +
\frac{Q(\lambda,\lambda_1)}{\delta_1}\lambda_{n+1} -
\frac{Q(\lambda,\lambda_{n+1})}{\delta_1}\lambda_1 \in \Lambda'.
$$
for $\lambda \in \Lambda$.
Hence $\Lambda' = \Z^{2n-2}$ and the restriction of
$Q$ onto $\Lambda'$ satisfies the hypotheses of Proposition
1.3.5. Using this procedure successively we thus construct a
basis $\lambda_1,\dots,\lambda_{2n}$ for $\Lambda$ such that in
this basis $Q$ takes the shape (\ref{6}).

Prove that $\delta_1$ divides $\delta_2$
by way of contradiction.
Let $\delta$ be the greatest common
divisor of $\delta_1$ and $\delta_2$ such that
$\delta < \delta_1$. Then there exist integers $k$ and $n$ such
that $\delta=k\delta_1+m\delta_2$. However,
$$
Q(k\lambda_1 +
m\lambda_2, \lambda_{n+1} + \lambda_{n+2}) = \delta < \delta_1
$$
which contradicts the choice of $\delta_1$.
Hence $\delta_1$ divides $\delta_2$.  We similarly
prove that $\delta_j$ divides
$\delta_{j+1}$ for every $j, 1 \leq j \leq n-1$.

It follows from the construction of
$\{\delta_j\}$ that this $n$-tuple is an invariant of $Q$.
This proves Proposition 1.3.5.

{\bf 1.3.6.}
{\sl Let $\omega$ be the Hodge form associated with a Hermitian
metric with constant coefficients on a Hodge torus
$M=\C^n/\Lambda$ and let $\lambda_1,\dots,\lambda_{2n}$ be a basis for
the lattice $\Lambda$ in which $Q$ takes the shape (\ref{6}).  Then the
set of the vectors $\lambda_1, \dots, \lambda_n$ is a complex
linear basis for $\C^n$.}

Proof of Proposition 1.3.6.
Since the coefficients of a Hermitian form are constant,
the form $Q$, $\R$-linearly-extended onto $\C^n$, coincides with
$\omega$ and determines a Hermitian metric by (\ref{2}--\ref{4}).

Consider the succession of vectors
$$
v_1 = \lambda_1,\ v_2 = \lambda_2
- \frac{Q(\sqrt{-1}\lambda_2,v_1)}{Q(\sqrt{-1}v_1,v_1)}v_1, \dots, \
v_k = \lambda_k - \sum_{j<k}
\frac{Q(\sqrt{-1}\lambda_k,v_j)}{Q(\sqrt{-1}v_j,v_j)}v_j, \
\dots .
$$
By (\ref{3}--\ref{4}), all the expressions
$Q(\sqrt{-1}\lambda_k,v_j)/Q(\sqrt{-1}v_j,v_j)$
are real.
It follows from (\ref{2}) and (\ref{6}) that
$(v_j,v_k) = 0$ for $j \neq k$.
Since $Q(v_k,\lambda_{n+k}) \neq 0$, the vectors
$v_1,\dots,v_n$ do not vanish and form a complex linear
basis for $\C^n$. Now Proposition 1.3.6 follows from the
coincidence of the linear spans of
$(v_1,\dots,v_n)$ and $(\lambda_1,\dots,\lambda_n)$ over $\C$.

Now we turn directly to proving the Riemann criterion.

Let $M=\C^n/\Lambda$ be a Hodge torus.
By Propositions 1.3.2--1.3.6,
there exist a basis ${\cal V}_0 =
(\lambda_1,\dots,\lambda_{2n})$ for $\Lambda$ and a
Hodge form $\omega$ with constant coefficients such that
in this basis $\omega$ is written as
\begin{equation}
\omega = \sum_{k=1}^n
\delta_k dy^k \wedge dy^{n+k}.
\label{7}
\end{equation}
Take a
complex basis ${\cal V}_{\C}$ and a real basis ${\cal V}_{\R}$ for
$\C^n$ as follows
$$
{\cal V}_{\C} =
(\frac{1}{\delta_1}\lambda_1,\dots,\frac{1}{\delta_n}\lambda_n),
$$
\begin{equation}
{\cal V}_{\R} =
(\frac{1}{\delta_1} \lambda_1,\dots,\frac{1}{\delta_n}\lambda_n,
\lambda_{n+1},\dots,\lambda_{2n}).
\label{8}
\end{equation}

Then $\Lambda$ takes the shape (\ref{5})
and we have the formula relating
${\cal V}_{\R}$ and ${\cal V}_{\C}$
\begin{equation}
\left(
\begin{array}{c}
z^1 \\ \cdots \\ z^n
\end{array}
\right)
=
\left(I_n, \ \Omega
\right)
\left(
\begin{array}{c}
x^1 \\ \cdots \\ x^{2n}
\end{array}
\right),
\label{9}
\end{equation}
where $I_n$ is the identity $(n \times n)$-matrix.

By (\ref{7}--\ref{9}), in the basis (\ref{8})
$\omega$ is written as
\begin{equation}
\omega = \sum_{k=1}^n dx^k
\wedge dx^{n+k}.
\label{10}
\end{equation}

Since $\omega$ is associated with the Hodge metric $h_{ij}$
by (\ref{1}), we infer from (\ref{10}) that
\begin{equation}
\omega=
\frac{\sqrt{-1}}{2\pi}h_{ij}(dx^i + \Omega^i_kdx^{n+k}) \wedge
(dx^j + \Omega^j_ldx^{n+l}) =  \sum_{m=1}^n dx^m \wedge dx^{n+m}.
\label{11}
\end{equation}
Since $h_{ij}$ is symmetric, we conclude that
(\ref{11}) is equivalent to
\begin{equation}
\frac{\sqrt{-1}}{2\pi}h_{kj}({\bar \Omega}^j_l - \Omega^j_l)
dx^k \wedge dx^{n+l} = \sum_{m=1}^n dx^m \wedge dx^{n+m}
\label{12}
\end{equation}
and
\begin{equation}
T_{lm} = T_{ml} \ \ \ {\mbox where}\ \ \
T_{lm} = h_{kj}\Omega^k_m{\bar \Omega}^j_l.
\label{13}
\end{equation}
It follows from (\ref{12}) that
\begin{equation}
T_{lm} = h_{kj}\Omega^k_l\Omega^j_m +
\frac{2\pi}{\sqrt{-1}}\Omega^m_l.
\label{14}
\end{equation}

Now we conclude that

1) the matrix $Im \ \Omega$ is positive definite
(this follows from (\ref{12}));

2) the matrix $\Omega$ is symmetric: $\Omega^l_m =
\Omega^m_l$ (this follows from (\ref{13}) and (\ref{14})).

Hence we prove that if $\C^n/\Lambda$
is a Hodge manifold then it admits a form given in the statement
of the Riemann criterion.

If $\C^n/\Lambda$
takes this shape, then, by (\ref{11}--\ref{14}),
the Hermitian form $h_{ij}dz^i \wedge d{\bar z}^j$ defined by
(\ref{12}) is a Hodge form and, consequently, this torus is a Hodge
torus.

The Riemann criterion is established.

{\bf 1.4. Polarizations and moduli of Abelian varieties.}

In this subsection of \S 1 we identify Hodge and Abelian tori.
The complete proof of the coincidence of these families is
given in \S 4.

A Hodge torus, or an Abelian variety,
$M=\C^n/\Lambda$ is {\it polarised}
with {\it polarization} $[\omega] \in H^2(M;\Z)$
if the cohomology class of a Hodge form $\omega$ is fixed.

By 1.3.3 and 1.3.5, the polarization $[\omega]$
uniquely defines the non-degenerate integer
skew-symmetric $2$-form
$$
Q : H^1(M) \times H^1(M) \rightarrow \Z
$$
taking the shape (\ref{6}) in a suitable basis for $\Lambda$.  Moreover
the $n$-tuple of positive integers $\delta_1,\dots,\delta_n$ is
uniquely determined by this form.  This $n$-tuple
$\delta_1,\dots,\delta_n$ is called a {\it polarization type}.  An
Abelian variety is {\it principally-polarised} if $\delta_1 = \dots =
\delta_n = 1$.

Two polarised Abelian varieties
$(M_1,[\omega_1])$ and $(M_2,[\omega_2])$ are
{\it equivalent} if there exists a bi-holomorphic
mapping $f :  M_1 \rightarrow M_2$ such that
$f^*([\omega_2]) = [\omega_1]$.  Given pairs
($n, \delta)$, the sets of equivalence classes of
$n$-dimensional Abelian varieties with polarization of type
$\delta$ are  complex manifolds
$\A_n^{\delta}$, the {\it moduli spaces of Abelian varieties}.

Denote by $\A_n$ the moduli spaces of $n$-dimensional
principally-polarised Abelian varieties.

The moduli spaces $\A_n^{\delta}$ are as follows.

Denote by ${\cal H}_n$ the {\it Siegel upper half-plane}
formed by symmetric $(n \times n)$-matrices with positive
definite imaginary parts.
It follows from the proof of the Riemann criterion and
Propositions 1.3.3 and 1.3.5  that every $n$-dimensional Abelian
variety with polarization of type $\delta$ is equivalent to the torus
$\C^n/\Lambda$ with $\Lambda$ of the shape (\ref{7}) and $\Omega
\in {\cal H}_n$. It is clear that two polarised Abelian varieties $M_1
= \C^n/\{\Delta_\delta\Z^n + \Omega_1\Z^n\}$ and $M_2 =
\C^n/\{\Delta_\delta\Z^n + \Omega_2\Z^n\}$ are equivalent if and only
if their lattices are adjoint by a transformation, in $GL(2n,\Z)$,
preserving (\ref{6}).

Precise the statement above. The group
$SL(\delta,\Z)$ is a subgroup of $GL(2n,\Z)$ formed by elements
$g$ such that
\begin{equation}
\left(
\begin{array}{cc} A & B \\ C & D \end{array} \right) \cdot \left(
\begin{array}{cc} 0 & \Delta_\delta \\ -\Delta_\delta & 0 \end{array}
\right) \cdot \left( \begin{array}{cc} A^* & C^* \\ B^* & D^*
\end{array}
\right)
=
\left(
\begin{array}{cc}
0 & \Delta_\delta \\
-\Delta_\delta & 0
\end{array}
\right)
\label{15}
\end{equation}
where
$$
g =
\left(
\begin{array}{cc}
A & B \\
C & D
\end{array}
\right)
$$
and $A, B, C$, and $D$ are integer $(n \times n)$-matrices.

The action of $Sp(\delta,\Z)$ on ${\cal H}_n$ generated by
automorphisms of lattices is as follows
\begin{equation}
\Omega \rightarrow (A\Omega + B\Delta_\delta)
\cdot (\Delta_\delta^{-1}C\Omega +
\Delta_\delta^{-1}D\Delta_\delta)^{-1}.
\label{16}
\end{equation}
We obtain the usual modular action of
$SL(2,\Z)$ on the complex upper half-plane
for $n=1$ and $\delta=1$.

{\bf 1.4.1.}
{\sl The moduli space of $n$-dimensional Abelian varieties
with polarization of type
$\delta$ is the quotient space of the Siegel upper half-plane
${\cal H}_n$ for the action (\ref{16}) of $Sp(\delta,\Z)$:
$$
\A^{\delta}_n = {\cal H}_n / Sp(\delta,\Z).
$$}

An {\it isogeny} of Abelian varieties $M$ and $M'$
is a holomorphic finite-sheeted covering
$\varphi :  M \rightarrow M'$.

{\bf 1.4.2.}
{\sl Let $M$ be a polarised Abelian variety. Then there exist
a principally-polarised Abelian variety $M'$ and an isogeny
$$
\varphi : M \rightarrow M'
$$
with degree $\Delta = \delta_1
\dots \delta_n$ where $(\delta_1,\dots,\delta_n)$ is the
polarization type of $M$.}

We prove Proposition 1.4.2. Let $M$ be an Abelian variety
$\C^n /\Lambda$ of the shape (\ref{7}) and let
$(\lambda_1,\dots,\lambda_{2n})$ be a basis for
$\Lambda$. Assume that a Hodge form $\omega$ takes the shape
(\ref{7}) in this basis.

Consider the lattice $\Lambda'$ spanned by the vectors of
the basis (\ref{8}).  This lattice contains
$\Lambda$ as a sublattice with index $\Delta$.
The projection
$$
\varphi :  \C^n /\Lambda \rightarrow
\C^n/\Lambda'
$$
is an isogeny of degree $\Delta$. Choose a form
$\omega'$ such that  $\varphi^*(\omega') = \omega$.
Then, by (\ref{10}),  it defines the principal polarization on $M'$.
This proves Proposition 1.4.2.

\vskip2.5mm

\begin{center}
{\bf \S 2. Line bundles on complex tori}
\end{center}

\vskip2.5mm

{\bf 2.1. Line bundles, Picard group, Chern classes, and
positive bundles.}

A {\it line bundle} $L$ on a complex manifold
$M$ is a vector bundle with a fibre $\C$ and holomorphic
coordinate transformation functions, i.e.,

1) a projection $p : L \rightarrow M$ is defined;

2) there exists an open cover ${\cal U}= \{U_{\alpha}\}$
of $M$ such that a bundle is trivial over
each element of this family:
$$
U_{\alpha} \times \C \stackrel{\xi_{\alpha}}{\approx}
p^{-1}(U_{\alpha}), \ \
p\cdot\xi_{\alpha}(z,v)=z \in U_{\alpha} \subset M;
$$

3) the coordinate transformations
$$
g_{\alpha\beta}(z) =  \xi_{\alpha}^{-1}\xi_{\beta}|_{p^{-1}(z)} :
U_{\alpha} \cap U_{\beta} \rightarrow \C^* = \C
\setminus \{0\}
$$
are holomorphic functions.

By definition, the coordinate transformations
satisfy the following conditions
$$
g_{\alpha\beta} \cdot g_{\beta\alpha} = 1, \ \ \
g_{\alpha\beta} \cdot g_{\beta\gamma} \cdot g_{\gamma\alpha} = 1.
$$
By these conditions, the family of functions
$g_{\alpha\beta}$ defining the line bundle uniquely
is a $1$-cocycle in $C^1({\cal U}, \O^*)$.
The sheaf $\O^*$ is defined as follows.
Let $\O^*(U)$ be a multiplicative group of non-vanishing
holomorphic functions on $U \subset M$.  It is said that
two cocycles $g$ and $g'$ are {\it cohomologous} if there exists a
family of sections $f_{\alpha} : U_{\alpha}
\rightarrow \O^*(U_{\alpha})$ such that
$g'_{\alpha\beta} = (f_{\alpha}/f_{\beta})g_{\alpha\beta}$.

Two line bundles are called {\it equivalent} if the cocycles
$g$ and $g'$ corresponding to these bundles are cohomologous.
Therewith the functions $f_{\alpha}$ deform the trivializations
as follows $\xi_{\alpha} \rightarrow \xi_{\alpha}/f_{\alpha}$.

The set of equivalence classes of line bundles on $M$ is
a commutative group, the {\it Picard group $Pic(M)$},
under the operations of
tensor product
\begin{equation}
L , L'
\rightarrow L \otimes L', \ \ \ g_{\alpha\beta},g'_{\alpha\beta}
\rightarrow g_{\alpha\beta}\cdot g'_{\alpha\beta} \label{17}
\end{equation}
and passing to the dual bundle
\begin{equation}
L \rightarrow L^*, \ \ \
g_{\alpha\beta} \rightarrow g^{-1}_{\alpha\beta}.
\label{18}
\end{equation}

We thus obtain the following proposition.

{\bf 2.1.1. Structure of the Picard group.}
{\sl The Picard group of line bundles on $M$ is isomorphic
to the first cohomology group of $M$ with coefficients in the sheaf
of germs of non-vanishing holomorphic functions:
$Pic(M) = H^1(M;\O^*)$.}

To a line bundle $L$ one assigns its Chern class $c_1(L)$.
Define it in terms of the coordinate transformations
$g_{\alpha\beta}$.  Let ${\cal U}$ be
an open cover of $M$ by simply-connected
domains in $M$ such that $L$ is trivial over each domain of this
cover and let $\{g_{\alpha\beta}\}$ be the corresponding coordinate
transformations.  Define functions
\begin{equation}
h_{\alpha\beta} = \frac{1}{2\pi\sqrt{-1}}
\log{(g_{\alpha\beta})},\ \ \log{(g_{\alpha\beta})} = -
\log{(g_{\beta\alpha})}.
\label{19}
\end{equation}
and construct a cocycle
$z_{\alpha\beta\gamma} \in C^2({\cal U};\Z)$
from these functions as follows
\begin{equation}
z_{\alpha\beta\gamma} = h_{\alpha\beta} + h_{\beta\gamma} -
h_{\gamma\alpha}.
\label{20}
\end{equation}
The cohomology class $c_1(L) \in H^2(M;\Z)$
realised by this cocycle is called the {\it (first) Chern
class} of $L$.

Since the choice of other branches of the logarithm
in the definition of $h_{\alpha\beta}$
($\log{(g_{\alpha\beta})} \rightarrow \log{(g_{\alpha\beta})} +
2\pi\sqrt{-1}k_{\alpha\beta}, k_{\alpha\beta} \in \Z$) is
equivalent to the change of the cocycle
$z_{\alpha\beta\gamma}$ to the cohomologous cocycle $z \rightarrow z
+ dk$, this definition is correct.

By (\ref{20}), the Chern class $c_1$ is a homomorphism of
the Picard group into $H^2(M;\Z)$:
\begin{equation}
c_1 (L \otimes L') =  c_1(L) + c_1(L'), \ \
c_1(L^*) = - c_1(L).
\label{21}
\end{equation}

Recall the invariant definition of $c_1$.
Consider the exact sequence of sheaves
$$
0 \rightarrow \Z \rightarrow \O
\stackrel{\exp}{\rightarrow} \O^* \rightarrow 0
$$
where $\Z$ is the constant sheaf, $\O$ is the sheaf of germs of
holomorphic functions, and $\exp$ is the exponential homomorphism
$f \rightarrow \exp{(2\pi\sqrt{-1}f)}$.
Assign to this sequence the following exact cohomology sequence
\begin{equation}
\dots \rightarrow H^k(M;\Z) \rightarrow H^k(M;\O)
\rightarrow H^k(M;\O^*) \rightarrow H^{k+1}(M;\Z) \rightarrow \dots.
\label{22}
\end{equation}
Consider the boundary homomorphism (\ref{22}) for $k=1$:
$$
\delta :  H^1(M;\O^*) = Pic(M) \rightarrow H^2(M;\Z).
$$
The first Chern class is the image of $L$ under this homomorphism:
$$
\delta(L) = c_1(L).
$$

A line bundle $L$ is called {\it positive}
if its first Chern class is realised by a Hodge form $\omega$,
namely, $c_1(L) =  [\omega]$.

{\bf 2.2. Chern classes of line bundles on tori.}

Let $M = \C^n/\Lambda$ be a complex torus and let
$\pi :  \C^n \rightarrow M$ be the universal covering
of $M$. For every line bundle $L$ on $M$
there exists a pull-back
$\pi^*(L)$, a line bundle on $\C^n$.
By the ${\bar \partial}$-lemma of Poincare,
$H^1(\C^n;\O^*) = 0$ (see, for instance, \cite{GH}).
Hence the bundle $\pi^*(L)$ is trivial,
$\pi^*(L) = \C^n \times \C$.

Fix a trivialization of $\pi^*(L)$.
Construct the automorphisms of a fibre
$$
e_{\lambda}(z) : \C \rightarrow \C, \ \
e_{\lambda}(z) \in \C^*
$$
by comparing the trivializations of $\pi^*(L)$ at the points
$z$ and $z+\lambda$ with $\lambda \in \Lambda$.
The set of all such automorphisms
is a set of
{\it multiplicators} of $L$, the family of functions
$\{e_{\lambda} \in \O^*(U)\}_{\lambda \in \Lambda}$ such that
\begin{equation}
e_{\lambda'}(z+\lambda)e_{\lambda}(z) =
e_{\lambda}(z+\lambda')e_{\lambda'}(z) = e_{\lambda+\lambda'}(z).
\label{23}
\end{equation}

Every set of functions
$\{e_{\lambda} \in \O^*(U)\}_{\lambda \in \Lambda}$
satisfying (\ref{23}) defines a line bundle
which is obtained from $\C^n \times \C$ by factorization
under the action of the group $\Lambda$
$$
(z,v) \rightarrow (z+\lambda,e_{\lambda}(z) \cdot v),
\ \ z \in \C^n, v \in \C.
$$

Two sets of multiplicators $e_{\lambda}$ and
${\tilde e}_{\lambda}$ define the same bundle if and only if
there exists a holomorphic function
$\varphi :  \C \rightarrow \C^*$ such that
\begin{equation}
{\tilde e}_{\lambda}(z) =
e_{\lambda}(z)\varphi(z+\lambda)\varphi^{-1}(z).
\label{24}
\end{equation}
The function $\varphi$ changes a trivialization of $\pi^*(L)$.

We find a formula which expresses the Chern class in terms of
multiplicators of a bundle. Fix a basis
$\lambda_1,\dots,\lambda_{2n}$ for the lattice
$\Lambda$ and denote by $x^1,\dots,x^{2n}$
the real coordinates on $\C^n$ referred to this basis.
Consider the cover of $M$ by the sets
$$
U_{\lambda} = \{ |x^k - \lambda^k| \leq \frac{3}{4} \}
$$
where $\lambda = (\lambda^1,\dots,\lambda^{2n}) \in \Lambda$.
By the definition of multiplicators , the coordinate
transformations are written as
\begin{equation}
g_{\alpha\beta}(z) =
e_{\alpha-\beta}(z+\beta).
\label{25}
\end{equation}
Define functions $f_{\alpha}(z)$ as follows
\begin{equation}
e_{\alpha}(z) = e^{2\pi\sqrt{-1}f_{\alpha}(z)}.
\label{26}
\end{equation}

We show how to construct from the cocycle (\ref{20})
a $2$-form realising it.
The nerve $N({\cal U})$ of the minimal subcover of the
cover $U_{\alpha}$ is homeomorphic to $M$
and its cohomologies with coefficients in $\Z$ coincide with
$H^*(M;\Z)$.  The formula (\ref{20}) gives the value of the cocycle
$z$ at a two-dimensional simplex
$(\alpha,\beta,\gamma) \in N({\cal U})$.

Take vectors $\tau,\sigma \in \Lambda$ and span over them the
parallelogram $\Pi_{\tau,\sigma}$ which is projected onto the torus
$T_{\tau,\sigma} \subset M$.
Given values of $z$ at two-dimensional simplexes,
we define the value of $c_1$ at this torus.

By (\ref{20}) and (\ref{24}--\ref{25}),
\begin{equation}
c_1([T_{\tau,\sigma}]) =
f_{\sigma}(z) + f_{\tau}(z+\sigma)
- f_{\tau}(z) - f_{\sigma}(z+\tau).
\label{27}
\end{equation}

{\bf 2.2.1.}
{\sl The formula (\ref{27}) expresses the Chern class $c_1$
of a line bundle on $M=\C^n/\Lambda$.  The Chern class
defines  an integer skew-symmetric $2$-form $Q$ on $\Lambda \times
\Lambda$ by $Q(\tau,\sigma) =
c_1([T_{\tau,\sigma}])$.  The form $\omega$ realising
$c_1(L)$, i.e., $[\omega] = c_1(L)$, is as follows
$$
\omega = \frac{1}{2}Q(\lambda_k,\lambda_l)dx^k \wedge dx^l.
$$}

The following proposition interprets the Chern class as
a class realised by the curvature form of a line bundle (\cite{C}).

{\it $(p,q)$-form} is a linear combination of differential forms
of the type $fdz^{i_1} \wedge
\dots \wedge dz^{i_p} \wedge d{\bar z}^{j_1} \wedge \dots \wedge
d{\bar z}^{j_q}$.  For a compact complex manifold $M$
we distinguish  $H^2_{1,1}(M;\Z)$,
a subgroup of $H^2(M;\Z)$, formed by elements
realised by $(1,1)$-forms.

{\bf 2.2.2.}
{\sl The first Chern class of a line bundle $L$ on a compact
complex manifold $M$ lies in $H^2_{1,1}(M;\Z)$, $c_1(L) \in
H^2_{1,1}(M;\Z)$.}

We prove Proposition 2.2.2 by applying the de Rham theorem of
comparison  between the sheaf cohomologies and de Rham's cohomologies.
The equivalence of (\ref{27}) and (\ref{29}) for tori is checked by
straightforward computations.

Take a Hermitian metric on the line bundle.
This means that we choose positive functions
$a_{\alpha} :  U_{\alpha} \rightarrow \R$,
the metric tensor, such that
\begin{equation}
a_{\alpha} \cdot |g_{\alpha\beta}|^2=a_{\beta} \
\ {\mbox {in}} \ \ U_{\alpha} \cap U_{\beta}.
\label{28}
\end{equation}

Let
$d = d' + d''$  be the decomposition of the
operator of inner derivative $d$ into its holomorphic and
anti-holomorphic parts where $d'$ and $d''$ transform
a $(p,q)$-form into  a $(p+1,q)$-form and a $(p,q+1)$-form. 
It is clear that $d'^2 = d''^2=0$.
By applying the operator  $d'$ to the logarithms of both sides of
(\ref{28}) and considering that the coordinate
transformations are holomorphic, obtain
$$
d' \log{(a_{\alpha})} + d'
\log{(g_{\alpha\beta})} = d' \log{(a_{\beta})}.
$$
Since $dd'\log{(g_{\alpha\beta})} = 0$, the form
$dd'\log{(a_{\alpha})}$ is correctly defined everywhere on $M$.

By comparing the de Rham cohomologies and
spectral cohomologies and (\ref{20}),
obtain
\begin{equation}
c_1(L)
= [\frac{1}{2\pi\sqrt{-1}} d''d'log{(a_{\alpha})}].
\label{29}
\end{equation}

This proves Proposition 2.2.2.

{\bf 2.3. Construction of line bundles
with prescribed Chern classes on tori.}

{\bf 2.3.1. Structure of the group $H^2_{1,1}(M;\Z)$.}
{\sl Let $M=\C^n/\Lambda$ be a complex torus. Then the group
$H^2_{1,1}(M;\Z)$ is isomorphic to the group of
skew-symmetric $\R$-bi-linear $2$-forms $Q$ on $\C^n
\times \C^n$ integer-valued on
$\Lambda \times \Lambda$ and
satisfying the following condition
\begin{equation}
Q(u,v) = Q(\sqrt{-1}u,\sqrt{-1}v).
\label{30}
\end{equation}}

By Propositions 1.3.2 and 1.3.3, it suffices to consider
$2$-forms $\omega \in H^2(M;\Z)$ with constant coefficients.
Assign to every $2$-form $\omega$ an integer
skew-symmetric bi-linear form $Q$ on $\Lambda \times \Lambda$.
Extend $Q$ to a $\R$-bi-linear form on
$\C^n \times \C^n$, preserving the previous notation. This
form coincides with $\omega$.

Put $Q = \alpha_{ij}dz^i \wedge dz^j + \beta{ij}dz^i \wedge
d{\bar z}^j + \gamma_{ij} d{\bar z}^i \wedge d{\bar z}^j$.
It is clear that $Q$ satisfies (\ref{30}) if and only if
$\alpha_{ij}=\gamma_{ij}=0$ which is equivalent to $\omega \in
H^2_{1,1}(M;\Z)$. This proves Proposition 2.3.1.

Let $\omega$ be a $2$-form on
$M=\C^n/\Lambda$ realising an element of $H^2_{1,1}(M;\Z)$ and
let $Q$ be associated with $\omega$ (see Proposition 1.3.3).

By (\ref{30}), the $\R$-bi-linear form
$$
H(u,v) = -\pi (Q(\sqrt{-1}u,v)+\sqrt{-1}Q(u,v))
$$
is Hermitian,
$$
H(u,v) = {\overline{H(v,u)}}, \ \
H(\alpha u,v) = \alpha H(u,v) ,
\alpha \in \C.
$$
If $\omega$ is a Hodge form with constant
coefficients then $H$ is the K\"ahler metric associated with
$\omega$.

Define functions $e_{\lambda}(z), \lambda
\in \Lambda, z \in \C^n$ as follows:
\begin{equation}
e_{\lambda}(z) = \alpha(\lambda)\exp{(H(z,\lambda) +
\frac{1}{2}H(\lambda,\lambda))}.
\label{31}
\end{equation}

{\bf 2.3.2.}
{\sl The set of functions (\ref{31}) is a set of multiplicators of a
line bundle $L$ on $M=\C^n/\Lambda$ if and only if the equality
\begin{equation}
\frac{\alpha(\lambda+\lambda')}{\alpha(\lambda)\alpha(\lambda')}
= \exp{(\pi\sqrt{-1}Q(\lambda,\lambda'))}
\ \ 
{\mbox {for \ every}}\
\lambda,\lambda' \in \Lambda
\label{32}
\end{equation}
holds. If (\ref{32}) is valid, then $\omega$ is the first Chern class
of $L$.}

Proposition 2.3.2 is derived from (\ref{23}) and (\ref{26}) by
straightforward computations.

The following proposition is a corollary of 2.3.2.

{\bf 2.3.3.}
{\sl For every element $[\omega] \in
H^2_{1,1}(M;\Z)$ there exists a line bundle
$L$ on $M=\C^n/\Lambda$ defined by multiplicators
(\ref{31}) such that $c_1(L) = [\omega]$.}

In order to prove this proposition it suffices to show
that (\ref{32}) are compatible for every skew-symmetric
form $Q$ integer-valued on $\Lambda \times \Lambda$.

As in 1.3.5, we construct a basis
$\lambda_1,\dots,\lambda_{2n}$ for $\Lambda$ such that
$Q$ takes the shape (\ref{6}) in this basis.
For every subspace $V_k= \R\{\lambda_k,\lambda_{n+k}\}$ we
define a function $\alpha_k$ as follows: $\alpha_k(x\lambda_k + y
\lambda_{n+k}) = \exp{(\pi\sqrt{-1}\delta_kxy)}$.  Now
define a function
$\alpha(v)=\alpha_1(v_1)\dots \alpha_n(v_n)$ where $v = v_1 +
\dots v_n , v_k \in V_k$.  We are left to notice that 
$\alpha$ satisfies (\ref{32}).
The proof of Proposition 2.3.3
is complete.

{\bf 2.4. Line bundles with $c_1=0$: the group
$Pic^0(M)$.}

{\bf 2.4.1.}
{\sl If the Chern class of a line bundle on a complex torus $M$
vanishes, then this bundle can be defined by constant
multiplicators.}

We prove this proposition. Assume $c_1(L)=0$.
By the definition of $c_1$ and (\ref{22}),
$[L] \in H^1(M;\O^*)$ lies in the image of the homomorphism
$H^1(M;\O) \rightarrow H^1(M;\O^*)$.
Indeed, the cocycle $\zeta \in H^1(M;\O)$ is defined by the sections
$\log{(g_{\alpha\beta})}$ and is transformed into the cocycle
$g_{\alpha\beta}$ defining $L$.

For every compact complex manifold

1) the group $H^1(M;\O)$ is isomorphic to the subgroup of
$H^1(M;\C)$ formed by the cohomology classes of closed $(0,1)$-forms;

2) the natural homomorphism $H^1(M;\C)
\rightarrow H^1(M;\O)$ generated by the embedding of sheaves
$\C \rightarrow \O$ is onto.
In terms of the de Rham cohomologies,
this homomorphism is a projection
$$
\xi = \xi^{1,0} + \xi^{0,1} \rightarrow
\xi^{0,1}
$$
where $\xi = \xi^{1,0} + \xi^{0,1}$ is the decomposition of a closed
form into the sum of a $(1,0)$-form and a $(0,1)$-form (\cite{C}).

Therefore $\zeta$ is cohomologous to the cocycle
$\log{(g_{\alpha\beta})}$ with constant coefficients.
This means that the line bundle is defined by constant
multiplicators. This proves Proposition 2.4.1.

Consider line bundles with constant multiplicators.
By Propositions 2.3.2  and 2.4.1, these are exactly line bundles with
$c_1=0$.

{\bf 2.4.2.}
{\sl The equivalence classes of line bundles with $c_1=0$ on
a torus $M$ are in a one-one correspondence with
points of the dual torus ${\hat M}= Pic^0(M)$.}

Let $L$ be a line bundle defined by the constant multiplicators
$e_{\lambda}$.  It is easy to construct a function
$\varphi(z)$, of the shape $\varphi(z) = \exp{(a_kz^k)}$, such that
$|\varphi(\lambda)| = |e_{\lambda}|^{-1}$.
By (\ref{24}), the constant multiplicators
${\tilde e}_{\lambda}(z) =
e_{\lambda}\varphi(z+\lambda)\varphi(z)^{-1}$ define the same
bundle $L$ and $|{\tilde e}_{\lambda}|=1$.
Such sets of multiplicators are in a one-one correspondence
with the homomorphisms
\begin{equation}
\psi :  \Lambda \rightarrow \C_1 = \{z \in \C :  |z|=1\}.
\label{33}
\end{equation}
This correspondence is defined by
\begin{equation}
\psi(\lambda) =
e_{\lambda}.
\label{34}
\end{equation}
The set of homomorphisms (\ref{33}) is in a one-one
correspondence with the torus
${\hat M} = Hom(\Lambda,\C_1)$.

Two distinct homomorphisms defines two distinct bundles.
Prove this by way of contradiction.
Assume that $\psi_1$ and $\psi_2$ define equivalent bundles.
Then, by (\ref{24}), there exists an entire function
$\varphi(z)$ such that
$\psi_1(\lambda) =
\psi_2(\lambda)\varphi(z+\lambda)\varphi^{-1}(z)$
for any $z \in \C^n$ and $\lambda \in \Lambda$.
Since
$|\psi_1(\lambda)| = |\psi_2(\lambda)| =1$, the function
$\varphi(z)$ is bounded and thus, by the Liouville theorem, it is
constant. Hence $\psi_1=\psi_2$.
The proof of Proposition 2.4.2 is complete.

By Proposition 2.4.2 and (\ref{34}),
the subgroup $Pic^0(M)$, formed by equivalence classes of bundles
with $c_1=0$,  coincides with ${\hat M} = Hom(\Lambda,\C_1)$ as a
group also. Since $H^0(M;\O^*)=0$, the other definition of
$Pic^0(M)$ follows from the exact sequence (24).

{\bf 2.4.3.}
{\sl $Pic^0(M) = H^1(M;\O) / H^1(M;\Z)$.}

The following theorem is derived from Propositions 2.3.2,
2.3.3, and 2.4.2.

{\bf 2.4.4. Appel--Humbert Theorem.}
{\sl 1) Every line bundle on the complex torus $M$
can be defined by multiplicators (\ref{31});
2) The sequence
\begin{equation}
0 \rightarrow {\hat M} = Pic^0(M) \rightarrow Pic(M)
\rightarrow H^2_{1,1}(M;\Z) \rightarrow 0,
\label{35}
\end{equation}
where the first homomorphism is an embedding and the second
assigns to a bundle its Chern class, is exact.}

{\bf 2.5. Positive line bundles on tori and their
sections.}

Let $M=\C^n/\Lambda$ be a torus endowed with a Hodge form $\omega$.
By Propositions 1.3.5 and 1.3.6, there exists a basis
${\cal V} = (\lambda_1,\dots,\lambda_{2n})$ for the lattice
$\Lambda$ such that

1) in this basis
$\omega$ takes the shape (\ref{7}), with $\delta_j \in \Z, \delta_j
> 0 $;

2) the set of vectors
$\lambda_1/\delta_1,\dots,\lambda_n/\delta_n$ is a complex basis for
$\C^n$  (denote by $z^1,\dots,z^n$ the coordinates referred to this
basis);

3) in terms of the coordinates $(z^1,\dots,z^n)$ the lattice
$\Lambda$ is written as $\Delta_{\delta} N_1 + \Omega N_2,
N_i \in \Z^n$ where the matrix $\Omega$ is symmetric and
positive definite.

By the Appel--Humbert theorem, the moduli space
${\cal L}(M,[\omega])$ of line bundles on
$M$ with $c_1=[\omega]$ is isomorphic to $Pic^0(M)$.
Indeed, if $L,L' \in {\cal L}(M,[\omega])$ there exists a unique
bundle $L'' \in Pic^0(M)$ such that $L' = L \otimes L''$.

However, $M$ acts on ${\cal L}(M,[\omega])$
by translations as follows. Put $L \in {\cal
L}(M,[\omega])$. Let $\tau_{\mu} :  M \rightarrow M$ be the
translation $z \rightarrow z+\mu$. Then the bundle
$\tau_{\mu}^*(L)$ is defined. In terms of multiplicators this
action is simply as follows
$e_{\lambda}(z) \rightarrow e'_{\lambda}(z) =
e_{\lambda}(z+\mu)$.

Let $L \in {\cal L}(M,[\omega])$ be defined
by the following multiplicators
\begin{equation}
e_{\lambda_k} = 1, \ \
e_{\lambda_{n+k}} = \exp{(-2\pi\sqrt{-1}z^k)}, \ \
k=1,\dots,n.
\label{36}
\end{equation}
Thus this bundle lies in ${\cal L}(M,[\omega])$.
Then the bundle
$\tau_{\mu}^*(L) \otimes L^* \in {\cal L}(M,0)$
is defined by the constant multiplicators
\begin{equation}
{\tilde e}^{(\mu)}_{\lambda_k} = 1, \ \
{\tilde e}^{(\mu)}_{\lambda_{n+k}} = \exp{(-2\pi\sqrt{-1}\mu^k)},
\ \
k=1,\dots,n,
\label{37}
\end{equation}
and by straightforward computations we derive from (\ref{37})
the following proposition.

{\bf 2.5.1.}
{\sl Assume that $L \in {\cal L}(M,[\omega])$. Then the mapping
\begin{equation}
\tau : M \rightarrow Pic^0(M), \ \
\mu \stackrel{\tau}{\rightarrow} \tau_{\mu}^*(L) \otimes L^* ,
\label{38}
\end{equation}
is a homomorphism onto the group $Pic^0(M)$. The kernel of this
homomorphism is isomorphic to the sublattice spanned by
$\{\delta_k^{-1}\lambda_k, \delta_k^{-1}\lambda_{n+k}\}$.
Every bundle $L \in {\cal L}(M,0)$ is written as
$\tau_{\mu}^*(L) \otimes L^*$.}

Denote by $\O(L)$ the sheaf of germs of holomorphic sections
of $L$. Global sections are in one-one correspondence with
entire functions $f(z)$ satisfying the periodicity conditions
\begin{equation}
f(z+\lambda)=e_{\lambda}(z)f(z)
\label{39}
\end{equation}
where
$e_{\lambda}$ are the multiplicators defining $L$.
The set of all global sections is the vector space $H^0(M;\O(L))$.

{\bf 2.5.2.}
{\sl Assume that a bundle $L$ is positive and
$c_1=[\omega]$. Then
$$
h^0(L) = \dim_{\C} H^0(M;\O(L)) =
\delta_1 \cdot \dots \cdot \delta_n
$$
where
$(\delta_1,\dots,\delta_n)$ is the polarization type of $\omega$.}

Proof of Proposition 2.5.2.

Since the translation by $\mu$ is homotopic to the identity
and the bundles in ${\cal L}(M,[\omega])$
are translates of one another,
it suffices to find $h^0(L)$ for an arbitrary bundle $L
\in {\cal L}(M,[\omega])$.  Consider the bundle with multiplicators
(\ref{36}).  Since $e_{\lambda_k} \equiv 1$, every section $f$
expands into the Fourier series in
$\exp{(2\pi\sqrt{-1}l_k\delta_k^{-1}z^k)}$:
\begin{equation}
f(z) =
\sum_{l \in \Z^n} a_l \exp{(2\pi\sqrt{-1}\sum_k l^k \delta_k^{-1}
z^k)}.
\label{40}
\end{equation}
By comparing the Fourier series
expansion of $f(z+\lambda_{n+k})$ and the Fourier series
expansion of $f(z)$ and considering (\ref{39}),
we derive
\begin{equation}
a_{l+d_k} = \exp{(2\pi\sqrt{-1} \langle l,
\Delta_{\delta}^{-1}\lambda_{n+k}
\rangle)}a_l
\label{41}
\end{equation}
where $\langle u,v \rangle = \sum_k u^k v^k $ and $d_k
\in \Z^n, \ d_k^k = \delta_k, d_k^j = 0$ for $j \neq k$.

By (\ref{41}), $f$ is completely determined by
$a_l$ for $l \in \Pi_{\delta} = \{0 \leq l^k < \delta_k\}$.
Thus we derive  $h^0(L) \leq \delta_1 \cdot
\dots \cdot \delta_n$.

Consider functions $\theta_{{\bar l}}, {{\bar
l}} \in \Pi_{\delta}$ of the shape (\ref{40}--\ref{41}) such that
\begin{equation}
a_l = \cases{0 , & for $l \in \Pi_{\delta}$ ,
\cr 1, & for $l = {{\bar l}}$.}
\label{42}
\end{equation}
By (\ref{40}--\ref{41}),
the series (\ref{40}) is uniquely defined for every
${{\bar l}} \in \Pi_{\delta}$.
All these series define entire functions.
We explain this in \S 3 for the theta function
of a principally-polarised Hodge torus.

It is clear that the functions $\theta_{{\bar
l}}$ are linearly independent as sections of $\O(L)$.
Hence, $h^0(L) = \delta_1 \cdot \dots \cdot
\delta_n$. This proves  Proposition 2.5.2.

The following proposition follows from Propositions 2.5.1
and 2.5.2.

{\bf 2.5.3.}
{\sl Let $M$ be a Hodge torus with polarization $[\omega]$.
Then the following statements are equivalent:

1) $[\omega]$ is the principal polarization;

2) the homomorphism (\ref{38})  $\tau : M \rightarrow {\hat M}$
is an isomorphism;

3) $h^0(L) = 1$ where $c_1(L) = [\omega]$.}

\vskip2.5mm

\begin{center}
{\bf \S 3. Theta functions}
\end{center}

\vskip2.5mm

{\bf 3.1. The theta function of a principally-polarised Hodge
torus.}

Every Hodge torus $M$ with the principal polarization
$[\omega]$ is represented as follows

1) $M=\C^n/\Lambda$ and
$\omega = \sum_k dx^k \wedge dx^{n+k}$
in the basis $(\lambda_1,\dots,\lambda_{2n})$ for
$\Lambda$;

2) the lattice $\Lambda$ is written as
$\{\Z^n + \Omega \Z^n\}$ in the basis
$(\lambda_1,\dots,\lambda_n)$ for  $\C^n$;

3) the $(n \times n)$-matrix $\Omega$ is symmetric and its
imaginary part $Im\,\Omega$ is positive definite.

This follows from Propositions 1.3.5 and
1.3.6 and the Riemann criterion.
In the sequel we mean by a principally-polarised Hodge  torus a
complex torus represented in this shape.

Let $L$ be the bundle defined by the following multiplicators
\begin{equation}
e_{\lambda_k}
\equiv 1, \ \ e_{\lambda_{n+k}} =
\exp{(-\pi\sqrt{-1}\Omega_{kk}-2\pi\sqrt{-1}z^k)}
\label{43}
\end{equation}
where  $\Omega_{ij} = \Omega^i_j$ and
$\lambda_{n+k} = \Omega e_k, \  k=1,\dots,n$.

By (\ref{27}), $c_1(L)=[\omega]$ and every bundle with the same
Chern class is a translate of $L$ (Proposition 2.5.1).
It thus suffices to construct the global sections
of this bundle only. By Proposition 2.5.2, the space of these
sections is one-dimensional.

Consider the formal series
\begin{equation}
\theta(z,\Omega) = \sum_{m \in \Z^n}
\exp{(\pi\sqrt{-1}\langle \Omega m,m\rangle + 2\pi\sqrt{-1} \langle
m,z \rangle)}.
\label{44}
\end{equation}
Since $Im\ \Omega$ is positive definite, for every compact subset
$U \subset \C^n$ there exists a constant $c(U)$ such that
$$
|\exp{(\pi\sqrt{-1}\langle \Omega m,m\rangle + 2\pi\sqrt{-1} \langle
m,z \rangle)}| < c(U)\exp{(-d|m|^2)}
$$
where $d$ is a positive
constant which depends on $\Omega$.  This implies that the series
(\ref{44}) converges uniformly on every compact set $U \subset
\C^n$ and defines an entire function.

{\bf Definition.} {\sl The theta function
$\theta(z,\Omega)$ is that defined by (\ref{44}).}

In the sequel we mean by $\theta(z)$ the function $\theta(z,\Omega)$.

It follows from (\ref{44}) that the theta function
satisfies the following periodicity conditions
$$
\theta(z+m,\Omega) = \theta(z,\Omega) ,
$$
\begin{equation}
\theta(z+\Omega m,\Omega) = \exp{(-\pi\sqrt{-1}\langle \Omega m,m
\rangle - 2\pi\sqrt{-1}\langle m,z\rangle)}\cdot\theta(z,\Omega),\ \
m \in \Z^n.
\label{45}
\end{equation}
This implies Proposition 3.1.1.

{\bf 3.1.1.}
{\sl $\theta(z,\Omega)$ is a global section of $L$ unique up to a
constant multiple.}

Since $M \rightarrow {\hat M}$ is an isomorphism,
every line bundle on $M$ with $c_1=[\omega]$
is uniquely expressed in the shape
$\tau^*_{\mu}(L)$ and the space of  its global sections is
generated by $\theta(z+\mu,\Omega)$.

{\bf 3.2. Theta functions with characteristics.}

Given a pair of vectors $a,b \in \R^n$,
$\theta[a,b](z,\Omega)$,
the {\it theta function with characteristics}
$a$ and $b$ is defined by
$$
\theta[a,b](z,\Omega)  \ ( = \theta \big[\begin{array}{c} a \\
b \end{array}\big](z, \Omega) )\ =
$$
\begin{equation}
\exp{(\pi\sqrt{-1}\langle \Omega a,a \rangle + 2\pi\sqrt{-1} \langle
a, z+b \rangle )} \cdot \theta(z+\Omega a + b,\Omega) =
\label{46}
\end{equation}
$$
\sum_{m \in \Z^n} \exp{(\pi\sqrt{-1}\langle \Omega
(m+a),(m+a) \rangle + 2\pi\sqrt{-1} \langle m+a, z+b \rangle )}.
$$
The analogues of the periodicity conditions
(\ref{45}) for theta functions with characteristics
are written as
$$
\theta[a,b](z+m,\Omega) = \exp{(2\pi\sqrt{-1}\langle a,m \rangle)}
\cdot \theta[a,b](z,\Omega),
$$
\begin{equation}
\theta[a,b](z+\Omega m,\Omega) =
\exp{(-\pi\sqrt{-1} \langle \Omega m,m \rangle -
2\pi\sqrt{-1}\langle m,z \rangle)} \times
\label{47}
\end{equation}
$$
\times \exp{( - 2\pi\sqrt{-1} \langle b,m \rangle)} \cdot
\theta[a,b](z,\Omega), \ \ m \in \Z^n.
$$

If the characteristics $a$ and $b$ are rational, then the
functions $\theta[a,b](z)$ determine sections of line bundles
on $M$.

Consider the tensor products of $L$.
Every bundle $L^d$ is defined by the multiplicators
$e^d_{\lambda_j}$ where $e_{\lambda_j}$ define $L$ and take the shape
(\ref{43}).

{\bf 3.2.1.}
{\sl Each of the families of functions

1) $\theta[a/d,0](d \cdot z, d \cdot \Omega), \ \
0 \leq a^i < d$ ,

2) $\theta[0,b/d](z, d^{-1} \cdot \Omega), \ \
0 \leq b^i < d$

\noindent
is a basis for global sections of $L^d$, i.e., a basis for
$H^0(M;\O(L^d))$.

These bases are related by the identity
\begin{equation}
\theta[0,b/d](z, d^{-1} \cdot \Omega) =
\sum_a \exp{(2\pi\sqrt{-1} \langle a,b \rangle)} \cdot
\theta[a/d,0](d \cdot z, d \cdot \Omega).
\label{48}
\end{equation}}

It follows from (\ref{47}) that these functions are sections of
$L^d$. For proving that these families are bases
it suffices to compare them with a basis $\{\theta_{{\bar l}}\}$
of the shape (\ref{42}) for sections of $L^d$.  Identity
(\ref{48}) is proved by comparing the Fourier series expansions of
the left-hand and right-hand sides of (\ref{48}).

{\bf 3.3. Modular transformations.}

In \S 1.4 the action of
$Sp(\Delta_{\delta},\Z)$ on the Siegel upper-half plane
${\cal H}_n$ is described. This action consists in changing
the matrix $\Omega$ for $\Omega'$
such that the Abelian varieties
$M=\C^n/\{\Delta_{\delta}\Z^n + \Omega \Z^n\}$ and
$\C^n/\{\Delta_{\delta}\Z^n + \Omega' \Z^n\}$ are
isomorphic. For the principal polarization
$\Delta_{\delta} = I_n$ we obtain the following action of
$Sp(2n,\Z)$:
$$
\Omega \rightarrow \Omega' =
(A\Omega+B)(C\Omega+D)^{-1},
$$
\begin{equation}
z \rightarrow z' = ((C\Omega+D)^{-1})^* \cdot z,
\label{49}
\end{equation}
$$
g =
\left(
\begin{array}{cc}
A & B \\ C & D
\end{array}
\right) \in Sp(2g,\Z).
$$

{\bf 3.3.1}
{\sl The theta function is
transformed by (\ref{49}) as follows
\begin{equation}
\theta[a',b'](z',\Omega') = C_0 \cdot
\exp{(\pi\sqrt{-1}\langle z, Tz \rangle)} \cdot
\theta[a,b](z,\Omega)
\label{50}
\end{equation}
with $C_0$ a constant independent of $z$ ,
$$
\left(
\begin{array}{c} a' \\ b' \end{array}
\right)
=
\left(
\begin{array}{cc}
D & -C \\ -B & A
\end{array}
\right)
\cdot
\left(
\begin{array}{c}
a \\ b
\end{array}
\right)
+
\frac{1}{2} \cdot
\left(
\begin{array}{c}
diag(C\cdot D^*) \\ diag(A\cdot B^*)
\end{array}
\right),
$$
$$
T = (C\Omega +D)^{-1} \cdot C.
$$}

The explicit formula for $C_0$ and
the proof of (\ref{50}) are given in \cite{Ig}.

{\bf 3.4. Addition theorems for theta functions.}

The most important theta function identities are
the Riemann addition theorems 3.4.1 and 3.4.2.

{\bf 3.4.1. Binary addition theorem of Riemann.}
{\sl
$$
\theta[a,c](z_1+z_2,\Omega) \cdot
\theta[b,d](z_1-z_2,\Omega) =
$$
\begin{equation}
\sum_{e \in \Z^n/2\Z^n}
\theta[\frac{a+b}{2}+\frac{e}{2}, c+d](2z_1,2\Omega)
\cdot \theta[\frac{a-b}{2}+\frac{e}{2}, c-d](2z_2,2\Omega)
\label{51}
\end{equation}
for $\theta : \C^n \rightarrow \C, z_1, z_2 \in \C^n$, and
$a,b,c,d \in \R^n$.}

It suffices to prove (\ref{51})
for $a=b=c=d=0$ to which the general case reduces by (\ref{46}).
In this case (\ref{51})
is written as
\begin{equation}
\theta(z_1+z_2,\Omega) \cdot
\theta(z_1-z_2,\Omega) = \sum_{e \in \Z^n/2\Z^n}
\theta[\frac{e}{2},0](2z_1,2\Omega) \cdot
\theta[\frac{e}{2},0](2z_2,2\Omega).
\label{52}
\end{equation}

The inverse of (\ref{52}) is the following addition theorem:
\begin{equation}
\theta(z_1+z_2,\Omega)
\cdot \theta(z_1-z_2,\Omega) = \sum_{e \in \Z^n/2\Z^n}
\theta[0,\frac{e}{2}](z_1,\frac{\Omega}{2}) \cdot
\theta[0,\frac{e}{2}](z_2,\frac{\Omega}{2}).
\label{53}
\end{equation}

The ternary addition theorem is a generalization of the
binary addition theorem of Riemann for the case of products
of four theta functions.

{\bf 3.4.2. Ternary addition theorem of Riemann.}
{\sl
$$
\theta[{\tilde a}_1,{\tilde b}_1]({\tilde z}_1) \cdot
\dots \cdot
\theta[{\tilde a}_4,{\tilde b}_4]({\tilde z}_4) =
$$
\begin{equation}
\frac{1}{2^n} \cdot
\sum_{c,d \in \frac{1}{2}\Z^n/2\Z^n}
\exp{(-4\pi\sqrt{-1}\langle d, {\tilde a_1} \rangle)}
\theta[a_1+c,b_1+d](z_1) \cdot \dots \cdot
\theta[a_4+c,b_4+d](z_4)
\label{54}
\end{equation}
for $\theta : \C^n \rightarrow \C, z_k \in \C^n$ and
$a_i,b_j \in \R^n$
where
${\tilde a}=aT, {\tilde b}=bT, {\tilde z}=zT$, and
$$
T =
\frac{1}{2}
\left(
\begin{array}{rrrr}
1 & 1 & 1 & 1 \\
1 & 1 & -1 & -1 \\
1 & -1 & 1 & -1 \\
1 & -1 & -1 & 1
\end{array}
\right).
$$}

Here the units in the formula for the matrix $T$ mean the
identity $(n \times n)$-matrices.

The proofs of the addition theorems are exposed, for
instance, in \cite{Mm2}. For $n=1$ the proofs are obtained by
straightforward comparing the Fourier series expansions
of the left-hand and right-hand sides of (\ref{52}--\ref{54}).
Generally, this requires repeating these arguments by
coordinates (\cite{Dubr1}).

{\bf 3.5. The theta divisor.}

A {\it divisor} on a variety $M$ is a formal integer linear
combination $D = \sum_k a_k V_k (a_k \in \Z)$ of finitely
many analytic hypersurfaces in $M$.
The set of divisors is the divisor group $Div(M)$
under the formal addition.  A divisor is called
{\it effective} if $a_k \geq 0$ for every $k$.

A {\it divisor $(f)$ of a function $f$} is a linear
combination of the sets of its zeros and poles
counted with multiplicities.
Therewith the sets of zeros and poles are taken with positive
and negative signs respectively.

The group $Div(M)$ is described in terms of sections of
sheaves as follows.
Let ${\cal M}^*$ be the multiplicative sheaf of germs of
meromorphic functions, on $M$, non-vanishing identically.
The sheaf $\O^*$ is a subsheaf of ${\cal
M}^*$.  Consider the space of sections of the sheaf ${\cal
M}^*/\O^*$.

{\bf 3.5.1.}
{\sl $Div(M) = H^0(M;{\cal M}^*/\O^*)$.}

Indeed, a global section of
${\cal M}^*/\O^*$ is defined by a family of local sections
$f_{\alpha} \in {\cal M}^*(U_{\alpha})$
such that $f_{\alpha}/f_{\beta} \in \O^*(U_{\alpha} \cap
U_{\beta})$.  Thus local divisors
$(f_{\alpha})$ are glued into a global
divisor which is identified with the section
$\{f_{\alpha}\}$ of ${\cal M}^*/\O^*$.

Consider the exact sequence of sheaves
\begin{equation}
0 \rightarrow \O^*
\rightarrow {\cal M}^* \rightarrow {\cal M}^*/\O^* \rightarrow 0.
\label{55}
\end{equation}
Assign for this sequence the exact cohomology sequence
which we need only extracting the following fragment
\begin{equation}
H^0(M;{\cal M}^*) \rightarrow H^0(M;{\cal
M}^*/\O^*) \stackrel{\delta}{\rightarrow} H^1(M;\O^*) = Pic(M).
\label{56}
\end{equation}
The homomorphism
\begin{equation}
\delta : Div(M) \rightarrow Pic(M) \label{57}
\end{equation}
assigns to a divisor $D$ an associated line bundle
as follows.
Let $D$ is defined by a section $\{f_{\alpha}\}$ of
${\cal M}^*/\O^*$. Then the family of functions
$g_{\alpha\beta} = f_{\alpha}/f_{\beta}$ defines the
coordinate transformations for the bundle $\delta(D) = [D]$.

Two divisors $D_1$ and $D_2$ are called {\it linearly equivalent},
$D_1 \approx D_2$, if $D_1 - D_2 \in
Ker\,\delta$. By (\ref{56}), $D_1 \approx D_2$ implies that
$D_1 - D_2$ is the divisor of zeros of a meromorphic
function on $M$.

The {\it theta divisor} $\Theta = (\theta)$ is the zero set of
the theta function (\ref{44}).

The following proposition is clear.

{\bf 3.5.2.}
{\sl 1) The line bundle $L = [\Theta]$ is defined by
multiplicators (\ref{43}) and the theta function generates
$H^0(M;\O(L))$.  2) $[\Theta + \mu] =
\tau_{-\mu}^*(L)$.}

To make our exposition complete, we state the following fact
which is proved by applying the Stokes theorem. The proof is
exposed, for instance, in \cite{GH}.

{\bf 3.5.3.}
{\sl Let $M$ be a compact complex
$n$-dimensional manifold and let $L$ be a line bundle on
$M$ such that $L=[D]$. Then $c_1(L)$ is dual
to the cycle $D \in H_{2n-2}(M;\Z)$
by the Poincare duality.}

\vskip2.5mm

\begin{center}
{\bf \S 4. Theta functions and mappings of tori into
projective spaces. Secants of Abelian varieties}
\end{center}

\vskip2.5mm

{\bf 4.1. Linear systems.}

Let $L = [D]$ be a line bundle, on a compact complex manifold,
associated with an effective divisor $D$. The divisor
$D$ is defined locally by the functions $f_{\alpha} \in
{\cal M}(U_{\alpha})$ where ${\cal M}$ is the sheaf of germs of
meromorphic sections of $L$. It is seen that these functions define
a meromorphic section $f_0$ of $L$.

Denote by $\L(D)$ the space of all meromorphic functions
$f$ on $M$ such that the divisor $(f) + D$
is effective, i.e., $(f)+D \geq 0$, and denote by $|D|$
the set of effective divisors linearly equivalent
to $D$. There exists a one-to-one correspondence
between the points of the projective space associated with
$\L(D)$ and the divisors of $|D|$. It is as follows:
\begin{equation}
D' \in |D| \longleftrightarrow [f] \in P\L(D) : (f) = D' - D.
\label{58}
\end{equation}
A section $f \cdot f_0$ of ${\cal M}$
is holomorphic if and only if  $f \in
\L(D)$.  Therefore the correspondence (\ref{58}) is extended
as follows
\begin{equation}
|D| \longleftrightarrow PH^0(M;\O([D])).
\label{59}
\end{equation}

A subset of effective divisors corresponding to
a linear subspace of $H^0(M;\O([D]))$
by (\ref{59}) is called a {\it linear system of
divisors}.  The simplest example is the
{\it complete linear system} $|D|$.

The set of {\it base points} $B(X)$ of a linear system $X$
is the intersection of all divisors of this system.
Assume that $E(X) \subset H^0(M;\O([D]))$ is a subspace
corresponding to $X$ by (\ref{59}). Then $B(X)$ comprises
the points of $M$ at which all sections in $E(X)$
vanish.

Given a linear system $X$,
the following mapping is defined up to projective
transformations of $\C P^k$:
\begin{equation}
\varphi : M
\setminus B(X) \rightarrow \C P^k \ \ : \ \ \varphi(x) = (f_1(x) :
\dots :  f_{k+1}(x))
\label{60}
\end{equation}
where $\{f_j\}$ is a basis for $E(X)$.

{\bf 4.2. The linear systems $|d\Theta|$. The Lefschetz theorem.}

Consider the mappings (\ref{60}) for the linear systems
$|d\Theta|$ where $d \geq 2$.

{\bf 4.2.1.}
{\sl The sets of base points of
$|d\Theta|$ are empty for $d \geq 2$.}

To prove this proposition it suffices to
find for every point $z \in M = \C^n/\{\Z^n +
\Omega \Z^n\}$ sections
$f_d \in H^0(M;\O(L^d))$ such that $f_d(z) \neq 0$.

By (\ref{45}), the function
$f^{\mu}_d(z) = \theta(z+\mu_1)\cdot \dots\cdot
\theta(z+\mu_d)$ defines a global section of $\O(L^d)$ for $\mu_1
+\dots+\mu_d=0$.  Since the theta divisor is an analytic
hypersurface in $M$, for every point $z \in M$
there exists a vector $\nu=(\nu_1,\dots, \nu_d)$
such that $\theta(z+\nu_k) \neq 0$ for any $k$.  Now it
suffices to take $f^{\nu}_d$ for $f_d$. The proof of
Proposition 4.2.1 is complete.

Denote by
\begin{equation}
\varphi_d : M = \C^n/\{\Z^n + \Omega \Z^n\}
\rightarrow \C
P^{d^n-1}
\label{61}
\end{equation}
the mapping (\ref{60}) for
the linear system $|d\Theta|$.
By Proposition 3.2.1, $\varphi_d$ can be defined by the
following formula
\begin{equation}
\varphi_d(z) =
(\theta[\frac{a_1}{d},0](d\cdot z,d\cdot \Omega) : \dots :
\theta[\frac{a_{d^n}}{d},0](d\cdot z,d\cdot \Omega)) \in \C P^{d^n-1}
\label{62}
\end{equation}
where $\{a_i\} = \Z^n/d\Z^n$.

{\bf 4.2.2. Lefschetz Theorem.}
{\sl The mapping $\varphi_d$ is an embedding for $d \geq 3$.}

For brevity we prove this only for $d=3$.
It is easily seen how from the proof of Proposition 4.2.1
to obtain the proof in the general case.

1) First, prove that the rank of $\varphi_d$ is maximal.

Take an arbitrary point $z_* \in M$.
Choose a basis $(\theta_0,\theta_1,\dots,\theta_N)$ for
$H^0(M;\O(L^2))$ such that $\theta_0(z_*)=1$ and
$\theta_1(z_*)=\dots=\theta_N(z_*)=0$.  Then the rank of
$\varphi_2$ at $z_*$ coincides with the
rank of the following matrix
$$ J(z_*) =
\left(
\begin{array}{ccc}
\frac{\partial \theta_1(z_*)}{\partial z^1}
& \dots & \frac{\partial \theta_1(z_*)}{\partial z^n} \\
\dots & \dots & \dots \\
\frac{\partial \theta_N(z_*)}{\partial z^1} & \dots
& \frac{\partial \theta_N(z_*)}{\partial z^n} \end{array} \right).
$$

Assume that $rank J(z_*) < n$. This means that there exist
non-trivial vectors $a_1,\dots,a_n$ such that
$$
a_1\frac{\partial \theta_k(z_*)}{\partial z^1} + \dots +
a_n\frac{\partial \theta_k(z_*)}{\partial z^n} = 0
$$
for any $k$.
Choose $\mu$ and $\nu$ such that
\begin{equation}
\theta(z_*+\nu) \cdot \theta(z_*-\mu-\nu) \neq 0.
\label{63}
\end{equation}
By (\ref{45}), the function
\begin{equation}
F(z,\mu,\nu) = \theta(z+\mu)\cdot \theta(z+\nu)\cdot
\theta(z-\mu-\nu)
\label{64}
\end{equation}
defines a section of
$L^d\,(d=3)$ for any $\mu,\nu$.
Hence this function is represented by a linear combination of
functions of the basis $\{\theta_k\}$ and the identity
\begin{equation}
a_1\frac{\partial
F(z_*,\mu,\nu)}{\partial z^1} + \dots + a_n\frac{\partial
F(z_*,\mu,\nu)}{\partial z^n} = 0
\label{65}
\end{equation}
holds. Define the function
$$
\xi(z) = a_1\frac{\partial
\log\theta(z)}{\partial z^1} + \dots + a_n\frac{\partial
\log\theta(z)}{\partial z^n}.
$$
It is seen that (\ref{65}) is equivalent to the following
equality
\begin{equation}
\xi(z+\mu) + \xi(z+\nu) + \xi(z-\mu-\nu) = 0.
\label{66}
\end{equation}
By (\ref{63}) and (\ref{66}), the function
${\tilde \xi}(z)= \xi(z+\mu)$ has no poles and
thus it is an entire function.
It follows from (\ref{25}) that ${\tilde \xi}(z)$
satisfies the following periodicity conditions
\begin{equation}
{\tilde \xi}(z+\lambda_k) =
{\tilde \xi}(z), \ \ {\tilde \xi}(z+\lambda_{n+k}) = {\tilde \xi}(z)
-2\pi\sqrt{-1}a_k.
\label{67}
\end{equation}
By (\ref{67}),
the derivatives $\partial {\tilde \xi}/\partial z^k$
of ${\tilde \xi}$ are $\Lambda$-periodic
and thus they are bounded.
Since these derivatives are entire functions,
they are constant and ${\tilde \xi}(z) = \sum_k b_k z^k
+ c$.  However, it follows from (\ref{67}) that
$b_k \equiv 0$ and, therefore,
${\tilde \xi}$ is constant. Since ${\tilde
\xi}(z+\lambda_{n+k}) - {\tilde \xi}(z) = -2\pi\sqrt{-1}a_k$,
we have $a_k \equiv 0$.

We conclude that the family of constants $\{a_k\}$ is
trivial and $J(z_*)$ has  maximal rank.
Since $z_*$ is  arbitrary, the rank of
$\varphi_d$ equals $n$ everywhere.

2) Prove that $\varphi_d(z_1) = \varphi_d(z_2)$
implies $z_1 - z_2 \in \Lambda$. This amounts to
the injectivity of $\varphi_d$.

Assume that $\varphi_d(z_1) = \varphi_d(z_2)$.
Hence there exists a non-zero constant $C$ such that
$\psi(z_1) = C\psi(z_2)$ for every section $\psi$
of $L^d$.  Consider a function
$F(z,\mu,\nu)$ of the shape (\ref{64}). Since this function defines a
section of $L^d$, the following identity
\begin{equation}
\frac{F(z_1,\mu,\nu)}{F(z_2,\mu,\nu)} = \frac{\theta(z_1+\mu)\cdot
\theta(z_1+\nu)\cdot \theta(z_1-\mu-\nu)} {\theta(z_2+\mu)\cdot
\theta(z_2+\nu)\cdot \theta(z_2-\mu-\nu)} = C
\label{68}
\end{equation}
holds.
For arbitrary $\mu$ it is possible to find $\nu$ such that
\begin{equation}
\theta(z_1+\nu)\cdot \theta(z_1-\mu-\nu) \cdot
\theta(z_2+\nu)\cdot \theta(z_2-\mu-\nu) \neq 0.
\label{69}
\end{equation}
By (\ref{68}--\ref{69}), the function
${\hat \xi} = \log{(\theta(z_1+z)/\theta(z_2+z))}$
is holomorphic.
As for ${\tilde \xi}$ it is shown that the
derivatives ${\hat \xi}(z)$ are constant. Therefore,
this function is linear:
${\hat \xi}(z) = 2\pi\sqrt{-1}\sum b_k z^k + c$.
It follows from (\ref{45}) that
$$
{\hat \xi}(z+\lambda_k) = {\hat \xi}(z) +
2\pi\sqrt{-1}\beta_k,
$$
$$
{\hat \xi}(z+\lambda_{n+k}) = {\hat
\xi}(z) + 2\pi\sqrt{-1}(z_2^k-z_1^k) + 2\pi\sqrt{-1}\alpha_k, \ \
\alpha_i, \beta_j \in \Z.
$$
Since ${\hat \xi}$ is linear,
we have $b_k = \beta_k$ and
$$
z_2^k-z_1^k = -\alpha_k +
\sum_l\Omega_{kl}\beta_l.
$$
The last equality is equivalent to $z_2 - z_1 \in \Lambda$.

The proof of the Lefschetz theorem is complete.
The following theorem is a corollary of the Lefschetz and
Chow theorems.

{\bf 4.2.3.}
{\sl A principally-polarised Hodge torus is an Abelian variety.}

In the case of a non-principally-polarised Hodge torus
$(M,[\omega])$ the space of sections of a line bundle
$L \rightarrow M$ with $c_1(L)=[\omega]$ is spanned by
the functions $\theta_{\bar l}$ (\ref{40}--\ref{42}).

Substituting any of the functions $\theta_{\bar l}$
instead of the theta function $\theta(z,\Omega)$ in the proofs of
Proposition 4.2.1 and the Lefschetz theorem  we may repeat all
arguments and prove the following theorem.

{\bf 4.2.4. Lefschetz Theorem (the general case).}
{\sl Let $(M,[\omega])$ be a Hodge torus, let $L$ be a line bundle
with $c_1(L)=[\omega]$, and let  $X_d$ be the linear system
corresponding to $H^0(M;\O(L^d))$. Then the set of base points
of $X_d$ is empty for $d \geq 2$ and this system defines by
(\ref{60}) an embedding of $M$ into a complex projective
space for $d \geq 3$.}

This theorem and the theorem of Chow imply the following criterion.

{\bf 4.2.5. Criterion for a complex torus to be algebraic.}
{\sl A complex torus is an Abelian variety if and only if
it is a Hodge torus.}

In the general case the Kodaira theorem is proved by
construction of a line bundle $L$ with $c_1(L)=[\omega]$ from a Hodge
form $\omega$ and embedding $M$ into  $\C P^N$ by a mapping of
the shape (\ref{60}) for a line system corresponding to $H^0(M;\O(L^d))$
for sufficiently large $d$ (\cite{C,GH}).

{\bf 4.3. The mapping $\varphi_2$ and Kummer varieties.}

Let  $M = \C^n/\Lambda$ be an Abelian variety with the principal
polarization $[\omega]$.  The reflection acts
on $M$ as follows
$$
\sigma : M \rightarrow M \ \ :
\sigma (z) = -z.
$$
This action explains why the claim of the
Lefschetz theorem fails for $d=2$.  The mapping $\varphi_2$
is defined by the section $\theta[e,0](d\cdot z, 2\Omega)$ of the
bundle $[2\Theta]$ where $2e \in \Z^n/2\Z^n$.  All these functions
are even and the rank of $\varphi_2$ vanishes at zero.
This also means that $\varphi_2$ splits as follows
$$
M \stackrel{\pi}{\longrightarrow} M/\sigma \stackrel{\Phi}
{\longrightarrow} \C P^{2^n-1}, \ \ \varphi_2 = \Phi \circ \pi
$$
with $\pi$ the natural projection.

The variety $K(M,\Theta) = \Phi(M/\sigma)$ is called
a {\it Kummer variety} and the mapping $\Phi$ is called
a {\it Kummer map}.

A polarised Abelian variety $(M,[\omega])$ is called
{\it reducible} if it is a direct product of two Abelian varieties
$(M_1,[\omega_1])$ and $(M_2,[\omega_2])$ of positive
dimensions and $[\omega] = [\omega_1] +
[\omega_2]$.

{\bf 4.3.1. Theorem on injectivity of the Kummer map.}
{\sl Let the principally-polarised Abelian variety
$M = \C^g/\{\Z^g + \Omega\Z^g\}$ be irreducible. Then
the Kummer map
\begin{equation}
\Phi :  z \rightarrow (\theta[n_1,0](2z, 2\Omega) :
\dots :  \theta[n_r,0](2z, 2\Omega)),
\label{70}
\end{equation}
where $\{2n_j\} = \Z^g/2\Z^g$ and $r=2^g$, is an embedding of
the variety with singularities $M/\sigma$ into $\C P^{2^n-1}$.}

We explain the proof of this theorem assuming Proposition
4.3.2 to be valid.

{\bf 4.3.2.}
{\sl Let an Abelian variety be principally-polarised and
irreducible. Then its theta divisor $\Theta$ is irreducible
and non-invariant under translations, i.e.,  $\Theta + \mu = \Theta$
implies $\mu \in \Lambda$.}

For brevity we omit the proof of this proposition. Notice that its
proof is based on the equivalence of reducibility of a
principally-polarised Abelian variety to  reducibility of its
theta divisor and the non-compatibility of existence of
non-trivial translations transforming variety into itself with a
principal polarization.

Assume that $\varphi_2(z_1) = \varphi_2(z_2)$. Then there exists
a non-zero constant $C$ such that
$$
\frac{\theta[e,0](2z_1,2\Omega)}{\theta[e,0](2z_2,2\Omega)} = C
$$
for any $e \in \frac{1}{2}\Z^n/2\Z^n$. It follows from
(\ref{52}) that
\begin{equation}
\frac {\theta(z+z_1)\cdot \theta(z-z_1)} {\theta(z+z_2)\cdot
\theta(z-z_2)} = C
\label{71}
\end{equation}
for any $z \in \C^n$.

Let $(z+z_1)$ be a non-singular point of the theta divisor
and let $\theta(z-z_1)$ does not vanish. Then, by (\ref{71}),
${\theta(z+z_2)\cdot \theta(z-z_2)} = 0$.

Suppose that
$\theta(z - z_2) = 0$. Then in a neighbourhood of $z$ the
function $\theta(z+z_1)/\theta(z-z_2)$ is holomorphic and,
since the theta divisor is irreducible,
$\Theta + (z_1+z_2) = \Theta$.
Since the theta divisor is non-invariant under translations,
$z_1 + z_2 \in \Lambda$. The case of $\theta(z + z_2) = 0$
is considered similarly. This completes the proof.

An important example of a Kummer variety is the quotient space of a
two-dimensional torus $\C^2/\Lambda$ for the reflection.
It contains 16 singular points corresponding to the sublattice
$\Lambda/2$. By resolution of singularities we obtain a
non-singular algebraic surface known as a {\it Kummer surface}, or
a {\it K3 surface}.

{\bf 4.4. Secants of Abelian varieties.}

Let $M = \C^n/\Lambda$ be an irreducible principally-polarised
Abelian variety and $K(M,\Theta) \subset \C P^{2^n-1}$ be its
Kummer variety.

{\bf Definition.}
{\sl N-secant of an Abelian variety $M$
is a $(N-2)$-dimensional plane $V \,(= \C P^{N-2})
\subset \C P^{2^n-1}$ which intersects the Kummer variety
$K(M,\Theta)$ at $N$ distinct points.}

Let $V$ be a $N$-secant passing through points
$\Phi(z_1),\dots,\Phi(z_N)$. Then there exist non-zero constants
$c_1,\dots,c_N$ such that
\begin{equation}
c_1\cdot \theta[e,0](2z_1,2\Omega) + \dots +
c_N\cdot \theta[e,0](2z_N,2\Omega) = 0
\label{72}
\end{equation}
for any $e \in \frac{1}{2}\Z^n/2\Z^n$.  The addition theorem
(\ref{52}) and (\ref{72}) imply Proposition 4.4.1.

{\bf 4.4.1.}
{\sl The following statements (i) and  (ii) are
equivalent:

(i) points $\Phi(z_1),\dots,\Phi(z_N)$  of
$K(M,\Theta)$ lie on the same $N$-plane and satisfy
the relation $c_1 \Phi(z_1) + \dots + c_N \Phi(z_N)=0$ ,

(ii) the theta functional identity
\begin{equation}
c_1\theta(z+z_1)\theta(z-z_1) + \dots +
c_N\theta(z+z_N)\theta(z-z_N) = 0, \ \ z \in \C^n
\label{73}
\end{equation}
holds.}

Consider conditions on the theta divisor
providing the existence of $N$-secants. It seems that they
were first considered in \cite{Weil}.
Denote by $\Theta_{\mu}$ the translate of the theta
divisor by $\mu$, i.e., $\Theta+\mu$.

{\bf 4.4.2.}
{\sl Assuming $\dim M \geq (N-1)$, the statements (i) and
(ii) are equivalent:

(i) there exist pairwise distinct points
$\mu_1,\dots,\mu_{N-1},x,y \in M$ such that the analytic
subset $\Theta_{\mu_1} \cap \dots \cap
\Theta_{\mu_{N-1}}$ is of codimension $N-1$ and
\begin{equation}
\Theta_{\mu_1} \cap \dots \cap \Theta_{\mu_{N-1}}
\subset \Theta_x \cup \Theta_y;
\label{74}
\end{equation}

(ii) the Kummer variety $K(M,\Theta)$ admits an $N$-secant and
(\ref{73}) holds with $c_N \neq 0$ for
$z_1 = (x+y)/2-\mu_1, \dots, z_{N-1}=(x+y)/2-\mu_{N-1}$, and
$z_N = (x-y)/2$.}

By Proposition 4.4.1, (ii) implies (i).
While proving that (ii) follows
from (i) we apply the Koszul complex of the intersection
$\Theta_{\mu_1} \cap \dots \cap
\Theta_{\mu_{N-1}}$ (\cite{Har}).  We consider only the
case of $N=4$ (\cite{BD1}) which illustrates how the
proof becomes more complicated for  $N>4$ and how it simplifies for
$N=3$ (\cite{Mm2}).

Consider the bundle $L = [2\Theta+x+y]$ and its sections
$s_0 = \theta(z-x)\theta(z-y),
s_1=\theta(z-\mu_1)\theta(z-(x+y-\mu_1)), \dots,
s_3=\theta(z-\mu_3)\theta(z-(x+y-\mu_3))$. Denote by $I$
the subsheaf $\O_M$ generated by the sections vanishing on
$\Theta_{\mu_1} \cap \dots \cap
\Theta_{\mu_{N-3}}$.

Consider the tensor product of the Koszul complex $(K_i,d_i)$
and  $\O(L)$:
$$
0 \rightarrow \O(L^{-2})
\stackrel{d_1}{\longrightarrow} \O(L^{-1}) \oplus \O(L^{-1}) \oplus
\O(L^{-1}) \stackrel{d_2}{\longrightarrow} $$ $$ \O \oplus \O \oplus
\O \stackrel{d_3}{\longrightarrow} I \otimes \O(L) \rightarrow 0
$$
where  $d_1(f)=(s_1f,s_2f,s_3f), d_2(f,g,h) = (c_2h-c_3g,
c_3f-c_1h, c_1g - c_2f)$,  and $d_3(f,g,h)=(s_1f+s_2g+s_3h)$.

Since $\dim \Theta_{\mu_1} \cap
\Theta_{\mu_2} \cap \Theta_{\mu_3} = \dim M - 3$,
this sequence is exact.

By the Kodaira vanishing theorem applied to the positive
bundle $L$ (\cite{C,GH}), $H^k(M;\O(L^{-m}) = 0$
for  $m > 0$ and  $k <\dim M$.  Successively applying
the exact cohomology sequence to the triples $0 \rightarrow
A= K_i/d_{i-1}(K_{i-1}) \rightarrow K_{i+1}\rightarrow
K_{i+1}/d_{i+1}(A)\rightarrow 0$ we arrive at the exact sequence
\begin{equation}
0\rightarrow
{\tilde A} \stackrel{d_{N-2}}{\longrightarrow} \O \oplus \dots \oplus
\O \stackrel{d_{N-1}}{\longrightarrow} I\cdot \O(L)
\label{75}
\end{equation}
where
$H^1(M;{\tilde A}) = 0$. It follows from the exact cohomology
sequence for (\ref{75}) that
$H^0(M;I\otimes \O(L)) = d_{N-1}(H^0(M;\O \oplus
\dots\oplus \O)) = d_{N-1}(\C^{N-1})$.  By (\ref{74}), the function
$s_0$ defines a section of $I\otimes \O(L)$ and, therefore, is a linear
combination of the sections $s_1,\dots,s_{N-1}$.  By (\ref{73}), this
is equivalent to existence of a secant. This proves Proposition
4.4.2.

\vskip2.5mm

\begin{center}
{\bf Chapter 2. Theta functions of Riemann surfaces}
\end{center}

\vskip2.5mm

\begin{center}
{\bf \S 5. Theta functions of Jacobi varieties}
\end{center}

\vskip2.5mm

{\bf 5.1. Riemann surfaces and their period matrices.}

A {\it Riemann surface} is a complex manifold of complex
dimension one. A compact surface $\Gamma$ is homeomorphic to a sphere
with handles. The number of handles is called the {\it genus} of
$\Gamma$.

A basis $a_1,\dots,a_g,b_1,\dots,b_g$ for $1$-cycles on a compact
surface $\Gamma$ of genus $g$ is called {\it canonical} if its
intersection form is as follows
$$
a_i \circ
a_j = b_i \circ b_j = 0, \ \ a_i \circ b_j = \delta_{ij}
$$
where
$\delta_{ij}$ is the Kronecker symbol.
The canonical basis is realised by contours such that by
cutting $\Gamma$ along these contours we obtain
a $4g$-gon ${\tilde \Gamma}$ with sides
$a_1, b_1, a_1^{-1}, b_1^{-1}, \dots, a_g,
b_g, a_g^{-1}, b_g^{-1}$ taken in order of traversal.

In the sequel we suppose that Riemann surfaces are
compact and the bases of $1$-cycles on surfaces are canonical.

Since a Hermitian metric on a one-dimensional complex manifold is
proportional to a Hodge metric, every Riemann surface is projective
algebraic, or an {\it algebraic curve}.

A {\it differential} is a $1$-form $\omega$ on
$\Gamma$.  A differential is called {\it  Abelian} (or holomorphic)
if it is represented in a neighbourhood of every point
by the form $\omega=f(z)dz$ where $f(z)$ is a holomorphic function.
In the same manner meromorphic differentials are defined.
It is clear that every Abelian differential
as well as the adjoint differential
${\bar \omega} = {\bar f}(z) d{\bar z}$
is closed.

{\bf 5.1.1. Reciprocity Law (for smooth differentials).}
{\sl Let $\omega$ and $\omega'$ be closed
differentials on a Riemann surface of genus $g$.
Then the identity
\begin{equation}
\int_{\Gamma} \omega \wedge \omega' = \sum^g_{k=1}
(\int_{a_k} \omega \cdot \int_{b_k} \omega' -
\int_{b_k} \omega \cdot \int_{a_k} \omega' )
\label{76}
\end{equation}
holds.}

In order to prove this proposition we apply
the Stokes theorem to a domain
${\tilde \Gamma}$ which is projected onto $\Gamma$
by the universal covering $\{z \in C : Im z > 0\} \rightarrow \Gamma$.
By applying this proposition to an Abelian differential $\omega$
and its adjoint differential ${\bar \omega}$ we infer
the following proposition.

{\bf 5.1.2.}
{\sl Let $\omega$ be a non-zero holomorphic differential
on a Riemann surface of genus $g$. Then the identity
\begin{equation}
Im \sum^g_{k=1} \int_{a_k} \omega \cdot \int_{b_k}
\omega < 0
\label{77}
\end{equation}
holds.}

The conjugation $\omega\rightarrow {\bar \omega}$
induces an involution on the group
$H^1(M;\C)$ decomposed into the direct sum of the subgroups
$H^{1,0}(M;\C)$ and $H^{0,1}(M;\C)$ generated by holomorphic and
anti-holomorphic differentials.
These subgroups permute under the conjugation and each of them is
isomorphic to $\C^g$. Considering (\ref{77}) we obtain

{\bf 5.1.3.}
{\sl The space of holomorphic differentials on a Riemann
surface of genus $g$ is isomorphic to $\C^g$ and every holomorphic
differential is determined uniquely by its $a$-periods.}

Hence to every canonical basis for $1$-cycles
there is assigned the {\it canonical basis for Abelian differentials}
$\omega_1,\dots,\omega_g$ which is defined uniquely by the formula
\begin{equation}
\int_{a_j} \omega_k
= \delta_{jk}.
\label{78}
\end{equation}
The {\it period matrix}, of Riemann surface $\Gamma$,
associated with this basis is constructed as follows
\begin{equation}
B(\Gamma)_{jk}  =
\int_{b_j} \omega_k.
\label{79}
\end{equation}
As one can see the period matrix
is determined from a canonical basis for $H_1(\Gamma)$.
Passage to a new canonical basis is achieved by
a transformation in $Sp(2g,\Z)$. Therewith the period matrix
is transformed by (\ref{16}).  By applying (\ref{77}) to the
differentials $\omega = \sum_k c_k \omega_k$, where $c_k \in \R$,
we obtain the following proposition.

{\bf 5.1.4.}
{\sl The period matrix of a Riemann surface
lies in the Siegel upper half-plane,
$B(\Gamma) \in {\cal H}_g$.}

Riemann surfaces of genus $g$ are
in a one-one correspondence with points of the
complex manifold ${\cal M}_g$, the moduli space of
Riemann surfaces of genus $g$.  It is known that
\begin{equation}
\dim_{\C}
{\cal M}_g = \cases{ 0 , & $g=0$; \cr 1, & $g=1$; \cr 3g-3, & $g
\geq 2$.}
\label{80}
\end{equation}
Hence the mapping
\begin{equation}
B :  {\cal M}_g \rightarrow \A_g = {\cal H}_g /
Sp(2g,\Z)
\label{81}
\end{equation}
is defined.
For a suitable compactification of ${\cal M}_g$, this mapping $B$
is algebraic. The same statement is valid for
the Prym mappings (\ref{99}) and (\ref{116}) which are defined
below (\cite{MFK}).

{\bf 5.2. Jacobi variety and Abel map.
Theta function of a Riemann surface.}

It follows from (\ref{81}) that to
every Riemann surface $\Gamma$
we may canonically assign
a principally-polarised Abelian variety,
its {\it Jacobi variety}:
\begin{equation}
J(\Gamma) = \C^g / \{\Z^g + B(\Gamma)\Z^g\}.
\label{82}
\end{equation}

Given a point $p_0 \in \Gamma$,
we consider the mapping
\begin{equation}
A : \Gamma \rightarrow
J(\Gamma) :  \ \ A(p) = (\int^p_{p_0} \omega_1,
\dots,\int^p_{p_0}\omega_g)
\label{83}
\end{equation}
which is called an {\it Abel mapping}. Its definition does not
depend on the choice of the path of integration $\gamma$.
Indeed, let us change $\gamma$ for $\gamma'$.
The contour $\gamma^{-1}\cdot\gamma'$ is closed and realises the
$1$-cycle $\zeta$.  Therewith $A(p)$ is changed for the sum of
$A(p)$ and the vector
$(\int_{\zeta}\omega_1,\dots,\int_{\zeta}\omega_g)$ which belongs to
the lattice $\{\Z^g+B(\Gamma)\Z^g\}$.

An Abel mapping extends linearly onto the divisor group of
$\Gamma$.

Let $\theta(z) = \theta(z,B(\Gamma))$ be the {\it theta
function} of a Riemann surface
$\Gamma$.  Construct the functions
$f(p) = \theta(A(p)-e)$.  Since their values depend on
paths of integrations in (\ref{83}),  these functions are
many-valued on $\Gamma$.
The differences between branches are given by (\ref{45}).
However, by (\ref{45}),
if $f$ does not vanish identically then
the set of zeros of $f(p)$ is a correctly defined divisor on $\Gamma$.

{\bf 5.2.1.}
{\sl There is a unique point $K_{\Gamma} \in
J(\Gamma)$ called the {\it ``vector of Riemann constants''}
such that one of the two possibilities takes place:

1) $A(\Gamma) \subset \Theta_e = \Theta+e$  and
$\theta(A(p)-e) \equiv 0$;

2) if $f(p)=\theta(A(p)-e)=0$ then $p \in \{p_1,\dots,p_g\}$ and
\begin{equation}
A(p_1)+\dots+A(p_g) + K_{\Gamma} = e.
\label{84}
\end{equation}
In the second case the divisor $D=p_1+\dots+p_g$
is uniquely determined by (\ref{84}).}

Since the second case holds for generic
$e \in J(\Gamma)$, (\ref{84}) determines almost everywhere
the inverse of
$$
A : S^g\Gamma\rightarrow J(\Gamma) , \ \
A(p_1,\dots,p_g)=A(p_1)+\dots+A(p_g)
$$
where $S^g\Gamma$ is the
$g$-th symmetric power of $\Gamma$.

It is clear that
if an effective divisor  $D$ of degree $g$
is the divisor of poles of non-trivial meromorphic function then
$A$ is non-invertible at $A(D)$.
Such divisors are called {\it special}.

{\bf 5.2.2.}
{\sl Let $\Sigma \subset S^g$ be the subset formed by special
divisors. Then the mapping $A$ is inverted by (\ref{84})
on $S^g \setminus \Sigma$.}

The theta divisor $\Theta$ is described by the following
theorem.

{\bf 5.2.3. Riemann Theorem on the theta divisor of a Jacobi
variety.}
\begin{equation}
\Theta = A(S^{g-1}) + K_{\Gamma}.
\label{85}
\end{equation}

{\bf 5.3. Abel and Riemann--Roch theorems.}

A divisor $D$ on a surface is a linear combination of finite
number of points of the shape $D = a_1 p_1 + \dots + a_n
p_n, a_k \in \Z, p_k \in \Gamma$.  The {\it degree} of $D$
is the following sum of coefficients
$\deg D = a_1 \dots + a_n \in \Z$.

Recall that
$[D] = \delta(D)$ is a line bundle constructed from $D$
(see \S 3.5, (\ref{57}). The proof of the
following theorem is given in \cite{GR}.

{\bf 5.3.1.}
{\sl Let $M$ be a projective algebraic variety.
Then every line bundle on $M$ has a meromorphic section and
hence the homomorphism (\ref{57})
$\delta :  Div(M) \rightarrow Pic(M)$ is a mapping onto.}

The Abel theorem enables us to linearise $\delta$.

{\bf 5.3.2. Abel Theorem.}
{\sl Assume that  $\deg D_1 = \deg D_2$.
Then $D_1-D_2$ is the divisor of a meromorphic
function on $\Gamma$ if and only if
$A(D_1) = A(D_2)$.}

By the exact sequence
(\ref{56}), Proposition 5.3.1, and the Abel theorem,
we have the following proposition.

{\bf 5.3.3. Structure of the group $Pic(\Gamma)$.}
{\sl The sequence
\begin{equation}
0 \rightarrow J(\Gamma) \rightarrow Pic(\Gamma)
\stackrel{\deg}{\rightarrow} \Z \rightarrow 0,
\label{86}
\end{equation}
is exact.}

The first homomorphism in (\ref{86}) is an embedding of a subgroup
generated by line bundles $[D]$ with $\deg D = 0$.
This sequence is an analogue of (\ref{35}).  Indeed,
$H^2(\Gamma;\Z)=\Z$ and it is easy to infer from
(\ref{19}--\ref{20}) that the homomorphism $\deg$
assigns to a line bundle its first Chern class.

The {\it divisor of a meromorphic differential} $\omega$ is the
formal sum $(\omega)$ of its zeros and poles taken with positive
and negative signs counted with multiplicities.
By the Abel theorem, all the divisors of meromorphic
differentials are linearly equivalent. This equivalence class is
called the {\it canonical class of a surface} and is denoted by
$C(\Gamma)$.  It is described by the following proposition.

{\bf 5.3.4.}
{\sl   1) $\deg C(\Gamma) = 2g-2$;
2) $A(C(\Gamma)) = -2K_{\Gamma}$.}

The dimensions $l(D)$ of the spaces
$H^0(\Gamma;\O([D])) = {\cal L}(D)$ (see \S 4.1) are given
by the Riemann--Roch theorem.

{\bf 5.3.5. Riemann--Roch Theorem.}
\begin{equation}
l(D) = \deg D - g + 1 + l(C(\Gamma) - D).
\label{87}
\end{equation}

The Riemann--Roch formula
for non-special effective divisors of degree $g$
degenerates into  $l(D)=1$.

{\bf 5.4. Trisecants of Jacobi varieties.}

{\bf 5.4.1. Fay Trisecant Formula (\cite{Fay1}).}
{\sl For every quadruple $p_1,p_2$ ,$p_3$, and $p_4$
of points of a Riemann surface $\Gamma$
there exist constants $c_1({\bar p}), c_2({\bar
p})$, and $c_3({\bar p})$ such that
$$
c_1 \cdot \theta(z + p_1 - p_3)\cdot \theta(z+ p_2 - p_4) +
c_2 \cdot \theta(z + p_1 - p_4)\cdot \theta(z+ p_2 - p_3) +
$$
\begin{equation}
c_3 \cdot \theta(z + p_1 + p_2 - p_3 - p_4)\cdot
\theta(z) = 0.
\label{88}
\end{equation}
This means that the following points
\begin{equation}
\Phi(\frac{p_1 + p_2 - p_3 - p_4}{2}),
\Phi(\frac{p_1 + p_3 - p_2 - p_4}{2}),
\Phi(\frac{p_1 + p_4 - p_2 - p_3}{2})
\label{89}
\end{equation}
of the Kummer variety
lie on the same line $\C P^1 \subset
\C P^{2^g-1}$ called a trisecant.}

For brevity here and below we denote $A(p_k)$ by
$p_k$.

The trisecant formula for a generic quadruple ${\bar p}$
follows from Propositions 4.4.2 and  5.2.1 and the Riemann theorem
on the theta divisor of a Jacobi variety.

In a generic case we can assume
$$
\dim \Theta_{-(p_1-p_3)} \cap \Theta_{-(p_1-p_4)} = g-2.
$$
In fact, this is always valid (\cite{Mm2}).
We are left to prove that
\begin{equation}
\Theta_{-(p_1-p_3)} \cap
\Theta_{-(p_1-p_4)} \subset \Theta_{-(p_1+p_2-p_3-p_4)} \cup \Theta.
\label{90}
\end{equation}
Assume that $\theta(u-p_3)=\theta(u-p_4)=0$ and $u=z+p_1$.  If
$\theta(u-p)=0$ for every $p \in \Gamma$ then $\theta(u)=0$.
Otherwise $\theta(u-p)=0$ for $p = p_3, p_4, q_1,
\dots, q_{g-2}$ where $q_j \in \Gamma$ and $u = p_1 + p_2 + q_1 +
\dots + q_{g-2} + K_{\Gamma}$. However in this case $u' =
z + p_1 + p_2 - p_3 - p_4 = p_2 + q_1 + \dots + q_{g-2} + K_{\Gamma}$
and, by the Riemann theorem,
$\theta(u')=0$.
The formula (\ref{90}) is established and it suffices to apply
Proposition 4.4.2.

For an arbitrary quadruple ${\bar p}$
the existence of a trisecant which passes through points (\ref{89}) is
obtained from a generic case by passing to a limit.  Explicit
formula for $c_k({\bar p})$ can be derived from (\ref{88}) (see also
\cite{Fay1}).

\vskip2.5mm

\begin{center}
{\bf \S 6. Theta functions of Prym varieties of double
coverings with two branch points}
\end{center}

\vskip2.5mm

{\bf 6.1. Prym varieties of ramified coverings.}

Let $\Gamma_0$ be a Riemann surface of genus $g$ and $q_0$ and
let $q_{\infty}$ be a pair of distinct points on $\Gamma_0$.  Denote
by $\gamma_0$ and $\gamma_{\infty}$ the contours bounding small
disks with centres in $q_0$ and $q_{\infty}$.

There is a one-to-one correspondence between
ramified double coverings $\Gamma\rightarrow \Gamma_0$
with branch points at $q_0$ and $q_{\infty}$ and
homomorphisms
$$
\rho :  H_1(\Gamma \setminus \{q_0,q_{\infty}\};\Z) \rightarrow
\Z_2
$$
such that $\rho([\gamma_0])=\rho([\gamma_{\infty}])=1$.
This correspondence is as follows.
Let $\gamma$ be a contour on $\Gamma_0$.
Put $\rho([\gamma])=1$
if the pull-back of $\gamma$ on the covering is non-closed.
Otherwise, put $\rho([\gamma])=0$.  It follows from this correspondence
that a covering of $\Gamma_0$ which is ramified at a pair of fixed
points is determined from these data in up to $2^{2g}$ possible ways.

Denote by ${\cal M}_{g,2}$ the moduli space of
Riemann surfaces of genus $g$ with pair of marked points and
denote by ${\cal DR}_g$ the moduli space of ramified double coverings
of surfaces of genus $g$ with two branch points.
The mapping
$\Gamma \rightarrow (\Gamma_0, q_0, q_{\infty})$ is
the $2^{2g}$-sheeted covering
${\cal DR}_g \rightarrow {\cal M}_{g,2}$.
It follows from
(\ref{80}) that $\dim {\cal DR}_g \ (\ = \dim {\cal
M}_{g,2} = \dim {\cal M}_g + 2)\ = 3g-1$ for $g \geq 2$.

Let $\pi : \Gamma\rightarrow \Gamma_0$ be a ramified covering
with branch points at $q_0$ and $q_{\infty}$.  Denote by
$\sigma :  \Gamma\rightarrow \Gamma$ the involution permuting
the branches of the covering.  It is clear that $\Gamma_0 =
\Gamma/\sigma$ and the involution $\Gamma\stackrel
{\sigma}{\longrightarrow}\Gamma$ uniquely determines the covering.

There exists a canonical basis for $1$-cycles
$a_1,\dots,a_{2g},b_1,\dots,b_{2g}$ on the surface
$\Gamma$ of genus $2g$ such that

\begin{equation}
\sigma (a_k) + a_{g+k} = \sigma(b_k) +
b_{g+k} = 0, \ \ 1 \leq k \leq g,
\label{91}
\end{equation}

\noindent
and the cycles $a_1,\dots,a_g,b_1,\dots,b_g$ are projected into the
canonical basis for cycles
$\{{\tilde a}_k,{\tilde b}_k\}$ on $\Gamma_0$ where ${\tilde
a}_k - \pi_*(a_k) = {\tilde b}_k - \pi_*(b_k)=0$.

Let $\omega_1,\dots,\omega_{2g}$ be a basis
for holomorphic differentials on $\Gamma$
and let ${\tilde \omega}_1,\dots{\tilde \omega}_g$
be a basis for holomorphic differentials on
$\Gamma_0$ which are related to the
canonical bases for cycles by (\ref{78}).
It follows from (\ref{91}) that
\begin{equation}
\sigma^*(\omega_k) = -\omega_{g+k} \ \ \
{\mbox {and}}\ \ \
\omega_k - \omega_{g+k} = \pi^*({\tilde \omega}_k)\ \ \
{\mbox {for}}\   1 \leq k \leq g.
\label{92}
\end{equation}

A differential $\omega$ is called a {\it Prym
differential} if $\sigma^*(\omega)=-\omega$.
The set of the forms $u_1 = \omega_1+ \omega_{g+1},\dots, u_g = \omega_g +
\omega_{2g}$ is a basis for holomorphic Prym differentials.
Consider the matrix of their $b$-periods
\begin{equation}
\Pi(\Gamma,\sigma)_{jk}\ (= \Pi_{jk}) \ =
\int_{b_k} u_j, \ \ 1 \leq j,k \leq g.
\label{93}
\end{equation}
This matrix is symmetric and, by (\ref{91}--\ref{92}),
relates to the period matrices of
$\Gamma$ and $\Gamma_0$ as follows
\begin{equation}
B(\Gamma) = \frac{1}{2} \left( \begin{array}{cc}
\Pi+B_0 &
\Pi-B_0 \\
\Pi-B_0 &
\Pi+B_0
\end{array}
\right)
\label{94}
\end{equation}
where
$$
B(\Gamma)_{jk} \ (= B_{jk}) \ = \int_{b_k}\omega_j, \ \
B(\Gamma_0)_{jk} \ (= B_{0jk})\ = \int_{\pi_*(b_k)}{\tilde \omega}_j.
$$

The principally-polarised Abelian variety
$$
Pr(\Gamma,\sigma) = \C^g / \{ \Z^g + \Pi\Z^g\}
$$
is called the {\it Prym variety} of the covering
$\Gamma\rightarrow\Gamma/\sigma$.

The Prym variety is defined for any involution
$\sigma$ on $\Gamma$ in an invariant manner as follows.
Let $A : \Gamma \rightarrow J(\Gamma)$ be the Abel mapping with
initial point ${\bar q}$ and let
$\sigma : \Gamma\rightarrow \Gamma$ be
an involution.  Then $J(\Gamma) = A(S^{n}\Gamma)$ where $n$ is the
genus of $\Gamma$ and the involution $\sigma$ induces the
following involution $\sigma :  J(\Gamma) \rightarrow J(\Gamma)$:
\begin{equation}
\sigma(\int^{q_1}_{\bar q} \omega +
\dots + \int^{q_{2g}}_{\bar q} \omega) =
\int^{\sigma(q_1)}_{\sigma(\bar q)} \omega + \dots +
\int^{\sigma(q_{2g})}_{\sigma(\bar q)} \omega.
\label{95}
\end{equation}
The {\it Prym variety} is an anti-invariant subvariety
\begin{equation}
Pr(\Gamma,\sigma) = \{ z \in J(\Gamma) : \sigma(z) = -z\}.
\label{96}
\end{equation}

In the present case it follows from (\ref{92})
that
\begin{equation}
\sigma(z_1,\dots,z_g,z_{g+1},\dots,z_{2g}) =
(-z_{g+1},\dots,-z_{2g},-z_1,\dots,-z_g)
\label{97}
\end{equation}
and the embedding
\begin{equation}
\varphi : Pr(\Gamma,\sigma) \rightarrow J(\Gamma), \ \
\varphi(z_1,\dots,z_g) = (z_1,\dots,z_g,z_1,\dots,z_g)
\label{98}
\end{equation}
is defined.
By (\ref{98}), derive the following proposition.

{\bf 6.1.1.}
{\sl Let $[\mu]$ be the principal polarization of
$J(\Gamma)$. Then the polarization
$\varphi^*([\mu])/2$ is principal.}

Hence the {\it Prym map}
\begin{equation}
Pr : {\cal DR}_g \rightarrow {\cal A}_g, \ \ \ \
(\Gamma,\sigma) \stackrel{Pr}{\longrightarrow}
(Pr(\Gamma,\sigma), \varphi^*([\mu])/2)
\label{99}
\end{equation}
is correctly defined.

Analogues of Proposition 6.1.1 are also valid for two classes of
double coverings, for  unramified and hyperelliptic
coverings.  In the first case the polarization, induced by an
embedding, is a multiple of the principal one. In the second case
a Prym variety coincides with the Jacobi variety of
a hyperelliptic Riemann surface which is a ramified double covering of
the $2$-sphere.

{\bf 6.2. Theta function of Prym variety.}

The Jacobi variety $J(\Gamma_0)$ is isomorphic to the subvariety
invariant under $\sigma$,
$J(\Gamma_0) = \{ z \in J(\Gamma) : \sigma(z)=z\}$, and it is
embedded into  $J(\Gamma)$ as follows
\begin{equation}
\pi^* :  J(\Gamma_0)\rightarrow J(\Gamma), \ \ \pi^*(z_1,\dots,z_g) =
(z_1,\dots,z_g,-z_1,\dots,-z_g).
\label{100}
\end{equation}
The projections
$$
\pi_1 : J(\Gamma) \rightarrow Pr(\Gamma,\sigma),\ \
\pi_2 : J(\Gamma) \rightarrow J(\Gamma_0)
$$
written in coordinates as
\begin{equation}
\pi_1(u,v) = u+v, \ \
\pi_2(u,v) = u-v, \ \
u,v \in \C^g
\label{101}
\end{equation}
are also given.
The composition
$Pr(\Gamma,\sigma) \times J(\Gamma_0)
\stackrel{\varphi + \pi^*}{\longrightarrow}
J(\Gamma)$
of the embeddings $\varphi$ and $\pi^*$
is an isogeny of degree $2^{2g}$. Notice that
$$
2e = \varphi\cdot \pi_1(e) +  \pi^*\cdot \pi_2(e).
$$
By (\ref{94}) and (\ref{101}), we have the following
proposition.

{\bf 6.2.1.}
{\sl 1)
\begin{equation}
\theta \big[
\begin{array}{cc}
 a & b \\
 c & d
\end{array}
\big] (z, B) =
\label{102}
\end{equation}
$$
\sum_{2e \in \Z^g/2\Z^g}
\theta \big[
\begin{array}{c}
(a+b)/2+e \\
c+d
\end{array}
\big] (\pi_1(z), 2\Pi) \cdot
\theta \big[
\begin{array}{c}
(a-b)/2+e \\
c-d
\end{array}
\big] (\pi_2(z), 2B_0)
$$
for all characteristics $a,b,c,d \in \R^g$.

\noindent
2) The formula
\begin{equation}
\frac{\theta(z,B)^2}
{\theta(\pi_1(z),\Pi)\cdot
\theta(\pi_2(z),B_0)}
\label{103}
\end{equation}
defines a meromorphic function on $J(\Gamma)$.}

The Abel--Prym mapping with initial point at $q_0$
is defined as follows
\begin{equation}
A_{Pr} :  \Gamma \rightarrow Pr(\Gamma,\sigma), \ \
A_{Pr}(q)_j =
\int^q_{q_0} u_j.
\label{104}
\end{equation}
Denote by $A_0$ and $A$
the Abel mappings for $\Gamma_0$ and $\Gamma$ with
initial points $\pi(q_0)$ and $q_0$.  By the
definitions of these mappings and (\ref{101}), we have
\begin{equation}
A_{Pr} = \pi_1\cdot A, \ \
A_0\cdot \pi =\pi_2\cdot A.
\label{105}
\end{equation}
Infer also that
\begin{equation}
C(\Gamma) = \pi^{-1}(C(\Gamma_0)) + q_0 +
q_{\infty}.
\label{106}
\end{equation}
By virtue of Proposition 5.2.1 applied to (\ref{103}),
Proposition 6.2.1, and (\ref{106}), the following proposition is
valid.

{\bf 6.2.2.}
{\sl Assume that $f(p) = \theta(A_{Pr}(p) - e,\Pi)$.
Then, for every $e \in Pr(\Gamma,\sigma)$, either of the two
possibilities holds

1) $f(p) \equiv 0$;

2) if $D$ is a divisor of zeros of the function $f(p)$
then
\begin{equation}
\varphi(e) = A(D) - A(q_0) -
A(q_{\infty}) + \pi^*(K_{\Gamma_0}),
\label{107}
\end{equation}
and this implies
\begin{equation}
D+\sigma(D)-q_0-q_{\infty} \approx C(\Gamma).
\label{108}
\end{equation}
Moreover a divisor $D$ satisfying (\ref{107}) is unique.}

Proposition 6.2.2 describes the inverse
(\ref{84}) of the Abel mapping restricted to the Prym
variety $\varphi (Pr(\Gamma,\sigma)) \subset J(\Gamma)$.

{\bf 6.3. Quadrisecants of Prym varieties.}

{\bf 6.3.1. Fay Quadrisecant Formula (\cite{Fay2}).}
{\sl For every quadruple $p_1, p_2, p_3,$ $p_4 \in \Gamma$
there exist constants $c_1({\bar p})$, $c_2({\bar p}), c_3({\bar
p})$, and $c_4({\bar p})$ such that
$$
c_1 \cdot \theta(z - p_1 - p_2)\cdot
\theta(z - p_3 - p_4) + c_2 \cdot \theta(z - p_1 -
p_3)\cdot \theta(z - p_2 - p_4) +
$$
\begin{equation} c_3 \cdot
\theta(z - p_1 - p_4)\cdot \theta(z - p_2 - p_3) + c_4 \cdot
\theta(z - p_1 - p_2 - p_3 - p_4)\cdot \theta(z) = 0.
\label{109}
\end{equation}
This means that the following points
$\Phi((p_1+p_2-p_3-p_4)/2),
\Phi((p_1+p_3-p_2-p_4)/2), \Phi((p_1+p_4-p_2-p_3)/2)$, and
$\Phi((p_1+p_2+p_3+p_4)/2)$
of the Kummer variety $K(Pr(\Gamma,\sigma),\Theta)$
lie on the same plane $\C P^2 \subset \C P^{2^g-1}$ called
a quadrisecant.}

Here $p_i$ stands for $A_{Pr}(p_i)$ and
$\theta(z)=\theta(z,\Pi(\Gamma,\sigma))$.

Here, as in \S 5.4, it suffices to prove that for
a generic quadruple of points the inclusion
\begin{equation}
X = \theta_{p_1+p_2} \cap \theta_{p_1+p_3} \cap
\theta_{p_1+p_4} \subset \theta \cup \theta_{p_1+p_2+p_3+p_4}
\label{110}
\end{equation}
holds and that $\dim X = g-3$.
Since the condition on the dimension of $X$ is not strong,
we sketch the proof of (\ref{110}).
We may assume that all the points $p_1,\dots,p_4$ are distinct.

Put $u = z - p_1$. It follows from Proposition 6.2.2 and
(\ref{110}) that one of the two possibilities holds:
1) $\theta(u-p) \equiv 0$, or
2) $p_2, p_3, p_4 \in D$ where $D$ is the divisor of zeros of the
function $f(p) = \theta(z-p)$.

By (\ref{105}), if $\theta(u-p) \equiv 0$ then
$z - p_1 - \sigma(p_1) = z$ and $\theta(z) =
\theta(u-\sigma(p_1)) = 0$.

We prove that if $\theta(u-p)$ does not vanish
identically  then $\theta(z-(p_1 + \dots + p_4))=0$.
Put $D_1 = D_0 - (p_2 + p_3 + p_4)$.

The set of divisors satisfying (\ref{108}) is a reducible
subvariety $M \subset S^{2g}\Gamma$.  The Abel mapping ${\hat A} : D
\rightarrow A(D-q_0-q_{\infty})+K_{\Gamma_0}$ is inverted on a
component $M_1 \subset M$ by (\ref{107}).  The set of divisors
$\pi(D)$ where $D \in M_1$ is a subvariety $M' \subset
S^{2g}\Gamma_0$.  The covering $M\rightarrow M'$ is ramified over the
subvariety $M'_0$ formed by divisors of zeros of Abelian
differentials.  It was proved (\cite{T91}) that
the end-points of the pull-back to $M_1$ of a closed path in $M'
\setminus M'_0$ differ on a divisor of the shape $D - \sigma(D)$
where $\deg D$ is even. Hence, conclude that $D_2 =
\sigma(p_2+p_3+p_4) + D_1$ does not lie in $M_1$ and, therefore,
the Abel mapping is non-invertible at ${\hat A}(D_2)$.

By Proposition 6.2.2 and the non-invertibility of the Abel
mapping at ${\hat A}(D_2)$, $\theta(u' - p) \equiv 0$ where $u' = {\hat
A}(D_2) = u - p_2 - p_3 - p_4$.  Thus derive $\theta(z - (p_1 + \dots+
p_4)) =0$.

The inclusion (\ref{110}) is established.

We thus outlined the proof of the quadrisecant formula
which is given in \cite{T91}. This proof, although similar
to the proof of the Fay trisecant formula given above,
shows the difference of the formulae from the geometric point of
view.

The formula (\ref{109}) was first obtained by Fay
(\cite{Fay2}) who applied the theory of the
Schottky--Jung relations.

\vskip2.5mm

\begin{center}
{\bf \S 7. Theta functions of Prym varieties of
unramified coverings}
\end{center}

\vskip2.5mm

{\bf 7.1. Prym varieties of unramified coverings.}

Let $\Gamma_0$ be a Riemann surface of genus $g$.
There is a one-to-one correspondence between
unramified double coverings
$\Gamma\rightarrow \Gamma_0$ and homomorphisms
\begin{equation}
\rho :  H_1(\Gamma; \Z)
\rightarrow \Z_2.
\label{111}
\end{equation}
This correspondence is as follows. Let $\gamma$ be a contour
on $\Gamma_0$. Put $\rho([\gamma])=1$ if its pull-back on the
covering is non-closed.  Otherwise, put $\rho([\gamma])=0$. It
follows from this correspondence that an unramified double covering is
determined by the surface $\Gamma_0$ up to $2^{2g+2}$ possibilities.
Moreover, for one of them, $\rho \equiv 0$,
the covering space has two components.

Denote by ${\cal DU}_g$ the moduli space of unramified double
coverings of surfaces of genus $g$. The mapping $\Gamma
\rightarrow \Gamma_0$ is a $2^{2g}$-sheeted covering
${\cal DU}_g \rightarrow {\cal M}_{g+1}$.  It follows from (\ref{80})
that $\dim {\cal DU}_g \ ( = \dim {\cal M}_g)\ = 3g-3$
for $g \geq 2$.

Let $\pi : \Gamma\rightarrow \Gamma_0$ be an unramified double covering.
Denote by $\sigma :  \Gamma\rightarrow \Gamma$
the involution permuting the branches of the covering.
It is clear that $\Gamma_0 = \Gamma/\sigma$ and the involution
$\Gamma\stackrel {\sigma}{\longrightarrow}\Gamma$ uniquely
determines the covering.

There exists a canonical basis
$a_0,b_0,a_1,\dots,a_{2(g-1)},b_1,\dots,b_{2(g-1)}$
for $1$-cycles on the surface $\Gamma$ of genus $2g-1$
such that
\begin{equation}
\sigma(a_0) - a_0 = \sigma(b_0) - b_0 = \sigma
(a_k) - a_{g-1+k} = \sigma(b_k) - b_{g-1+k} = 0
\label{112}
\end{equation}

\noindent
for $1 \leq k \leq g-1$
and the cycles $a_0,b_0,a_1,\dots,a_{g-1},b_1,\dots,b_{g-1}$
are projected into a canonical basis
$\{{\tilde a}_k,{\tilde b}_k\}$ for cycles on $\Gamma_0$ where
${\tilde a}_0 = \pi_*(a_0), {\tilde b}_0 = \pi_*(b_0)$ and ${\tilde
a}_k - \pi_*(a_k) = {\tilde b}_k - \pi_*(b_k)=0$ for
$1 \leq k \leq g-1$.

In this basis the homomorphism $\rho$ is simply written as
$$
\rho(a_j) = 0, \ \rho(b_j) = 0 \ \
{\mbox {for}} \ j \neq 0, \ \ \rho(b_0) = 1.
$$

Let $\omega_0, \omega_1, \dots, \omega_{2(g-1)}$ be a basis for
holomorphic differentials on $\Gamma$ and let ${\tilde \omega}_0$,
${\tilde \omega}_1, \dots, {\tilde\omega}_(g-1)$ be a basis for
holomorphic differentials on $\Gamma_0$ which are related
to the canonical bases for cycles by (\ref{78}).  It follows from
(\ref{112}) that $\sigma^*(\omega_0)=\omega_0, \ \sigma^*(\omega_k) =
\omega_{g-1+k}$ and $\omega_k + \omega_{g-1+k} = \pi^*({\tilde
\omega}_k)$ for $1 \leq k \leq g-1$.

Similarly as in the case of ramified coverings a differential
$\omega$ is called a {\it Prym differential} if
$\sigma^*(\omega)=-\omega$.  The set of the forms $u_k =
\omega_k+\omega_{g-1+k}$ where $1 \leq k \leq g-1$ is a basis
for holomorphic Prym differentials.

The Prym variety of the covering $\Gamma\rightarrow
\Gamma_0$ is defined by (\ref{96}) as in the case of ramified
coverings and is as follows
$$
Pr(\Gamma,\sigma) =
\C^{g-1}/\{\Z^{g-1} + \Pi(\Gamma,\sigma)\Z^{g-1}\}
$$
where
$\Pi(\Gamma,\sigma)$ is given by (\ref{93}).
However, the relation between the Jacobi and Prym varieties
is different:
\begin{equation}
B(\Gamma) =
\left( \begin{array}{ccc} T_0 & T_1 & T_1 \\ T_1^* & (\Pi+T_2)/2 &
(\Pi-T_2)/2 \\ T_1^* & (\Pi-T_2)/2 & (\Pi+T_2)/2 \end{array} \right)
\label{113}
\end{equation}
where
$$
B(\Gamma_0) = \left(
\begin{array}{cc}
T_0/2 & T_1 \\
T_1 & T_2
\end{array}
\right),\ \
T_0 = \int_{b_0}\omega_0, \ \
T_{1k} = \int_{b_k}\omega_0.
$$

The involution $\sigma : J(\Gamma)\rightarrow J(\Gamma)$
is written as
\begin{equation}
\sigma(z_0,z_1,\dots,z_{g-1},z_g,\dots,z_{2(g-1)}) =
(z_0,z_g,\dots,z_{2(g-1)},z_1,\dots,z_{g-1}),
\label{114}
\end{equation}
and the embedding of the Prym variety into $J(\Gamma)$
is given by
\begin{equation}
\varphi :
Pr(\Gamma,\sigma)\rightarrow J(\Gamma), \ \
\varphi(z_1,\dots,z_g) =
(0,z_1,\dots,z_{g-1},-z_1,\dots,-z_{g-1}).
\label{115}
\end{equation}
By (\ref{115}), derive the following proposition.

{\bf 7.1.1.}
{\sl Let $[\mu]$ be the principal polarization of $J(\Gamma)$. Then
$\varphi^*([\mu])/2$ is the principal polarization of the Prym
variety $Pr(\Gamma,\sigma)$.}

Thus as in the case of ramified coverings
the {\it Prym map} is defined by
\begin{equation}
Pr :  {\cal DU}_g \rightarrow {\cal A}_{g-1}, \ \ \ \
(\Gamma,\sigma) \stackrel{Pr}{\longrightarrow} (Pr(\Gamma,\sigma),
\varphi^*([\mu])/2).
\label{116}
\end{equation}

{\bf 7.2. Theta function of a Prym variety and the
quadrisecant formula for unramified coverings.}

The Jacobi variety $J(\Gamma_0)$ is the
$\sigma$-invariant subvariety of
$J(\Gamma)$ embedded as follows
\begin{equation}
\pi^* :  J(\Gamma_0)\rightarrow J(\Gamma), \ \
\pi^*(z_0,z_1,\dots,z_{g-1}) =
(2z_0,z_1,\dots,z_{g-1},z_1,\dots,z_{g-1}).
\label{117}
\end{equation}
The projections
$$
\pi_1 : J(\Gamma) \rightarrow
Pr(\Gamma,\sigma),\ \
\pi_2 : J(\Gamma) \rightarrow J(\Gamma_0),
$$
are defined by
\begin{equation}
\pi_1(u,v,w) = v-w, \ \
\pi_2(u,v,w) = (u,v+w), \ \
u \in \C, v,w \in \C^g.
\label{118}
\end{equation}
The composition
$Pr(\Gamma,\sigma) \times J(\Gamma_0)
\stackrel{\varphi + \pi^*}{\longrightarrow}
J(\Gamma)$
of the embeddings $\varphi$ and $\pi^*$
is an isogeny of degree $2^{2g-1}$ and
\begin{equation}
2e = \varphi\cdot \pi_1(e) +  \pi^*\cdot \pi_2(e).
\label{119}
\end{equation}
State the analogue of Proposition 6.2.1.

{\bf 7.2.1.}
{\sl For any $a_0,c_0 \in \R, a,b,c,d \in \R^{g-1}$
\begin{equation}
\theta \big[
\begin{array}{ccc}
a_0 & a & b \\
c_0 & c & d
\end{array}
\big] (z, B) =
\label{120}
\end{equation}
$$
\sum_{2e \in \Z^{g-1}/2\Z^{g-1}}
\theta \big[
\begin{array}{cc}
a_0 & (a+b)/2+e \\
c_0 & c+d
\end{array}
\big] (\pi_1(z), 2\Pi) \cdot
\theta \big[
\begin{array}{c}
(a-b)/2+e \\
c-d
\end{array}
\big] (\pi_2(z), 2B_0).
$$}

Since we have shown above how to prove the secant formulae, we
give the Beauville--Debarre formula omitting its proof.
For the same reason we omit the analogue of
Proposition 6.2.2 for unramified coverings.

Define the Abel--Prym mapping
$$
A_{Pr} : \Gamma \rightarrow Pr(\Gamma,\sigma)
$$
by
\begin{equation}
A_{Pr}(p)_j = \frac{1}{2}\int^p_{\sigma(p)} u_j.
\label{121}
\end{equation}

The following formula was obtained in \cite{BD1}.

{\bf 7.2.2. Beauville--Debarre Quadrisecant Formula.}
{\sl For every quadruple $p_1, p_2, p_3,$ $p_4 \in
\Gamma$ the following points
$\Phi((p_1+p_2-p_3-p_4)/2), \Phi((p_1+p_3-p_2-p_4)/2),
\Phi((p_1+p_4-p_2-p_3)/2)$ and $\Phi((p_1+p_2+p_3+p_4)/2)$
of the Kummer variety
$K(Pr(\Gamma,\sigma),\Theta)$ lie on the same complex projective
plane in $\C P^{2^g-1}$.}

\vskip2.5mm

\begin{center}
{\bf Chapter 3. Jacobi varieties and soliton equations
(the Riemann--Schottky problem, the Novikov conjecture,
and trisecants)}
\end{center}

\vskip2.5mm

{\bf \S 8. Baker--Akhieser functions and rings of commuting
differential operators}

\vskip2.5mm

{\bf 8.1. Finite-zone operators and Baker--Akhieser functions.}

A {\it Baker--Akhieser function}, on a surface $\Gamma$ of genus
$g$, corresponding to a point $q \in \Gamma$,
a local parameter $k^{-1}$ in a neighbourhood of $q$
($k(q) = \infty$), a polynomial $Q(k)$, and a divisor $D =
p_1 + \dots+p_g$ is a function $\psi (p)$ such that

1) $\psi$ is meromorphic on $\Gamma \setminus q$ and
$(\psi) \geq -D$ (see \S 4.1);

2) the function $\psi(p)\cdot \exp{(-Q(k))}$ is analytic in
a neighbourhood of $q$.

Denote by ${\cal S} = \{\Gamma, q, k^{-1}, Q(k), D\}$
the set of ``spectral  data'' and denote by $\Lambda
({\cal S})$ the space of all Baker--Akhieser functions corresponding
to these data.
The following proposition is a corollary of the Riemann--Roch
theorem.

{\bf 8.1.1.}
{\sl Let the divisor of poles of $\psi$ be non-special.
Then for a generic polynomial $Q$ the divisor of zeros of $\psi$ is
also non-special and $\dim \Lambda ({\cal S})$ = 1.}

This definition is naturally generalised for $n$-point or
vector functions. Baker--Akhieser functions were introduced by
Krichever (\cite{Kr1}) and now it is the main tool of the method of
finite-zone integration of soliton equations.

In the simplest case when $\Gamma$ is a hyperelliptic
surface $w^2 = E^{2g+1} + \dots$, $q$ is the point at  infinity,
$Q(k) = kx$, and $k = \sqrt{E}$,
the Baker--Akhieser function is the Bloch eigenfunction
of the following operator
\begin{equation}
L = \frac{d^2}{dx^2} + u (x).
\label{122}
\end{equation}
To be precise, under the assumption
that the operator (\ref{122}) is periodic, i.e., $u(x+T)=u(x)$,
we call $\psi(x,p)$ a {\it Bloch function} if
$$
L\psi=E(p)\cdot \psi;
$$
$$ \psi(x+T,p) = \mu(p)\cdot \psi(x,p).
$$
Here $p$ belongs to a double-sheeted covering of the
$E$-plane which , thus, parametrises the set of the common
eigenfunctions of the operator (\ref{122}) and of
the translation operator $f(x) \rightarrow f(x+T)$.

If a covering $(E,\mu) \rightarrow E$ has finitely many
branch points, then it is shown that this number is odd and
the covering space is a hyperelliptic surface
$$
w^2 = E^{2g+1} + \dots
$$
on which the Bloch function $\psi(x,p)$ is defined.
The Bloch function has the following asymptotic expansion
$\psi(x,p)
\sim \exp{(\sqrt{E}x)}\cdot (1 + O( E^{-1/2}))$
as $E \rightarrow \infty$ (\cite{DMN}).

{\bf Definition.}
{\sl The operator (\ref{122}) is called finite-zone ($g$-zone) if
its Bloch function is defined on a Riemann surface of finite genus
(of genus $g$).}

Generally, an operator is called {\it
finite-zone} if it has an eigenfunction being  a Baker--Akhieser
function which can be multi-point and vector-valued.

{\bf 8.2. Finite-zone solutions to the Korteweg--de Vries
and Kadomtsev--Petviashvili equations.}

The Korteweg--de Vries equation (KdV)
\begin{equation}
u_t = \frac{1}{4}(6uu_x + u_{xxx})
\label{123}
\end{equation}
is represented by an {\it $L, A$-pair},
the condition of commutation of two differential operators:
\begin{equation}
[L, -\frac{\partial}{\partial t}+A] = L_t + [L,A] = 0
\label{124}
\end{equation}
where
\begin{equation}
L = -\frac{\partial^2}{\partial x^2} + u(x,t), \ \
A = \frac{\partial^3}{\partial x^3} -
\frac{3}{2}u(x,t)\frac{\partial}{\partial x} -
\frac{3}{4}\frac{\partial u(x,t)}{\partial x}.
\label{125}
\end{equation}

There exists an infinite family of
differential operators $A_{2n+1}$ of the shape
$$
A_{2n+1} = \frac{\partial^{2n+1}}{\partial x^{2n+1}} +\dots,
$$
such that

1) it contains the operator $A=A_3$ (\ref{125});

2) the commutation conditions (\ref{124}) on
$A=A_{2n+1}$ generate
an infinite family of commuting flows on the space of potentials
$u(x)$. These flows are described by the non-linear equations of the
KdV type, the {\it KdV hierarchy} of non-linear equations
(\cite{L1}).

The stationary solutions have the following prominent
properties.

{\bf 8.2.1.} {\sl (Novikov (\cite{Nov1}),
Dubrovin and Novikov (\cite{DN}) }
{\sl The potential $u(x)$ of the operator (\ref{122})
is finite-zone if and only if it is a stationary solution
to an equation of the KdV hierarchy
\begin{equation}
[L, A_{2n+1} + c_{2n-1}\cdot
A_{2n-1} + \dots + c_1\cdot A_1] = 0.
\label{126}
\end{equation}
The operator $L$ has $g$ zones if (\ref{126}) holds for
$n=2g+1$ and fails for less values of $n$.
The finite-zone operators transform into
finite-zone ones under the KdV flows.}

In 1974 in the articles by Novikov, Dubrovin, Matveev, Its,
and Lax the class of finite-zone operators was introduced and studied
(\cite{Nov1,DN,IM,L15}). In 1975 some of these results were obtained
also in \cite{MM}.  Another approach to the periodic problem for
KdV was considered by Marchenko (\cite{Mar}). A detailed survey of
these papers was done in \cite{DMN}. Recently Shabat and
Veselov have found a very interesting characterization of finite-zone
operators (\ref{122}) in terms of the Darboux
transformations (\cite{VS}). For brevity we just explain the
procedure of constructing finite-zone solutions to the KdV
equation (Theorems 8.2.2 and 8.2.3).

Let $\Gamma = \{w^2 = E^{2g+1} + \dots\}$ be a hyperelliptic
surface, let $q=(E=\infty)$ be the point at infinity, and let
$k = \sqrt{E}$. Put $Q(k) = kx$.
For a non-special divisor $D = p_1 +
\dots +p_g$ of degree $g$ the space $\Lambda ({\cal S})$
is one-dimensional where ${\cal S} = \{\Gamma,q,k^{-1},Q,D\}$.

A set $\Sigma  = \{\Gamma,q,k^{-1},D\}$ is called
the {\it ``spectral data''} of the Schr\"odinger operator (\ref{122}).

{\bf 8.2.2.}
{\sl Generally, when $\dim \Lambda ({\cal S}) = 1$,
there exists a unique operator $L$ of the shape (\ref{122}) such that
$L\psi=E\psi$ where  $\psi \in \Lambda ({\cal S})$ and $\psi(x,k)
= \exp{(kx)} (1 + O(k^{-1}))$.}

Prove this proposition. Put
\begin{equation}
\psi(x,k) = \exp{(kx)}\cdot (1 + \frac{\xi(x)}{k} + O(\frac{1}{k^2})).
\label{127}
\end{equation}
Then for $u=-2\xi_x$ the function $(L-k^2)\psi$
does not belong to $\Lambda ({\cal S})$
and $\dim \Lambda ({\cal S}) = 1$ implies
$(L-k^2)\psi = \lambda\cdot \psi$ with $\lambda$ a constant.
However, $\exp{(-kx)}\cdot (L-k^2)\psi \rightarrow 0$ as
$E\rightarrow \infty$. This implies that
$\lambda=0$ and  $L\psi = E\psi$ for $u = -2\xi_x$.

In fact, we obtain the solution to the inverse problem, i.e., the
procedure of reconstructing $L$ from its ``spectral data''
${\cal S}$:
\begin{equation}
{\cal S} \rightarrow u(x) = -2\xi_x
\label{128}
\end{equation}

Finite-zone solutions to the KdV equation are obtained
as follows. Change the spectral data ${\cal S}$ by
substituting   $Q_1(k) = kx+k^3t$ for
$Q(k)=kx$ and construct from the new spectral data
${\cal S'}$ the function $\psi(x,t,k)$ such that
\begin{equation}
\psi(x,k) = \exp{(kx+k^3t)}\cdot (1 +
\frac{\xi(x,t)}{k} + O(\frac{1}{k^2})).
\label{129}
\end{equation}
Put
\begin{equation}
u(x,t) = -2\xi_x.
\label{130}
\end{equation}
The function $f=(-\partial_t+A)\psi$ belongs to
$\Lambda({\cal S'})$ and, since
$\exp{(-kx-k^3t)}\cdot (-\partial_t+A)\psi \rightarrow
0$ as $E\rightarrow\infty$, this function vanishes
identically.  Hence $[L,-\partial_t+A]\psi=0$ but the operator
$[L,-\partial_t+A]$ is the operator of multiplication by
the scalar function $f(x,t)
= (u_t - (6uu_x + u_{xxx})/4)$ and $\psi$
is non-trivial. We conclude that $f$ vanishes
identically and thus $u(x,t) = -2\xi_x$
satisfies the KdV equation (\ref{123}).
This proves the following proposition.

{\bf 8.2.3.}
{\sl Given a set of spectral data of generic
position, the finite-zone Schr\"odinger operator (\ref{127}--\ref{128})
and the finite-zone solution to the KdV equation (\ref{129}--\ref{130})
are constructed uniquely.}

It was shown in \cite{IM} that this solution is expressed by the
theta function formula
\begin{equation}
u(x,t) =  2
\partial^2_x \log{\theta(Ux+Wt+Z,\Omega)}
\label{131}
\end{equation}
where $U$ and $W$ are constant vectors, $\Omega =
B(\Gamma)$, and $Z = -A(D)-K_{\Gamma}$.
The inference of this formula is now routine and is
explained, for instance, in the surveys \cite{DMN,Kr2,Dubr1}.
The function $\psi$ is the Bloch eigenfunction of $L$ and together with
the ``spectral data'' is explicitly constructed from it
(\cite{DMN}).

The Kadomtsev--Petviashvili equation (KP)
\begin{equation}
\frac{3}{4}u_{yy} =
\frac{\partial}{\partial x} (u_t - \frac{1}{4}(6uu_x + u_{xxx})) = 0
\label{132}
\end{equation}
is the commutation condition
$$
[-\frac{\partial}{\partial_y} + L,
-\frac{\partial}{\partial_t} + A] = 0
$$
where
$$
L = \frac{\partial^2}{\partial x^2} + u(x,y,t), \ \
A = \frac{\partial^3}{\partial x^3} + \frac{3}{2}
u\frac{\partial}{\partial x} + w.
$$
This equation is a two-dimensional generalization of the KdV
equation and degenerates into it if $u(x,y)$ does not
depend on $y$.

{\bf 8.2.4.} {\sl (Krichever (\cite{Kr1,Kr2}))}
{\sl Let an arbitrary Riemann surface $\Gamma$ of genus $g$,
a point $q \in \Gamma$, a local parameter $k^{-1}$ in a
neighbourhood of $q$,  a non-special effective divisor $D = p_1 +
\dots + p_g$ of degree $g$,  and  a polynomial $Q(k) = kx + k^2y +
k^3t$ form the set of spectral data ${\cal S} = \{ \Gamma, q, k^{-1},
Q(k), D\}$.  Generally, when $\dim \Lambda ({\cal S})=1$, the
function $\psi \in \Lambda ({\cal S})$ is defined uniquely by the
condition $\exp{(-Q(k))}\cdot \psi(x,y,t,k)
\rightarrow 1$ as $k \rightarrow \infty$.
The operators $L$ and $A$ of the shape (\ref{122})
satisfying
\begin{equation}
(-\frac{\partial}{\partial y} + L)\psi =
(-\frac{\partial}{\partial t} + A)\psi = 0
\label{133}
\end{equation}
are reconstructed uniquely from $\psi$.
The reconstruction formulae are written as
\begin{equation}
u = -2\frac{\partial \xi_1}{\partial x}, \ \
w = 3\xi_1 \frac{\partial \xi_1}{\partial x} + 3 \frac{\partial^3
\xi_1}{\partial x^3} - 3 \frac{\partial \xi_2}{\partial x}
\label{134}
\end{equation}
where
\begin{equation}
\psi(x,k) =
\exp{(kx+k^2y+k^3t)}\cdot (1 + \frac{\xi_1(x,t)}{k} +
\frac{\xi_2(x,t)}{k^2} + O(\frac{1}{k^3})).
\label{135}
\end{equation}
Moreover the function $u(x,y,t)$ is a solution to the KP
equation written as
\begin{equation}
u(x,y,t) = 2\partial^2_x
\log{\theta(Ux+Vy+Wt+Z,\Omega)}
\label{136}
\end{equation}
where $U, V$, and $W$ are constant vectors, $\Omega =
B(\Gamma)$, and $Z=-A(D)-K_{\Gamma}$.}

The proof of this theorem is similar to the proofs of Theorems
8.2.2 and 8.2.3 given above.  It is essential that for the
operator (\ref{122})
the direct problem of constructing its Bloch function does not possess
an effective solution. The direct problem is solved effectively
for the matrix Dirac operator and in the middle of the 70s
this enables construction of finite-zone solutions to the non-linear
Schr\"odinger equation (Its and Kotlyarov) and the sine--Gordon
equation (Kozel and Kotlyarov). These solutions are expressed in terms
of theta functions of hyperelliptic surfaces (see the survey
\cite{DMN}).  An axiomatization of the analytic properties of Bloch
functions in the definition of Baker--Akhieser functions enables
construction of a wide class of solutions without solving the
direct problem. The KP equation was the first one to which this
scheme, the Krichever method, was applied.  Moreover, its finite-zone
solutions are constructed from all Riemann surfaces.

{\bf 8.3. Commutative rings of ordinary linear
differential operators.}

By Theorem 8.2.1, the problem of describing commutative rings of
ordinary linear differential operators is related to that
for finite-zone operators and was posed by
Novikov. Its effective solution for commuting operators
of rank 1 was obtained by Krichever (\cite{Kr1,Kr3}) who
repeat in this case some results of Burchnall and Chaundy
(\cite{BC}).

Here we mean by operators ordinary linear differential
operators.

A pair of commuting operators $L_1$ and
$L_2$ has {\it rank $l$} if for a generic pair of eigenvalues
$(\lambda_1, \lambda_2)$ the space of common
eigenfunctions corresponding to these eigenvalues is
$l$-dimensional. In this case $l$ is the greatest common
divisor of the orders of $L_1$ and $L_2$.

As before we mean by ``spectral data'' $\Sigma$
the set formed by
a Riemann surface $\Gamma$ of genus $g$,
a point $q \in \Gamma$, a local parameter $k^{-1}$ in a
neighbourhood of $q$, and a generic effective divisor
$D$ of degree $g$. Let $\psi(x,p)$ be the Baker--Akhieser function
corresponding to $\Sigma$ and normalised uniquely by the
following condition: $\exp{(-kx)}\cdot
\psi(x,p) \rightarrow 1$ for $p\rightarrow q$.

{\bf 8.3.1.}
{\sl Let $f$ be a meromorphic function on $\Gamma$ with a single pole
$q$ and let  $Q(k)$ be a polynomial such that $f(k) = Q(k) + O(1)$.
Then an operator
$L_f$ such that $L_f \psi(x,p) = f(p)\psi(x,p)$
is constructed uniquely from $f$.}

The proof of Theorem 8.3.1 is similar to the proofs of Theorems
8.2.2 and 8.2.4.
For the function $f$ with a pole of multiplicity $2$ at $q$ this is
exactly Theorem 8.2.2.
The operator
$L_f$ is reconstructed effectively
by formulae analogous to (\ref{130}) and (\ref{134}).

Denote by $M(q) = \cup_d {\cal L}(dq)$ the algebra of meromorphic
functions on $\Gamma$ having a single pole at $q$.

{\bf 8.3.2.}
{\sl The procedure $f \rightarrow  L_f$
sets an isomorphism of the algebra $M(q)$ to
a commutative algebra of ordinary differential
operators. Generic pairs, of operators from this algebra, are of
rank 1.}

The converse statement is also valid.

{\bf 8.3.3.}
{\sl Every maximal commutative algebra $R$ of rank 1 of ordinary
differential operators is isomorphic to the algebra
$\{L_f\} \approx M(q)$ for suitable spectral data $\Sigma$.}

The construction of $\{L_f\}$ is straightforward.
First, observe that the following proposition is valid.

{\bf 8.3.4.}
{\sl Let two ordinary differential operators
$L_1$ and $L_2$ commute. Then there exists a polynomial
$Q(x,y)$ such that $Q(L_1,L_2)=0$.}

Now the bundle of common eigenfunctions of $L_1$ and $L_2$ is
constructed on the Riemann surface
$\Gamma=\{Q(x,y)=0\}$ in terms of formal series.  This bundle is
one-dimensional and it is a Baker--Akhieser function
$\psi(x,p)$. For  completeness, we need to consider
Riemann surfaces with simplest singularities but with non-singular
marked points $q$. This situation was analysed separately in
[14.III]. The surface $\Gamma\setminus q$ is defined
in invariant form as the spectrum of the ring  $R$,  $\Gamma \setminus
q = Spec R$.

The detailed proofs of these results are exposed in (\cite{Kr1,Kr3}).
Algebraizations of this construction and, in particular, its
generalization to the case of an arbitrary field of coefficients were
given in  \cite{Dr,Mm3}.  An effective description for
pairs of commuting operators of rank $l >1$ is still not available
(see the surveys \cite{KN2,Nov2,PW}).  Krichever
and Novikov found the generic pair of commuting operators of rank 2
corresponding to a surface of genus $1$ (\cite{KN1,KN2}, see also
\cite{Grin1,NG}) and using their methods Mokhov described pairs of
rank $3$ corresponding to a surface of genus $1$ (\cite{Mokh1,Mokh2}).
The case of matrix ordinary differential operators was examined
by Grinevich (\cite{Grin2}).

\vskip2.5mm

{\bf \S 9. The Novikov conjecture in the Riemann--Schottky
problem and trisecants of Jacobi varieties}

\vskip2.5mm

{\bf 9.1. Effectivization of theta function formulae
for finite-zone solutions.}

The vectors $U, V$, and $W$ and the constant $C$ entering into the theta
functional formulae (\ref{131}) and (\ref{136}) for solutions to the
KdV and KP equations are determined by spectral data via
transcendent equations. Novikov proposed an effective procedure for
writing out these formulae. Namely, he proposed to find
the relations for a Riemann matrix $\Omega =
B(\Gamma) \in \A_g$ and the vectors $U, V$, and $W$ by
straightforward substitution of the theta functional formulae
into the KdV and KP equations. Therewith a system of algebraic
equations for $U, V$, and $W$ is obtained. The compatibility
conditions on this system gives relations for the Riemann matrix.

For $g = 2, 3$
there are no relations for the Riemann matrix. This is a
corollary of the following theorem.

{\bf 9.1.1. Torelli Theorem.}
{\sl The mapping $B : {\cal M} \rightarrow  \A_g$
(\ref{81}) assigning the Jacobi variety to a Riemann surface is
an embedding.}

A geometric proof of this theorem is exposed, for instance,
in \cite{GH}, and a soliton proof is given in \cite{Dubr1}. This
theorem implies that $\dim A_g = \dim {\cal
M}_g$ for $g = 2,3$ and thus a generic matrix
$\Omega \in {\cal H}_g$ is the period matrix of a Riemann
surface. Indicate a single constraint on Jacobi varieties.

{\bf 9.1.2. Mertens Theorem.}
{\sl The Jacobi variety of a Riemann surface is
irreducible, i.e., it is not a direct product of
Abelian varieties of positive dimension.}

The irreducibility condition is expressed effectively in terms
of theta constants, the values of the theta function and its
derivatives at zero.

{\bf 9.1.3.} {\sl (\cite{Sasaki})}
{\sl A principally-polarised Abelian variety
$\C^g / \{\Z^g + \Omega\Z^g\}$ is irreducible if and only
if the rank of the matrix
\begin{equation}
\left( \begin{array}{cccccc} {\hat
\theta}_{11}[n_1,0] & \dots & {\hat \theta}_{jk}[n_1,0] & \dots &
{\hat \theta}_{gg}[n_1,0] &
{\hat \theta} [n_1,0] \\
\dots & \dots & \dots & \dots & \dots & \dots \\
{\hat \theta}_{11}[n_r,0] &
\dots &
{\hat \theta}_{jk}[n_r,0] &
\dots &
{\hat \theta}_{gg}[n_r,0] &
{\hat \theta} [n_r,0]
\end{array}
\right)
\label{137}
\end{equation}
is maximal, i.e., equals $g(g+1)/2+1$. Here
\begin{equation}
{\hat \theta}[\alpha,\beta] =
\theta[\alpha,\beta](0,2\Omega), \ \
{\hat \theta}_{jk}[\alpha,\beta] =
\frac{\partial^2 {\hat \theta}[\alpha,\beta](0,2\Omega)}
{\partial z^j \partial z^k},
\label{138}
\end{equation}
and $\{2n_j\} = \Z^g/2\Z^g, r=2^g$.}

The programme of effective construction of two- and three-zone solutions
to the KdV and KP equations was realised by Dubrovin. This programme
is based on finding relations for $U, V$, and $W$.

{\bf 9.1.4.} {\sl (Dubrovin (\cite{Dubr1}))}
{\sl Let an Abelian variety
$\C^g/\{\Z^g+\Omega\Z^g\}$ be irreducible.

1) The formula (\ref{131}) defines for every $Z$ a solution
to the Korteweg--de Vries equation if and only if
the relations
\begin{equation}
\partial^4_U {\hat \theta}[n,0] -
\partial_U \partial_W {\hat \theta}[n,0] + d {\hat \theta}[n,0] = 0,
\label{139}
\end{equation}
hold for every $n \in \frac{1}{2}\Z^g/2\Z^g$
with $d$ a constant;

2) The formula (\ref{136}) defines for every $Z$ a solution
to the Kadomtsev--Petviashvili equation
if and only if the relations
\begin{equation}
\partial^4_U {\hat \theta}[n,0] - \partial_U \partial_W
{\hat \theta}[n,0] +
\partial ^2_V {\hat \theta}[n,0] +
d {\hat \theta}[n,0] = 0,
\label{140}
\end{equation}
hold for every $n \in \frac{1}{2}\Z^g/2\Z^g$
with $d$ a constant.}

Dubrovin's proof is exposed in detail in the survey
(\cite{Dubr1}) and we just notice that  in order
to obtain (\ref{139}) and (\ref{140}) it suffices to substitute
the theta formulae into the non-linear equations and to apply the
binary addition theorem of Riemann and Theorem 3.2.1.  As it was
noticed later this substitution leads to the Hirota equations
(\ref{196}), for the theta function, reducing to (\ref{139}) and
(\ref{140}) by the theorems mentioned above.  Convenience of the
formulae (\ref{139}) and (\ref{140}) consists in possibility of
applying both to construction of explicit solutions and problems
of algebraic geometry.

By using (\ref{139}) and (\ref{140})
Dubrovin obtained a direct procedure of construction of
two- and three-zone solutions to the KP equation from an arbitrary
irreducible
principally-polarised Abelian variety of dimension $g=2$ or $g=3$.
Since all Riemann surface of genus $2$ are hyperelliptic,
in the case of the KdV equation this procedure works for
$g=2$ only. For $g=3$ it suffices to complete this procedure
by an effective distinguishing the Jacobi varieties of hyperelliptic
surfaces. This was also done in \cite{Dubr1}.

{\bf 9.2. The Novikov conjecture in the Riemann--Schottky
problem.}

For $g \geq 4$ the subvariety $J_g = B({\cal
M}_g)$ has a positive codimension, in $\A_g$, equal to
$(g-2)(g-3)/2$.  Thus the procedure of construction of solutions
to the KP equation by using the Dubrovin effectivization formulae
(\ref{140}) needs at the first step an effective description of the
subvariety $B({\cal M}_g)$.  Until the early 80s this was one of
the most prominent problems of algebraic geometry.

{\bf 9.2.1. Riemann--Schottky Problem.}
{\sl Find equations distinguishing the closure
of the set of Jacobi varieties, ${\bar J}_g =
{\overline{B({\cal M}_g)}}$, in the moduli space $\A_g$ of
principally-polarised Abelian varieties.}

Before the Novikov conjecture a non-trivial
relation was known only for $g=4$. By the end of the 19th century it
was found by Schottky who showed that this relation distinguishes the
subvariety, of codimension $1$ in $\A_4$, possibly
reducible and containing ${\bar J}_g$ (\S 10.1).  In 1981 Igusa
proved that this subvariety is irreducible and thus the Schottky
relation solves the Riemann--Schottky problem for $g=4$ (\cite{Ig2}).

By the end of the 70s Novikov made the following conjecture
initiated by the Krichever construction of finite-zone solutions to
the KP equation from an arbitrary Riemann surface (see \S 8.2).

{\bf 9.2.2. Novikov Conjecture.}
{\sl An irreducible principally-polarised Abelian variety
$M = \C^g/\{\Z^g+\Omega\Z^g\}$
is the Jacobi variety of a Riemann surface if and only if
there exist vectors $U \in \C^g \setminus \{0\}, V, W
\in \C^g$ such that the function
$$
{\cal U}_Z(x,y,t) =
2\frac{\partial^2}{\partial x^2} \log \theta(Ux+Vy+Wt+Z)
$$
is a solution to the KP equation for every $Z \in \C^g$.}

According to this conjecture the relations solving the
Riemann--Schottky problem are exactly the relations distinguishing
the matrices $\Omega \in {\cal H}_g$ for which
(\ref{140}) are compatible. Rather soon a partial proof
of this conjecture was obtained.

{\bf 9.2.3.} {\sl (Dubrovin (\cite{Dubr2}))}
{\sl The compatibility condition on
(\ref{140}) distinguishes the subvariety (probably
reducible) of dimension  $3g-3 = \dim {\cal M}_g$ in $\A_g$.}

The Dubrovin theorem as well as the Schottky theorem for $g=4$
solves the Riemann--Schottky problem locally but for all
$g \geq 4$.  In order to obtain a global solution it
suffices to prove that this subvariety is irreducible.

The complete of proof of the Novikov conjecture as well as the
solution to the Riemann--Schottky problem was obtained by
Shiota (\cite{Shi1}).  The important papers of
Arbarello and De Concini (\cite{AC1}) and Mulase
(\cite{Mul1}) preceded his work. We explain all these
results in  \S 10.

{\bf 9.3. The Fay trisecant formula and finite-zone solutions
to soliton equations.}

Before proceeding to discussion of the Riemann--Schottky problem,
we expose the method of constructing finite-zone solutions to
the KP equation by using the Fay trisecant formula. This
method was introduced by Mumford (\cite{Mm2}) in the early 80s. It
shows how the KP equation arises as a degeneration of a purely
algebro-geometric identity. We expose this method following
the paper \cite{AC2} where it was applied to the
Riemann--Schottky problem.

Denote by ${\vec \theta}(z)$ the vector
$(\theta[n_1,0](z,2\Omega), \dots , \theta[n_r,0](z,2\Omega))$
and denote by $u \wedge v$ the product of vectors $u,v \in \C^{2^g}$
in the Grassmann algebra generated by $\C^{2^g}$.

Put $\Omega=B(\Gamma)$ and
take  $p_1, p_2, p_3, p_4  \in \Gamma$
where $\Gamma$ is a Riemann surface of genus  $g$.
By the Fay trisecant formula,
\begin{equation}
{\vec \theta}(p_1+p_2-p_3-p_4) \wedge {\vec
\theta}(p_1+p_3-p_2-p_4) \wedge {\vec \theta}(p_1+p_4-p_2-p_3) = 0
\label{141}
\end{equation}
where $p_k$ stands for $A(p_k)$ and
$A$ is the Abel mapping (\ref{83}).

We may assume that $p_1$ is the initial point of
the Abel mapping, $A(p_1)=0$.  Take a local parameter
$k^{-1}$ in a neighbourhood of $p_1$ such that
$k(p_1)=\infty$ and expand the Abel mapping in a neighbourhood of
$p_1$ into the series in the powers of $k^{-1}$:
$$
A : k^{-1} \rightarrow \sum_{m=1}^{\infty}
\frac{U_m}{k^m}
$$
where $U_m \in \C^g$.  Denote by $D_m$
the derivative in the direction $U_m$
$$
D_m = U_m^1
\frac{\partial}{\partial z^1} + \dots + U_m^g
\frac{\partial}{\partial z^g}.  $$

Introduce the formal operator $T(k)$
\begin{equation}
T(k^{-1}) = \exp{(\sum_{m=1}^{\infty} \frac{D_m}{k^m})} =
\sum_{j \geq 0} \frac{S_j(D_1,\dots)}{k^j}.
\label{142}
\end{equation}

Consider the limit of
(\ref{141}) as $p_4\rightarrow p_1$:  $k^{-1}(p_4) = \alpha
\rightarrow 0$.  Expand the left-hand side of
(\ref{141}) into the series in the powers of
$\alpha$.  The first non-trivial term corresponds to
$\alpha$
\begin{equation}
{\vec \theta}(p_2-p_3) \wedge D_1{\vec \theta}(p_2-p_3)
\wedge {\vec \theta}(p_2+p_3) = 0.
\label{143}
\end{equation}

Consider now the limit of (\ref{143}) as $p_3\rightarrow p_1$:
$k^{-1}(p_3)=\beta  \rightarrow 0$ and expand the left-hand side
of (\ref{143}) into the series in the powers of $\beta$.
The first non-trivial term corresponds to
$\beta^2$:
\begin{equation}
{\vec \theta}(p_2) \wedge D_1{\vec \theta}(p_2)
\wedge (D^2_1+D_2){\vec \theta}(p_2+p_3) = 0.
\label{144}
\end{equation}

Mumford indicates (\cite{Mm2}) that the identity
(\ref{144}) first appeared in the book of Fay
(\cite{Fay1}) (Corollary 2.13).

Consider the resulting limit of
$p_2\rightarrow p_1$,  $k^{-1}(p_2)=\varepsilon \rightarrow 0$, and
expand the left-hand side of (\ref{144}) into the series in the
powers of $\varepsilon$.  Every term of this series defines
a linear relation for values of the functions $\theta[n,0](z,2\Omega)$
and their partial derivatives at zero. Moreover  the coefficients of
this linear relation do not depend on the characteristic $n$.
Symbolically, these equations are written as
\begin{equation}
T(\varepsilon) \circ ({\vec \theta}(z) \wedge
D_1{\vec \theta}(z) \wedge
(D^2_1+D_2){\vec \theta}(z))|_{z=0} \equiv 0.
\label{145}
\end{equation}

Notice that we implicitly used the facts that
${\vec \theta}(u) = {\vec \theta}(-u)$
and that the derivatives of odd order of
$\theta[n,0](z,2\Omega)$ vanish at $z=0$.

As it was noticed in \cite{AC1}, (\ref{145})
reduces to the Dubrovin effectivization equations of the shape
(\ref{139}--\ref{140}).  Renormalising $D_j$,
if need be, rewrite the system
\begin{equation}
\Delta_j({\vec \theta}(z) \wedge D_1{\vec \theta}(z)
\wedge (D^2_1+D_2){\vec \theta}(z))|_{z=0} = 0 \ \
j \leq n
\label{146}
\end{equation}
where $\Delta_j = S_j(D_1,\dots)$, as
\begin{equation}
(\Delta_j\cdot D_1 - \Delta_{j-1}\cdot
(D^2_1+D_2) + \sum_{m=3}^j d_{m+1} \Delta_{j-m}) {\vec \theta}(z)
|_{z=0} = 0, \ \ j \leq n
\label{147}
\end{equation}
with $d_1, d_2, \dots $ constants.
Assign to every derivative $D_j$ its degree as follows:
\begin{equation}
\deg D_j = j.
\label{148}
\end{equation}
This agrees with the definition of $D_j$
as the derivative with respect to the $j$-th term of the
expansion of the Abel mapping. Then the principal term of
(\ref{147}) is written as
\begin{equation}
(\Delta_j\cdot D_1 - \Delta_{j-1}\cdot
(D^2_1+D_2)){\vec \theta}(z)|_{z=0} = 0.
\label{149}
\end{equation}
In particular, since the principal terms are invariant under
changes of the local parameter $k^{-1}$, these terms prevail.

{\bf 9.3.1.} {(\cite{Mm2,AC2})}
{\sl Equations (\ref{145}) and (\ref{147})
obtained from the Fay trisecant formula by
the degeneration $p_2, p_3, p_4 \rightarrow p_1$
are the effectivization equations for the theta
functional solutions to non-linear equations in the KP hierarchy.
The solutions to these equations are given by
\begin{equation}
{\cal U}_Z(x,t_1,t_2,\dots) = 2\frac{\partial^2}{\partial x^2} \log
\theta(U_1x+U_2t_1+U_3t_2+\dots + Z, B(\Gamma)).
\label{150}
\end{equation}
The simplest equation of this family coincides with the KP
equation and the finite-zone
solutions to it are described by
\begin{equation}
2\frac{\partial^2}{\partial x^2} \log
\theta(U_1x+\sqrt{3}U_2y+3U_3t + Z, B(\Gamma)).
\label{151}
\end{equation}
If $\Gamma$ is a hyperelliptic surface, $p_1$ is a
fixed point of the hyperelliptic involution $\sigma$, and a local
parameter $k^{-1}$ is inverted by the involution, $\sigma(k)=-k$,
then $U_2 = 0$ and the hierarchy  (\ref{151})
is the KdV hierarchy.}

The hierarchy of sine--Gordon equations is also obtained by
degenerating the trisecant formula. The finize-zone solutions
to the sine-Gordon equation
\begin{equation}
\frac{\partial^2 u}{\partial x\partial y} = \sin{u}
\label{152}
\end{equation}
are constructed from hyperelliptic Riemann surfaces.
Let $p_1$ and $p_2$ be different branch points of the hyperelliptic
covering $\Gamma\rightarrow \C P^1$.
The sine--Gordon equation and its finite-zone solutions are
obtained from the trisecant formula as $p_3
\rightarrow p_1$ and $p_4\rightarrow p_2$.

In neighbourhoods of $p_1$ and $p_2$
we introduce local parameters $k_1^{-1}$ and $k_2^{-1}$
which are inverted by the hyperelliptic involution,
$\sigma(k_j)=-k_j$.  Denote by $U_j$
the derivative of the Abel mapping
with respect to $k^{-1}_j$
at the point $p_j$ and denote by $D_j$ the derivative in
the direction $U_j$.

Consider now the limit of (\ref{141}) as
$p_4\rightarrow p_2$:  $k_2^{-1}(p_3)=\alpha \rightarrow 0$, and
expand the left-hand side of (\ref{141}) into the series in the
powers of $\alpha$. The term corresponding to $\alpha$ is
\begin{equation}
{\vec \theta}(p_1-p_3) \wedge
D_2{\vec \theta}(p_1-p_3) \wedge {\vec \theta}(p_1+p_3-2p_2) = 0.
\label{153}
\end{equation}
Put $k_1^{-1}(p_3) = \beta \rightarrow 0$ and expand the left-hand
side of (\ref{153}) into the series in the powers of $\beta$.
We obtain the series which starts with the linear term with the
coefficient
\begin{equation}
{\vec \theta}(0) \wedge D_1D_2{\vec
\theta}(0) \wedge {\vec \theta}(2\delta) = 0
\label{154}
\end{equation}
where $\delta = p_1-p_2$.  By (\ref{154}),
there exist constants $c_1, c_2$, and $c_3$ such that
$$
c_1 \cdot \theta[n,0](0,2\Omega)\cdot \theta[n,0](2z,2\Omega) + c_2
\cdot \theta[n,0](2\delta,2\Omega)\cdot \theta[n,0](2z,2\Omega) +
$$
\begin{equation}
c_3 \cdot D_1D_2\theta[n,0](0,2\Omega) \cdot
\theta[n,0](2z,2\Omega) = 0
\label{155}
\end{equation}
where
$\Omega=B(\Gamma)$.
Inspecting the Fay identity more thoroughly, we may obtain that
all the constants $c_j$ do not vanish. Apply the Riemann addition
theorem (\ref{52}) to  (\ref{155}) and derive
\begin{equation}
c_1  + c_2 \frac{\theta(z-\delta)\cdot \theta(z+\delta)}{\theta(z)^2} +
\frac{c_3}{2} D_1D_2 \log{\theta(z)} = 0
\label{156} 
\end{equation}
where $\theta(z)=\theta(z,B(\Gamma))$.

Since $p_1$ and $p_2$ are branch points of the hyperelliptic
covering, it is easy to show that $\delta$ is a half-period:
$\delta = (m_1 + B(\Gamma)m_2)/2, m_1,  m_2 \in \Z^g$
where $g$ is the genus of $\Gamma$.  Take (\ref{156}) for
$z=u$ and $z=u+\delta$ and subtract one from
another. After renormalization of constants we obtain
\begin{equation}
D_1D_2
\log{\frac{\theta(u+\delta)^2}{\theta(u)^2}} = C'\cdot
(\frac{\theta(u+2\delta)\theta(u)}{\theta(u+\delta)^2} -
\frac{\theta(u+\delta)\theta(u-\delta)}{\theta(u)^2}).
\label{157}
\end{equation}
Considering that $\delta$ is a half-period we
arrive at the following conclusion.

{\bf 9.3.2.} {\sl (\cite{Mm2})}
{\sl The function
\begin{equation}
{\cal U}_Z =
2\sqrt{-1} \log
{\frac{\theta(U_1x+U_2y+\delta+Z,B(\Gamma))}
{\theta(U_1x+U_2y+Z,B(\Gamma))}} -
2\pi(m,U_1x+U_2y+Z+\frac{\delta}{2})
\label{158}
\end{equation}
satisfies the sine--Gordon equation
$$
\frac{\partial^2 {\cal U}}{\partial x \partial y} = C\cdot \sin {\cal
U},
$$
where $C$ is a non-vanishing constant, for any $Z \in \Z^g$.}

The solution (\ref{158}) was found in  \cite{KozK}.
The identity (\ref{156}) which itself ``contains'' both
the sine--Gordon equation and
its finite-zone solutions was obtained by Fay in
\cite{Fay1} (Proposition 2.10).  The equations effectivising the
theta functional formulae for solutions to the sine--Gordon equation
were obtained in \cite{BE,DNat}.

Notice that, since , first, we considered not all possible
degenerations and, second, we just considered only first
non-trivial terms of vanishing identically series as while obtaining
(\ref{143}) and (\ref{144}), we gave evidently not all equations
which are contained in the trisecant formula.

\vskip2.5mm

{\bf \S 10. The Riemann--Schottky problem}

\vskip2.5mm

{\bf 10.1. The Schottky--Jung relations.}

Let $\Gamma_0$ be a Riemann surface of genus
$g$.  Consider an unramified double covering
$\Gamma \rightarrow \Gamma_0 = \Gamma/\sigma$.
As in \S 7, assign to this covering
a canonical basis for cycles such that the
homomorphism (\ref{111}) is written as
$$
\rho(a_j) = 0, \ \rho(b_j) = 0 \ \
{\mbox {for}} \ j \neq 0, \ \ \rho(b_0) = 1.
$$
For brevity denote $B(\Gamma)$ by $\tau$ and denote
$Pr(\Gamma,\sigma)$ by $\pi$.

We just consider the case when all characteristics have order $2$,
i.e. , if $\varepsilon = (\varepsilon^1, \dots,\varepsilon^{2n})$
is a characteristic then $2\varepsilon^j \in \{0, 1\}$.
In this case a characteristic
$(\varepsilon_1,\varepsilon_2)$, where $\varepsilon_j \in \R^n$,
of a theta function of $n$ variables is called
{\it even} if $4 \sum_j \varepsilon^j_1\varepsilon^j_2$ is even
and is called {\it odd} otherwise.

{\bf 10.1.1. Schottky--Jung Proportionality Relations.}
{\sl For even characteristics
$(\varepsilon_1,\varepsilon_2)$
the ratio
\begin{equation} \frac
{\theta \big[
\begin{array}{c}
\varepsilon_1 \\ \varepsilon_2
\end{array}
\big] (0, \pi)^2}
{\theta \big[
\begin{array}{cc}
0 & \varepsilon_1  \\
0 & \varepsilon_2
\end{array}
\big] (0, \tau) \cdot
\theta \big[
\begin{array}{cc}
0 & \varepsilon_1  \\
1 & \varepsilon_2
\end{array}
\big] (0, \tau)} = k^2(\pi,\tau)
\label{159}
\end{equation}
is a constant independent of the characteristic.}

This theorem was proved by Schottky for
$g=4$ in 1888 but in general the relations were introduced
in 1909 in the joint paper of Schottky and Jung
(\cite{SJ}, see also the survey \cite{F}) where they were
conjectured.  The proof of these relations was obtained by Farkas
and Rauch (\cite{FR}).

The scheme of applying them to the Riemann--Schottky problem
is as follows:

1) for a Riemann surface $\Gamma_0$ construct an unramified
double covering as it was shown above;

2) take a homogeneous identity met by all theta
functions, for instance, the ternary addition theorem of Riemann,
(\ref{54}), and apply it to the Prym theta
functions $\theta(z,\pi)$;

3) using (\ref{159}),
substitute into this identity for the values of
$\theta[\alpha,\beta](z,\pi)$ at zero the values
of theta functions of $\Gamma$ with characteristics
proportional to the first ones.

The resulting identity is not valid for all theta functions
since it is geometric rather than analytic due to its origin.

The first non-trivial identity was found by Schottky for $g=4$.
Its modern inference is given in
\cite{Clem}.  Igusa proved that this identity distinguishes ${\bar
J}_4$ exactly (\cite{Ig2}). Subsequently, it was proved that the
Schottky--Jung relations are sufficient for obtaining a local
solution to the Riemann--Schottky problem.

{\bf 10.1.2.} {\sl (van Geemen (\cite{vG}))}
{\sl For any $g \geq 4$ the Schottky--Jung relations
define in ${\cal A}_g$ a subvariety, of dimension
$3g-3$, containing ${\bar J}_g$ as one of irreducible
components.}

As Donagi showed, these subvarieties are generally reducible
and contain another but ${\bar J}_g$ components even
for $g=5$ (\cite{Don}).

{\bf 10.2. Geometry of the theta divisor of a Jacobi
variety and the Riemann--Schottky problem.}

All approaches to the Riemann--Schottky problem preceding the
Novikov conjecture were based on the description for the theta
divisor of a Jacobi variety given by Riemann. In two of these
approaches, this reveals rather brightly.

The first is based on investigation of
$Sing \Theta$, the set of singular points of the theta
divisor, which was also described by Riemann for Jacobi varieties.
By his theorem on singularities (see \cite{GH}), this set is rather
big. Its examination leads to the following result.

{\bf 10.2.1.} {\sl (Andreotti and Mayer (\cite{AM}))}
{\sl The condition $\dim Sing \Theta \geq g-4$ defines
in ${\cal A}_g$ a subvariety, of dimension $3g-3$, possibly
reducible and containing ${\bar J}_g$ as one of irreducible
components.}

In much the same ways as in the case of the Schottky--Jung
relations, these subvarieties are reducible in general:  just for $g=4$
such subvariety contains two components (Beauville, (\cite{Beau})).
Notice that the condition $\dim Sing \Theta = g-3$
defines in ${\bar J}_g$ exactly the Jacobi varieties of
hyperelliptic surfaces.

Another approach, going back to the papers of Lie and Wirtinger,
uses the representation of the theta divisor in the shape of
translated sum of $g-1$ examples of the embedded surface $\Gamma$:
\begin{equation}
\Theta = C  +
\dots + C  + \kappa
\label{160}
\end{equation}
where $C = A(\Gamma)$ is a surface embedded into  $J(\Gamma)$ and
$\kappa = K_{\Gamma}$.  A subvariety, of an Abelian
variety, of the shape (\ref{160}), is called a {\it translation
hypersurface}.  The investigation of such subvarieties was
started by Lie. It follows from the theorem of Lie and Wirtinger
that

{\bf 10.2.2.}
{\sl A principally-polarised Abelian variety is the Jacobi
variety of a non-hyperelliptic surface if and only if in a
neighbourhood of one of its non-singular points it has two distinct
representations of the shape (\ref{160}).}

The problem how to describe the hypothesis of Theorem 10.2.2
by equations was discussed by Poincare and Tchebotarev
but for the moment the complete solution to it is not obtained
(see the survey \cite{Little}).

{\bf 10.3. Trisecants of Jacobi varieties.}

One more approach, introduced by Gunning (\cite{G,G1}), is based
on the Fay trisecant formula.

Before exposing the Gunning theorem, we
mention the {\it Matsusaka--Hoyt criterion}
(\cite{Mats2,Hoyt}) which is, roughly speaking, as follows.
Let $M$ be an Abelian variety of complex dimension $g$, let
$X\subset M$ be an effective divisor, and let $C$ be a $2$-cycle
in $M$.  Denote by $D$ the Poincare duality operator
$D : H_{2g-k}(M;\Z) \rightarrow H^k(M;\Z)$. Put $\omega = D(X)$.
Then the following statements are equivalent:

a) $(\omega^g,[M]) = g !$ and $\omega^{g-1} = (g-1)! D(C)$;

b) $M$ is a direct product of the Jacobi
varieties of Riemann surfaces
$\Gamma_1, \dots , \Gamma_n$, and
$$
X = \cup_j (J(\Gamma_1) \times \dots \times J(\Gamma_{j-1}) \times \Theta_j
\times J(\Gamma_{j+1}) \times \dots \times J(\Gamma_n)) ,
$$
$$
C = \cup_j (0 \times \dots \times 0 \times \Gamma_j \times 0 \times \dots
\times 0)
$$
where $\Theta_j$ is the theta divisor of $J(\Gamma_j)$.

The Matsusaka--Hoyt criterion
also describes the closure ${\bar J}_g$ as
the set of direct products of Jacobi varieties.

Gunning's reasoning is based on the following
rigidity of trisecants of Jacobi varieties.
Let $M = \C^g/\{\Z^g + \Omega\Z^g\}$ be a principally-polarised
Abelian variety.  For every triple of points
$\alpha, \beta,\delta \in \C^g$ denote by
$V_{2\alpha,2\beta ,2\gamma}$ a set of vectors $t \in \C^g$ such that
the vectors $\Phi (t/2+\alpha), \Phi (t/2+\beta)$, and $\Phi
(t/2+\gamma)$ are coplanar, the Kummer map is defined by
(\ref{62}) and (\ref{70}).

{\bf 10.3.1.} {\sl (\cite{Fay1,G0})}
{\sl Assuming $M$ to be the Jacobi variety of $\Gamma$, we have
\begin{equation}
V_{2\alpha,2\beta
,2\gamma} = \{ t \in \C^g :  t + \alpha+\beta +\gamma \in {\tilde
\Gamma}\}
\label{161}
\end{equation}
where ${\tilde \Gamma}$ is a pre-image of $\Gamma = A(\Gamma) \subset M$
under the projection $\C^g \rightarrow M$.}

The trisecant identity (\ref{88}--\ref{89})
corresponds exactly to the case of $\alpha=p_2, \beta
=p_3, \gamma=p_4$, and $t = (p_1-p_2-p_3-p_4)$.  This implies
the following nice procedure of reconstructing the Riemann surface
$\Gamma \subset J(\Gamma)$ from a triple of
generic points of $\Gamma \subset J(\Gamma)$.

Gunning proved that if for an irreducible principally-polarised
Abelian variety $M$ and points $\alpha, \beta ,
\gamma \in M$ the set $V_{2\alpha, 2\beta ,2\gamma}$ is not empty,
then under some additional conditions we have that $V_{2\alpha,
2\beta ,2\gamma}$ contains a component $\Gamma$
such that $M = J(\Gamma)$. The following theorem was proved by
using the Matsusaka--Hoyt criterion.

{\bf 10.3.2. Theorem of Gunning.} {\sl (\cite{G})}
{\sl An irreducible principally-polarised Abelian variety
$M$ of dimension $g$ is the Jacobi variety of a Riemann
surface if and only if there exist points $\alpha,\beta ,\gamma \in
\C^g$ such that

1) they represent distinct points of
$M = \C^g/\{\Z^g + \Omega\Z^g\}$;

2) there are no complex multiplications
$F : M \rightarrow  M$ satisfying
$F(\alpha-\gamma)=F(\beta -\gamma)=0$;

3) $\dim_{(-\alpha-\beta)}V_{2\alpha,2\beta,2\gamma} > 0$.}

A complex multiplication is an endomorphism of $M$ induced by a
linear mapping $t \rightarrow ct+d$ where $t, d\in \C^g$ and
$c \in \C$. The validity of the hypothesis 3 implies existence of a
surface $\Gamma$ as a component of $V_{2\alpha, 2\beta, 2\gamma}$
passing through $(-\alpha-\beta)$.

The Gunning criterion was subsequently strengthened
by Welters (\cite{Welt0,Welt1,Welt2}) and Debarre (\cite{D2}).
Welters also derived its infinitesimal analogue.
To avoid plunging into the scheme theory
we formulate it roughly:  if $M$ is an irreducible
principally-polarised Abelian variety and
for some point $\alpha$ ``a
limit subvariety'' $\lim_{\beta ,\gamma\rightarrow
\alpha}V_{2\alpha,2\beta ,2\gamma}$ '' is of positive dimension
somewhere, then $M$ is the Jacobi variety of a Riemann surface.
This criterion was translated into the language of equations by
Arbarello and De Concini (\cite{AC1}):

{\bf 10.3.3.}
{\sl An irreducible principally-polarised Abelian variety $M$
of dimension $g$ is the Jacobi variety of a Riemann surface
if and only if there exist vectors  $U_1 \neq 0, U_2, \dots \in \C^g$
such that these vectors and the theta function of $M$
satisfy (\ref{145}).  Moreover just
it suffices to take the equations
\begin{equation}
T(\varepsilon) \circ ({\vec
\theta}(z) \wedge D_1{\vec \theta}(z) \wedge (D^2_1+D_2){\vec
\theta}(z))|_{z=0} \equiv \ 0 \ mod \ \varepsilon^{N+1}
\label{162}
\end{equation}
for $N = [(3/2)^g\cdot g!]$.}

This and Theorem 9.3.1 implies

{\bf 10.3.4.} {\sl (Arbarello and De Concini (\cite{AC1}, and Mulase
(\cite{Mul1}))}
{\sl An irreducible principally-polarised
Abelian variety $M$ of dimension $g$ is the Jacobi variety of
a Riemann surface if and only if
the functions (\ref{150}) satisfy the first $N(g) (<\infty)$ equations
of the KP hierarchy for any $Z \in \C^g$.}

We avoid exposing the formula for $N(g)$ because the Arbarello-De
Concini upper bound has already been given in Theorem 10.3.3
and Mulase (\cite{Mul1}), the author of the ``soliton'' proof of
Theorem 10.3.4, pointed out merely that such bound exists.  If the
upper bound for $N(g)$ is decreased to $N(g)=2$, then the
Novikov conjecture is proved. This was done by Shiota (see \S
10.5).

{\bf 10.4. On the soliton proof of Theorem 10.3.4.}

The higher KP equations are the equations met by the
function $u(x,t_1, \dots) = -2\xi_x$ (\ref{128}) with the spectral
data for the KP equation with
$Q(k) = kx + k^2t_1 + k^3t_2 +\dots$ (here we treat $y$
as a temporal variable $t_1$).

Give a strong definition of this hierarchy.
Let $R$ be the ring of smooth functions of $x, t_1, \dots$.
Introduce the algebra $E$ of formal
pseudo-differential operators
\begin{equation}
A = \sum_{-\infty}^{N<\infty} p_n
\left(\frac{d}{dx}\right)^n, \ \ p_n
\in R
\label{163}
\end{equation}
with multiplication defined by the Leibniz
rule
\begin{equation}
\left(\frac{d}{dx}\right)^k \cdot f = \sum_{j\geq 0}
\left(\begin{array}{c} k \\ j\end{array}\right) \frac{d^k f}{dx^k}
\cdot
\left(\frac{d}{dx}\right)^{k-j}, \ \ k \in \Z,\ f \in R.
\label{164}
\end{equation}
The order of the operator (\ref{163}), $ord A$,
equals $d$ if $p_d \neq 0$ and $p_j = 0$ for $j \geq d$.  Denote
by $E^-$ the subalgebra formed by operators of order $ <0$.
For every operator $A$ a decomposition $A = A^+ + A^-$ into the sum
of the differential operator $A^+$ and the operator $A^- \in E^-$
is defined.

Let $H$ be the affine subspace, of $E$,
formed by the operators
$L = d/dx + \psi$ where $\psi \in E^-$.  The {\it
Kadomtsev--Petviashvili hierarchy} is the following system of
evolution equations for $H$ with respect to
$t_1, \dots$:
\begin{equation}
\frac{\partial}{\partial t_n} L = [ B_n, L], \ \
B_n = (L^n)^+.
\label{165}
\end{equation}
The following proposition is derived from the definition of these
flows by simple reasonings.

{\bf 10.4.1.}
{\sl The KP hierarchy is equivalent to the following system
of the Zakharov--Shabat equations:
\begin{equation}
\frac{\partial B_m}{\partial t_n} - \frac{\partial B_n}{\partial t_m}
= [B_n, B_m].
\label{166}
\end{equation}
The vector fields on $H$ defined by the higher KP equations
commute.}

For the finite-zone solutions to the KP hierarchy
the action of the KP flows (\ref{168}) on  potentials
$u(x,t_1, \dots)$ of the operators $L$ is linearised by
(\ref{150}). Here we treat
the spatial coordinate $y$
in (\ref{150}) as the first among the temporal, $t_1$.

Put
\begin{equation}
L = \frac{d}{dx} + \frac{{\cal U}}{2}\left(\frac{d}{dx}\right)^{-1} +
\sum_{k=2} u_k\left(\frac{d}{dx}\right)^{-k}.
\label{167}
\end{equation}
If $L$ satisfies (\ref{165}), then all the functions
$u_k$ are expressed as polynomials in the function
${\cal U}$ and its derivatives.  Thus it is natural to call
${\cal U}$ a solution to the KP hierarchy.
Precisely this is meant in the statement of the Novikov
conjecture, which also implies
\begin{equation}
\frac{\partial L}{\partial \tau} = 0 \ \ {\mbox {for}} \ \
\frac{\partial {\cal U}}{\partial \tau} = 0
\label{168}
\end{equation}
where $\tau = c_1t_1 + \dots + c_nt_n, c_j \in \C$.

Now it becomes clear how to derive Theorem 10.3.4 from the
Burchnall--Chaundy--Krichever theory (\S 8.3).  Let
$\theta(z,\Omega)$ be
a theta function of $g$ variables and  the function
\begin{equation}
{\cal U}_Z(x,t_1,t_2, \dots) =
2\frac{\partial^2}{\partial x^2} \log \theta(Ux+V_1t_1+V_2t_2+\dots +
Z, \Omega))
\label{169}
\end{equation}
be a solution to the first $N(g)=2g+1$ equations of the KP hierarchy
for any $Z \in \C^g$.
Then there exist the following linear combinations
\begin{equation}
c_1V_1 + \dots + c_{2m}V_{2m} +
V_{2m+1} = d_1V_1 + \dots + d_{2n-1}V_{2n-1} + V_{2n} = 0
\label{170}
\end{equation}
such that $2m+1$ is relatively prime with $2n$.
Take the linear operators
\begin{equation}
B' = c_1B_1 + \dots + c_{2m}B_{2m} + B_{2m+1},  \
B'' = d_1B_1 + \dots + d_{2n-1}B_{2n-1} + B_{2n}.
\label{171}
\end{equation}
Put $\tau' = c_1t_1 + \dots + c_{2m}t_{2m} + t_{2m+1}$ and
$\tau'' = d_1t_1 + \dots + d_{2n-1}t_{2n-1} + t_{2n}$.
By (\ref{165}) and (\ref{168}) , we have
\begin{equation}
\frac{\partial L}{\partial \tau'} -
\frac{\partial L}{\partial \tau''} = [B',B''] = 0.
\label{172}
\end{equation}
We conclude from (\ref{172}) that, since the operators
$B'$ and $B''$ commute and are of coprime orders,
they are included into a maximal commutative algebra of
differential operators ${\cal A}$ of rank $1$. It follows from
Theorem 8.3.3 that this algebra is obtained by the
Burchnall--Chaundy--Krichever construction from a Riemann surface
$\Gamma$ and additional spectral data. It is left to prove that
$\Omega = B(\Gamma)$.

From the fact that the algebra
${\cal A}$ is maximal it is inferred
that the vectors
$U,V_1$, $V_2, \dots$ span  the tangent space to
$J(\Gamma)$. This was strictly proved in \cite{Shi1}.
Now, since the formal solution (\ref{169}) is
induced from the finite-dimensional torus $\C^g / \{\Z^g +
\Omega\Z^g\}$, we conclude that
$\Gamma$ is smooth. Since at the same time due to the Krichever
construction the finite-zone solutions (\ref{136}) and
(\ref{150}) to the KP hierarchy are constructed from the ``spectral
data'' of ${\cal A}$, it is easy to conclude that these solutions are
given exactly by (\ref{169}).

We omit the rigorous proofs (\cite{Mul1,Shi1}) because we
have explained the most important thing, i.e., how a Riemann surface
arises in the ``soliton'' proof.

{\bf 10.5. The Shiota theorem (the proof of the Novikov
conjecture).}

Shiota derived from Theorem 10.3.4 the proof of the Novikov
conjecture showing that if the function
\begin{equation}
{\cal U}_Z(x,t_1,t_2) =
2\frac{\partial^2}{\partial x^2} \log \theta(Ux+V_1t_1+V_2t_2 + Z,
\Omega)
\label{173}
\end{equation}
satisfies the first two equations of the KP hierarchy (after
the substitution $t_1 \rightarrow
y, t_2 \rightarrow t$ this function satisfies the KP equation
(\ref{132})) then we may successively construct
vectors $V_3, V_4, \dots$ such that
the function (\ref{169}) satisfies the higher KP equations
(\cite{Shi1}).  The proof of Shiota is analytical. In (\cite{AC2})
Arbarello and De Concini simplified his proof.

{\bf 10.5.1. Shiota Theorem.}
{\sl Let $M = \C^g /\Lambda$ be an irreducible
principally-polarised Abelian variety and
$\theta$ be its theta function.
Then (\ref{140})  are solvable
for $U, V, W \in \C^g$ with $U \neq 0$
if and only if $M$ is the Jacobi variety of a
Riemann surface.}

{\bf 10.6. On new approaches to describing ${\bar J}_g$.
The Welters trisecant conjecture.}

Although the Riemann--Schottky problem was completely solved by
Shiota, it still attracts attention of researchers.
We indicate above five
approaches to it: the Schottky--Jung relations, the geometry of $Sing
\Theta$, translation surfaces, trisecants, the Novikov conjecture,
and only the latter one led to a success. Theorem 10.3.3, of Welters,
Arbarello, and De Concini, points at the deep relation
between trisecants and soliton equations and is a forcible argument
for the following conjecture.

{\bf 10.6.1. Welters Conjecture.}
{\sl An irreducible principally-polarised Abelian variety
is a Jacobi variety if and only if it, or, more
precisely, its Kummer variety, admits a trisecant.}

Although the statement of the conjecture is exceptionally
strong it has already gained serious confirmations.

{\bf 10.6.2.} {\sl (Beauville and Debarre (\cite{BD2}))}
{\sl If an irreducible principally-polarised Abelian variety
admits a trisecant, then $Sing \Theta \geq g-4$.}

In view of Theorem 10.2.1 (of Andreotti and Mayer) this
means that the hypothesis of the Welters conjecture solves the
Riemann--Schottky problem locally, it
distinguishes a subvariety, in
${\cal A}_g$, which has ${\bar J}_g$ as one of its components.

Leaning upon Theorem 10.6.2 Debarre showed that

{\bf 10.6.3.} {\sl (Debarre (\cite{D1}))}
{\sl If an irreducible principally-polarised Abelian variety is the
Prym variety of an unramified double covering and admits a
trisecant, then it is the Jacobi variety of a Riemann surface.}

The set of such Prym varieties is a
$3g$-dimensional family and the Prym map
$Pr : {\cal DU}_{g+1} \rightarrow {\cal A}_g$ is of maximal rank
at a generic point (see \S 12). Hence ${\overline{Pr({\cal
DU}_{g+1})}}={\cal A}_g$ from $g=4,5$.  This implies that

{\bf 10.6.4.} {\sl (Debarre (\cite{D1}))}
{\sl The Welters conjecture is valid for four- and
five-dimensional varieties.}

Note two more approaches to describing ${\bar J}_g$:

1) A nice geometric property of ${\bar J}_g$ has been found recently
by Buser and Sarnak: for large $g$ the subvariety
${\bar J}_g$ lies in a small neighbourhood of the boundary of
${\cal A}_g$ (\cite{BS}).
Methods used in its proof are far from ones discussed in this
survey;

2)  Recently Buchstaber and Krichever have introduced a programme of
soliton description for ${\bar J}_g$ based on functional equations,
distinguishing Baker--Akhieser functions, and addition theorems for
theta functions. For the present , this programme is still
not completed (\cite{BK0,BK}).

Other geometric approaches to the Riemann--Schottky
problem, related the previous and based on study
of Kummer varieties, were considered in \cite{Mun,vG2}.

\vskip2.5mm

\begin{center}
{\bf Chapter 4. Prym varieties and soliton equations}
\end{center}

\vskip2.5mm

{\bf \S 11. Soliton equations and Prym theta functions}

\vskip2.5mm

{\bf 11.1. Finite-zone two-dimensional Schr\"odinger
operators and the Veselov--Novikov equation.}

A direct generalization, of the scheme of finite-zone
integration of the KdV equation, for $2+1$-equations is
hindered by some circumstances.

The KdV equation is integrated in terms of the theta
function of a surface $\Gamma$ parametrising the Bloch functions of
one-dimensional Schr\"odinger operator.  For $n$-dimensional
Schr\"odinger operator $\Delta +u$ the global Bloch
function is defined on an $n$-dimensional manifold $M^n$
as follows.  Assume for simplicity that
a potential  $u(x^1,\dots,x^n)$ is $\Z^n$-periodic.
Let $(\mu_1,\dots,\mu_n) \in \C^n$ and let $E \in \C$.
The condition of existence of a Bloch function
$\psi(x^1,\dots,x^n)$ with multiplicators
$\mu$ ($u(x^1,\dots,x^j+1,\dots,x^n) = \mu_j \cdot
u(x^1,\dots,x^j,\dots,x^n)$) and with eigenvalue
$E$ is written as the condition for the kernel of
the operator
$$
L_{\nu,E} = \Delta + \sum_j\nu_j
\partial_{x^j} + (u + \sum_j\nu_j^2 -E)
$$
to be non-zero in the space $W^2_2(\R^n/\Z^n)$
where $\mu_j = \exp{\eta _j}$,
This is equivalent to a property of the kernel of the operator
\begin{equation}
L_{\nu,E} \circ \Delta^{-1}
:  L_2(\R^n/\Z^n) \rightarrow L_2(\R^n/\Z^n).
\label{174}
\end{equation}
to be non-zero.
The operator (\ref{174}) is of the shape $1 + P(\nu,E)$ where $P$
is a compact operator being a polynomial in $\nu$ and $E$.
By the Keldysh theorem, non-invertibility of such
operator is equivalent  to  vanishing of its determinant defined for
such operators and this determinant is an entire function of $\nu$ and
$E$.  The manifold $M^n$ is the set of zeros of this determinant.

A natural analogue of the finite-zone condition is the
condition for $M^n$ to be  projective algebraic.
However, the scheme of integration whereby a non-linear equation
determines a deformation of a Bloch eigenfunction
does not work. A Bloch function on a Riemann surface defines a
line bundle from which an operator is reconstructed and which is
deformed by the equation (see \S 8.3).
This scheme works in the case of Riemann surfaces because their
Picard groups are non-discrete (see Theorem 5.3.3) and does not
work for $n=2$ because the Picard groups of 
``almost all'' two-dimensional
algebraic varieties are discrete.

We may obviate this difficulty for $n=2$, for instance, by
fixing an energy level $E = const$.  In this event $M^2$ is
changed for a one-dimensional complex surface, the {\it ``spectrum'' of
$L$ on the energy level  $E$}. This works for the KP equation with $L =
-\partial_y + \partial^2_{x} + u(x,y)$ (\cite{Kr89}). The fact that the
derivative with respect to $y$ enters the operator
$L$ linearly is not accidental.  Manakov showed that if two differential
operators of several variables form a non-trivial $L,A$-pair, then
derivatives with respect to all variables but of single one enter
linearly into each of the operators (\cite{Man}).

Manakov introduced a generalization of $L,A$-pairs,
the {\it $L,A,B$-triples}:
\begin{equation}
\frac{\partial L}{\partial t} = [L,A] + BL.
\label{175}
\end{equation}
Such deformations preserve the ``spectrum'' of $L$
on the zero energy level, $E=0$.

A natural candidate for a deformed operator  $L$
is the two-dimensional Schr\"odinger operator
\begin{equation}
\frac{\partial^2 }{\partial x^2}
+ \frac{\partial^2}{\partial y^2} + V(x,y).
\label{176}
\end{equation}
The development of the theory of two-dimensional
Schr\"odinger operators finite-zone on one energy level was
started in 1976 by Dubrovin, Krichever, and Novikov (\cite{DKN1})
by solving the inverse problem for these operators.

{\bf 11.1.1.}
{\sl 1) (Dubrovin, Krichever, and Novikov (\cite{DKN1}))
Let $\Gamma$ be a Riemann surface of genus $g$,
let $q_1$ and $q_2$ be a pair of distinct points of $\Gamma$,
let $k^{-1}_1$ and $k^{-1}_2$ be local parameters
defined in some neighbourhoods of $q_1$ and $q_2$, and let $D$ be an
effective non-special divisor of degree $g$.  The two-point
Baker--Akhieser function $\psi(z,{\bar z},p)$ corresponding to these
spectral data and with asymptotic expansions
\begin{equation}
\psi(z,{\bar z},p) = \exp{(k_1z)}\cdot (1 + \frac{\xi^+_1}{k_1}+
O(\frac{1}{k_1^2})) , \ p \rightarrow q_1 ,
\label{177}
\end{equation}
$$
\psi(z,{\bar z},p) = c(z,{\bar z})\cdot \exp{(k_2{\bar z})}
\cdot (1 + O(\frac{1}{k_1})) , \ p \rightarrow q_2
$$
is unique and presents an
eigenfunction of the operator
\begin{equation}
L = \partial{\bar \partial} + A(z,{\bar z})\cdot {\bar \partial} +
V (z,{\bar z})
\label{178}
\end{equation}
with the zero eigenvalue
\begin{equation}
L\psi=0.
\label{179}
\end{equation}
The coefficients of the operator are reconstructed by
the formulae
\begin{equation}
A(z,{\bar z}) = - \partial \log c(z,{\bar z}), \ \
V(z,{\bar z}) = - {\bar \partial} \xi^+_1.
\label{180}
\end{equation}
The potential $V(z,{\bar z})$ is written as
\begin{equation}
V(z,{\bar z}) = 2\partial {\bar \partial} \log
\theta({\tilde U}^+z + {\tilde U}^- {\bar z} + {\tilde Z}(D),
B(\Gamma)).
\label{181}
\end{equation}

2) (Cherednik (\cite{Cher})) Moreover, if
an anti-holomorphic involution
$\tau : \Gamma \rightarrow \Gamma$ acts such that
$\tau(q_1) = q_2, \tau(k_1) = -{\bar k}_2$, and $ \tau(D) + D - q_1
-q_2 \approx C(\Gamma)$, then the operator $L$ is self-adjoint.

3) (Veselov and Novikov (\cite{VN1}))
If there exists a holomorphic involution
$\sigma : \Gamma\rightarrow
\Gamma$ such that $\sigma(q_j)=q_j, \sigma(k_j)=-k_j$, and
\begin{equation}
\sigma(D) + D - q_1 - q_2 \approx C(\Gamma),
\label{182}
\end{equation}
then the operator $L$  is potential ,  $A \equiv 0$, and the
potential is expressed in terms of the Prym theta function of the
ramified covering $\Gamma \rightarrow \Gamma/\sigma$
as
\begin{equation}
V(z,{\bar z}) = 2\partial{\bar \partial} \log
\theta(U^+z + U^- {\bar z} + Z(D), Pr(\Gamma,\sigma)).
\label{183}
\end{equation}
These conditions for the operator to be real and to be potential are
compatible for $\sigma\tau=\tau\sigma$ and $\tau(D)=D$.
The equation (\ref{182}) is solvable for smooth Riemann surfaces
only if the covering $\Gamma\rightarrow \Gamma/\sigma$ has just two
branch points.}

Here $\partial = \partial /\partial z$ and ${\bar\partial} =
\partial /\partial {\bar z}$ where $z = x + \sqrt{-1}y \in \C$.
The inference of (\ref{181}) and (\ref{183}) is
analogous to that of (\ref{131}).

For smooth periodic operators the conditions for operator to
be potential found by Veselov and Novikov are necessary.
Different proofs of this were given in \cite{Kr89} and in the paper
of the author (Functional. Anal. Appl. {\bf 24}:1 (1990), 76--77).

The set  of the equations describing ``isospectral'' deformations of
potential operators is the {\it Veselov--Novikov hierarchy of
equations}.

{\bf 11.1.2.} {\sl (Veselov and Novikov (\cite{VN2}))}
{\sl The Baker--Akhieser function corresponding to the same
spectral data as in 11.1.1.3 and with the asymptotic expansions
\begin{equation}
\psi(z,{\bar z},t',t'',p) = \exp{(k_1z + \sum k_1^{2m+1}t'_m)}
\cdot (1 + \frac{\xi^+_1}{k_1}+
O(\frac{1}{k_1^2})) , \ p \rightarrow q_1 ,
\label{184}
\end{equation}
$$
\psi(z,{\bar z},t',t'',p) = c(z,{\bar z})\cdot
\exp{(k_2z + \sum k_2^{2n+1}t''_n)}\cdot (1 + O(\frac{1}{k_2})) , \ p
\rightarrow q_2
$$
is unique and satisfies the equations
\begin{equation}
(\frac{\partial }{\partial t'_m} - A'_m) \psi =
(\frac{\partial }{\partial t''_n} - A''_n) \psi = 0,
\label{185}
\end{equation}
where
$$
A'_m = \partial^{2m+1} +
a'_{2m-1} \partial^{2m-1} + \dots, \ \
A''_n = {\bar \partial}^{2n+1} +
a''_{2n-1}{\bar \partial}^{2n-1} + \dots.
$$
The potentials $V(z,{\bar z},t',t'')$ reconstructed by
(\ref{180}) satisfy the non-linear equations represented by
$L,A,B$-triples
\begin{equation}
\frac{\partial L}{\partial t'_m} =  [L, A'_m] +
B'_mL, \ \ \frac{\partial L}{\partial t''_n} =  [L, A''_n] + B''_nL.
\label{186}
\end{equation}
The simplest of them are
\begin{equation}
\frac{\partial V}{\partial t'_1} = \partial ^3 V + \partial (uV), \ \
{\bar \partial }u = 3\partial V,
\label{187}
\end{equation}
and
\begin{equation}
\frac{\partial V}{\partial t''_1} = {\bar \partial} ^3 V +
{\bar \partial} (wV), \ \
\partial w = 3{\bar \partial} V.
\label{188}
\end{equation}
The coordinated deformation (\ref{187})
and (\ref{188}) with respect to
$t = t'_1 = t''_1$
\begin{equation}
\frac{\partial V}{\partial t}
= \partial ^3 V + {\bar \partial }^3 V + \partial (uV) + {\bar
\partial }({\bar u}V), \ \
\label{189}
\end{equation}
preserves the class of real operators. Solutions to these
equations are given by
\begin{equation}
V(z,{\bar z},t',t'') =
\label{190}
\end{equation}
$$
2\partial{\bar \partial} \log \theta(U^+z + U^- {\bar z} + W^+_1 t'_1
+ \dots + W^-_1 t''_1 + \dots+ Z(D), Pr(\Gamma,\sigma)).
$$}

The inference of (\ref{190}) is analogous to that of (\ref{136}).

The Veselov--Novikov equation (\ref{190}) reduces to the KdV
equation in the case when $V$ does not depend on $y$ and
thus it is a generalization of the KdV equation different from the KP
equation.

{\bf 11.2. Infinite-dimensional Lie algebras in
soliton theory and the hierarchies BKP and DKP.}

Another definition of the KP hierarchy and its
generalizations different from given \S 10.5 was introduced in the
series of papers of Date, Jimbo, Kashiwara, and Miwa basing on the
ideas of Sato (\cite{DJKM}).  We explain it in the form close to
the methods discussed above.

Finite-zone and soliton solutions to the KP equation
are written as
\begin{equation}
u = 2\partial^2_x \log \tau(x_1,x_2,\dots).
\label{191}
\end{equation}
In the finite-zone case the function
$\tau$ is in fact a theta function which in terms of the
physical variables $x=x_1, y=x_2, t=x_3$
is written as
\begin{equation}
\tau(x_1,x_2,\dots) = \theta(U_1x_1 + U_2x_2 + \dots
+Z, \Omega)
\label{192}
\end{equation}
and the solution $\tau$ to the whole hierarchy
belongs to the space of formal series
$R = \C[[x_1,x_2,\dots]]$ in infinite number of variables.

Consider the general case. Let a group $G$ act
on a vector space $V$, let $v_0 \in V$, and let
an operator $A$ commuting with the group action be defined.
Assume that $A(v_0) = 0$.
Then
\begin{equation}
A(g\cdot v_0)=0 \ \ {\mbox {for \ \ any}}
\ \ g \in G.
\label{193}
\end{equation}

Consider the group $G = GL(\infty)$ acting diagonally on
$V = R \otimes R$. Let $v_0 = 1 \otimes
1$ and let the operator $A$ be of the shape
\begin{equation}
A \cdot \tau(x') \tau(x'') = Res_{z=0} \left( \exp \sum_{j \geq 1}
z^j (x'_j-x''_j)\right) \times
\label{194}
\end{equation}
$$
\left( \exp -\sum_{j \geq 1} \frac{z^j}{j}
\left( \frac{\partial }{\partial
x'_j}- \frac{\partial }{\partial x''_j}\right)\right)\cdot
\tau(x')\tau(x'').
$$

The operator $A$ is an example of vertex operators found by
physicists in string theory. In terms of quantum field
theory this system is expressed simpler.

Let $F$ be the $\Z$-graded fermion Fock space.
The Clifford algebra
$Cl$ generated by $\psi_j$ and $\psi^*_k$ and satisfying the relations
$\psi_j\psi^*_k+\psi^*_k\psi_j = \delta_{jk}$ and
$\psi^*_j\psi^*_k+\psi^*_k\psi^*_j =$
$\psi_j\psi_k+\psi_k\psi_j = 0$ acts on $F$.
The inclusion $gl(\infty) \rightarrow Cl$:
$E_{jk} \rightarrow \psi_j\psi^*_k$, where $E_{jk}$ is a matrix
with a single non-zero element, the unit at the intersection of the
$j$-th row with $k$-th column, is defined. The operator
$A$ on $F \otimes F$ is of the shape $\sum_j \psi_j \otimes
\psi^*_j$ and $v_0 = |0\rangle \otimes
|0\rangle$.  It is easy to check that $A$ and the action of
the algebra $gl(\infty)$ commute.
$F^{(0)}$ is transformed into $R$ by bosonization and we obtain the
system described above.

{\bf 11.2.1.} {\sl (\cite{DJKM})}
{\sl The $GL(\infty)$-orbit of $v_0$
is the set of all non-zero solutions to the equation $Av=0$.}

The detailed description for the Fock space $F$ and
bosonization is given in \cite{Kac}.

Write out (\ref{194}) in detail. Introduce the
variables $x_j = (x'_j + x''_j)/2, y_k = (x'_j - x''_j)/2$ in terms
of which (\ref{194}) is written as
\begin{equation}
\sum_{j \geq 0}
S_j(2y)S_{j+1}(-{\tilde \partial}_y) \tau(x+y)\tau(x-y) = 0
\label{195}
\end{equation}
where ${\tilde \partial}_y = (\frac{\partial }{\partial y_1},
\frac{1}{2}\frac{\partial }{\partial y_2},
\frac{1}{3}\frac{\partial }{\partial y_3}, \dots)$, and
$S_j$ are the elementary Schur polynomials (see (\ref{142})).

The equations (\ref{195}) are expressed by
{\it Hirota bi-linear equations}.
Let  $P(x_1,x_2,\dots)$ be a polynomial in finitely many
variables, then the Hirota bi-linear equation (\cite{Hir})
is as follows
\begin{equation}
P(\frac{\partial }{\partial y_1},
\frac{\partial }{\partial y_2},\dots)
f(x_1+y_1,\dots)g(x_1-y_1,\dots)|_{y=0} = 0.
\label{196}
\end{equation}
Expand (\ref{195}) into the Taylor series in  $y_1,
y_2, ...$.  Every coefficient of this expansion vanishes and this
condition is of the shape of a Hirota equation. For instance,
the vanishing of the coefficient at $y_3$ implies
\begin{equation}
(D^4_1 - 4D_1D_3 + 3D^2_2) \tau\cdot \tau = 0.
\label{197}
\end{equation}
Renormalising $x_1, x_2$, and $x_3$ in (\ref{197}),
changing $\tau$ for $\theta$ and applying the binary addition theorem
of Riemann we obtain the Dubrovin effectivization equation for the KP
equation (\ref{191}) with $d=0$.

From a realization of representation of
infinite-dimensional Lie algebras by using Clifford algebras,
Date, Jimbo, Kashiwara, and Miwa showed how in this situation
the KdV hierarchy corresponding to the algebra $A^{(1)}_1$ appears
and derived new hierarchies, natural analogues to the KP hierarchy.

In [19.IV] they introduced the {\it BKP hierarchy}
corresponding to the central extension of the algebra
$so(\infty)$, the algebra $b_{\infty}$. The BKP equations
have the shape (\ref{166}) for
\begin{equation}
B_{2n+1}
= \frac{\partial^{2n+1}}{\partial x^{2n+1}} + \sum_{k=1}^{2n-1}
b_{nk}\frac{\partial^k }{\partial x^k} , \ \ B_{2n} \equiv 0.
\label{198}
\end{equation}
Existence of a commutation representation made it possible to
find finite-zone solutions by using Baker--Akhieser functions
[19.V].  These solutions are of the shape
(\ref{191}--\ref{192}). Like the solutions to the
Veselov--Novikov equations, they are expressed in the Prym theta
functions of double coverings with two branch points.
The first equations of the BKP hierarchy in the Hirota shape
are as follows
\begin{equation}
(D^6_1 - 5D^3_1D_3 - 5D^2_3 + 9 D_1D_5)
\tau\cdot \tau = 0,
\label{199}
\end{equation}
$$
(D^8_1 + 7D^5_1D_3 -
35D^2_1D^2_3 -21 D^3_1D_5 - 42D_3D_5 + 90D_1D_7)\tau\cdot \tau=0.
$$

The CKP hierarchy relates to the algebra $sp(\infty)$
[19.VI] and also admits a commutation representation. Its
finite-zone solutions are constructed from double coverings of Riemann
surfaces with $2n > 2$ branch points.  The Prym varieties of
such coverings are not principally-polarised and for the
moment explicit theta functional formulae for finite-zone solutions
are not written.

In [19.IV] the $m$-component BKP hierarchies
related to central extensions of the algebra
$so(m\cdot \infty)$ were also introduced.  Solutions to
equations of this hierarchy are functions of infinitely many
variables divided into $m$ classes of the shape
$(x^{(1)}_1, x^{(1)}_2, \dots)$,..., $(x^{(m)}_1, x^{(m)}_2,\dots)$.
The simplest example is the hierarchy corresponding to $m=2$.
In this case the first non-trivial equations are
\begin{equation}
{\hat D}(D^3_1-D_3)\tau\cdot \tau = 0,
\label{200}
\end{equation}
$$
{\hat D}(D_1^5  + 5D_3 D_1^2 - 6D_5) \tau\cdot \tau = 0
$$
where $D_j = D^{(1)}_j$ and ${\hat D} = D^{(2)}_1$.  In
[19.VII] this hierarchy is also called BKP-II. In the
sequel they became to call it the {\it DKP hierarchy} (\cite{JM}).  No
commutation representation for the DKP equations was found and thus
its solutions were not discussed.

Notice that another approach to constructing hierarchies
of non-linear equations related to Prym theta functions was
considered in \cite{Nat}.

The ideas of papers \cite{DJKM} in the sequel were developed in
different directions among which we mention two ones.
Theorem 11.2.1 going back to Sato may be treated simply as follows:
tau-functions of the KP hierarchy are
parametrised by points of infinite-dimensional Grassmann space,
$\C[[x_1,x_2,\dots]]/GL(\infty)$.  Investigation of cell
decomposition, of this space, related to decomposition of all
possible solutions into different classes (finite-zone, rational,
soliton, etc.) was undertaken in \cite{SW}.
The other interesting problem on canonical correspondence to
every loop group a hierarchy of soliton equations was
investigated in \cite{KW}.

{\bf 11.3. The Landau--Lifschitz equation.}

In \cite{DJKM2} the methods of \cite{DJKM} were applied to the
Landau--Lifschitz equation (the LL equation)
known in solid state physics:
\begin{equation} S_t=[S \times S_{xx}] + [S\times JS]
\label{201}
\end{equation}
where $S$ is a three-dimensional vector,
$S_1^2+S_2^2+S_3^2=1$, and $J=diag (J_1,J_2,J_3)$.

Its integrability by methods of soliton theory was established
by Borovik and Sklyanin by the end of the 70s. Hirota equations
for (\ref{201}) were found in \cite{Hir2} on base of the ansatz
introduced by the physicists Bogdan and Kovalev:
$$
S_1=(fg^*+gf^*)/(ff^*+gg^*), $$ $$ S_2=-i(gf^*-fg^*)/(ff^*+gg^*),
$$
$$
S_3=(ff^*-gg^*)/(ff^*+gg^*).
$$
These equations are
\begin{equation}
D_1(f\cdot f^*+g\cdot g^*)=0,
\label{202}
\end{equation}
\begin{equation}
(D_2-D_1^2)(f\cdot f^*-g\cdot g^*)=0,
\label{203}
\end{equation}
\begin{equation}
(D_2-D_1^2+\lambda)f\cdot g^*+\mu g\cdot f^*=0,
\label{204}
\end{equation}
and
\begin{equation}
(D_2-D_1+\lambda)g\cdot f^*+\mu f\cdot g^*=0
\label{205}
\end{equation}
where
$$
\lambda+\mu=J_3-J_1, \lambda-\mu=J_3-J_2.
$$

As it was pointed out in \cite{DJKM2}
finite-zone solutions to the LL equation are expressed
in terms of theta functions of the Prym varieties of double
coverings of a special type. We describe it in \S 12.3.

\vskip2.5mm

{\bf \S 12. Methods of the finite-zone integration in the theory of
Prym map}

\vskip2.5mm

{\bf 12.1. Analogues of the Torelli theorem for Prym maps.}

We considered above two Prym maps corresponding to double
coverings unramified (\ref{116})  and ramified at two branch points
(\ref{99}).

Generally, analogues of the Torelli theorem are not valid for
either of them.

For unramified coverings the situation is rather clear.

{\bf 12.1.1.}
{\sl 1) The Prym map $Pr : {\cal DU}_g
\rightarrow {\cal A}_{g-1}$ has maximal rank at a generic point.
For $g \leq 5$ the rank equals
$g(g-1)/2 < \dim {\cal DU}_g$ and $\overline{Pr({\cal DU})_g)} =
{\cal A}_{g-1}$, and for $g \geq 6$ the rank equals $3g-3 = \dim
{\cal DU}_g$.

2) (Donagi and Smith (\cite{DonS})) For $g=6$ the Prym map is
a $27$-sheeted covering at a generic point and $\overline{Pr({\cal
DU}_6)} = {\cal A}_5$.

3) (Kanev (\cite{Kan}), Friedman and Smith (\cite{FS})) For
$g \geq 7$ the Prym map is injective at a generic point.}

The Donagi tetragonal construction shows that for any $g$
there exists principally-polarised Abelian varieties being
the Prym varieties of several non-equivalent coverings.

The case of coverings ramified at a pair of points was
not practically examined. By soliton methods it was proved
that

{\bf 12.1.2.} {\sl (\cite{T90})}
{\sl The Prym map $Pr : {\cal DR}_g \rightarrow
{\cal A}_g$ has maximal rank at a generic point.}

This proposition was proved as an auxiliary one
while studying an analogue of the Novikov conjecture in
the Riemann--Schottky problem (see \S 12.4).

{\bf 12.2. Effectivization of theta functional formulae
for solutions to the Veselov--Novikov equation.}

Effectivization of theta functional formulae for finite-zone
solutions to the Veselov--Novikov equation was
developed in \cite{T90} by using the methods of (\cite{Dubr1}).
The results were announced by the author in
Soviet Math. Dokl. 32 (1985), 843--846.

{\bf 12.2.1.} {\sl (\cite{T90})}
{\sl Let $\theta$ be a theta function of an irreducible
principally-polarised Abelian variety of dimension $g$.
If for every $Z \in \C^g$ the formula
\begin{equation}
{\cal U}_Z = 2\partial {\bar \partial } \log \theta(Uz + V{\bar z} + Wt +
Z) + C
\label{206}
\end{equation}
gives a solution to the ``holomorphic'' half-part of the
Veselov--Novikov equation
\begin{equation}
\frac{\partial {\cal U}_Z}{\partial t}  = \partial^3
{\cal U}_Z + \partial (v {\cal U}_Z),
\label{207}
\end{equation}
$$
v = 6\partial^2 \log \theta(Uz + V{\bar z}+ Wt + Z) + d,
$$
then the relations
\begin{equation}
(D_2D_3  - 4D_1^3D_2 - dD_1D_2 - 3CD_1^2 - a) {\hat \theta}[n,0] = 0
\label{208}
\end{equation}
hold,
with $a, d, C$ constants and $D_1 =
\partial_U, D_2 = \partial_V, D_3 = \partial_W$.
For the ``anti-holomorphic'' half-part the relations have the same
shape with $U$ and $V$ permuted.}

Notice that for $a=d=C=0$ (\ref{208})
transforms into the Hirota equation for the first equation
of the BKP hierarchy (\ref{200}).

From every double ramified covering with two branch points
a solution to the Veselov--Novikov equation is constructed.
Hence, for every Prym variety there exist vectors
$U, V \in \C^g \setminus \{0\}, W \in \C^g$ and constants
$a, d, C$ such that (\ref{208}) hold.
This led Novikov to an analogue of his conjecture for Prym
varieties which we introduce here and discuss in \S 12.4.

{\bf 12.2.2. Analogue of the Novikov conjecture.}
{\sl An irreducible principally-polarised Abelian variety
is the Prym variety of a double covering with two branch
points if and only if (\ref{208}) are solvable
with $U, V \neq 0$.}

The relations (\ref{208}) have the following property.

{\bf 12.2.3.} {\sl (\cite{T90})}
{\sl If an Abelian variety is irreducible and vectors
$U$ and $V$ are linearly independent, then by
(\ref{208}) the vector $W$ and the constants $a, c, d$
are reconstructed uniquely up to the transformations
\begin{equation}
U \rightarrow \lambda U, V \rightarrow \mu V,
W \rightarrow \lambda^3 W + \alpha U,
\label{209}
\end{equation}
$$
a \rightarrow \lambda^3\mu
a, C \rightarrow \lambda\mu C, d \rightarrow \lambda^2d+\alpha.
$$}

In particular, for  $g=2$ we obtain a direct procedure for
constructing solutions to the Veselov--Novikov equation
from an arbitrary irreducible principally-polarised Abelian variety
and an arbitrary pair of linearly independent vectors.
The constant $C$ is reconstructed by two different ways, i.e., by
the effectivization equations for ``holomorphic'' and
``anti-holomorphic'' half-parts of the Veselov--Novikov equation.
The condition of coincidence of solutions implies non-trivial
relations for theta constants in the case $g=2$.

{\bf 12.2.4.} {\sl (\cite{T87})}
{\sl Let $\theta$ be a theta function of two variables
and let
\newline
$(a^{11}_n \ a^{12}_n \ a^{22}_n \ a_n)$ be
the inverse matrix of
$({\hat \theta}_{11}[n,0] \ {\hat
\theta}_{12}[n,0] \ {\hat \theta}_{22}[n,0] \ {\hat \theta}[n,0])$:
$$
a^{pq}_n{\hat \theta}_{kl}[n,0] = \delta^p_k\delta^q_l, a_n{\hat
\theta}_{kl}[n,0] = a^{pq}_n{\hat \theta}[n,0] = 0, a_n{\hat
\theta}[n,0] = 1.
$$
Then the relations
$$
\sum_n (a^{km}_n {\hat \theta}_{mmmm}[n,0] + 2a^{kk}_n {\hat
\theta}_{kmmm}[n,0]) = 0,
$$
\begin{equation}
\sum_n (a^{km}_n {\hat \theta}_{kmmm}[n,0] -
a^{mm}_n {\hat \theta}_{mmmm}[n,0]) + 3a^{kk}_n
{\hat \theta}_{kkmm}[n,0]) = 0,
\label{210}
\end{equation}
$$
\sum_n (a^{11}_n {\hat \theta}_{1112}[n,0] -
a^{22}_n {\hat \theta}_{1222}[n,0]) = 0
$$
hold where $1 \leq k,m \leq 2$, $n \in 1/2
(\Z^2/2\Z^2)$.}

{\bf 12.3. Quadrisecant formulae and the Veselov--Novikov, BKP,
and Landau--Lifschitz equations.}

Denote by ${\vec \theta}(z)$ the vector
$(\theta[n_1,0](z,2\Omega),\dots,\theta[n_r,0](z,2\Omega))$
and denote by $u \wedge v$ the product of vectors
$u,v \in \C^{2^g}$ in the Grassmann algebra.

Let $\Gamma\rightarrow \Gamma/\sigma$ be a double covering with
branch points at $q_1$ and $q_2$.
By the Fay quadrisecant formula 6.3.1 for an arbitrary
quadruple $p_1, \dots, p_4 \in \Gamma$
the identity
\begin{equation}
{\vec\theta}(p_1+p_2+p_3+p_4) \wedge {\vec\theta}(p_1+p_2-p_3-p_4)
\wedge
\label{211}
\end{equation}
$$
{\vec\theta}(p_1+p_3-p_2-p_4)
\wedge {\vec\theta}(p_1+p_4-p_2-p_3) = 0
$$
holds where as in \S 6 we mean
$A_{Pr}(p_j)$ by $p_j$.

Introduce local coordinates $k^{-1}_j$ in neighbourhoods
of points $q_j$ and consider the expansions of the Abel--Prym mapping
into series
\begin{equation}
k_1^{-1} \rightarrow A_{Pr}(q_1) + \sum_m \frac{U_m}{k_1^m}, \ \
k_2^{-1} \rightarrow A_{Pr}(q_2) + \sum_n \frac{V_n}{k_2^n}.
\label{212}
\end{equation}
If the involution $\sigma$ inverts parameters then
$U_{2m}=V_{2n}=0$. Denote by $D_j$ the derivative in
the direction $U_j$ and denote by ${\hat D}_k$ the
derivative in the direction  $V_k$. Also introduce the
operator $T(k)$
\begin{equation}
T(k_1^{-1}) =
\exp{(\sum_{m=1}^{\infty} \frac{D_m}{k^m})}= 1+ \sum_{n\geq 1}
\frac{\Delta_n}{k_1^n}.
\label{213}
\end{equation}
Take $q_1$ as the initial point of the Abel--Prym mapping,
$A_{Pr}(q_1) = 0$. Since $q_2$ is a fixed point of the involution
$\sigma : \Gamma \rightarrow \Gamma$ then we may assume
that $A_{Pr}(q_2) = 0$.

Consider the degeneration (\ref{211}) corresponding to the case
$p_4 = q_1$ and $p_2, p_3 , p_4 \rightarrow q_1$.
We omit details and just give the final result.

{\bf 12.3.1.} {\sl (\cite{T91})}
{\sl The equations
\begin{equation}
T(k_1^{-1}) \circ
({\vec \theta}(z) \wedge
\Delta_1 {\vec \theta}(z) \wedge
\Delta^2_1 {\vec \theta}(z) \wedge
(\Delta_1\Delta_2 - \Delta_3){\vec \theta}(z))|_{z=0} \equiv 0
\label{214}
\end{equation}
are obtained from the Fay quadrisecant formula
(\ref{211}) by the degeneration $p_1, p_2$, $p_3$, $p_4
\rightarrow q_1$.
By expanding the left-hand side of (\ref{214}) into the series in
the powers of $k_1^{-1}$ we obtain a hierarchy of equations
in the Hirota shape. The principal members of this hierarchy,
possibly after change of a local parameter, are written as
\begin{equation}
(\Delta_n\Delta_1 - \Delta_{n-1}\Delta^2_1 +
\Delta_{n-2}(\Delta_1\Delta_2-\Delta_3))\theta\cdot \theta=0.
\label{215}
\end{equation}
Its members for $n=5,7$ coincide with (\ref{199}) and
the hierarchy (\ref{214}) coincides with the BKP hierarchy.}

{\bf 12.3.2.} {\sl (\cite{T91})}
{\sl The equations
\begin{equation}
T(k_1^{-1}) \circ
({\vec \theta}(z) \wedge
D_1 {\vec \theta}(z) \wedge
{\hat D}_1 {\vec \theta}(z) \wedge
D_1{\hat D}_1{\vec \theta}(z))|_{z=0} \equiv 0,
\label{216}
\end{equation}
are obtained from the Fay quadrisecant formula (\ref{211})
by the degeneration $p_1, p_2$, $p_3
\rightarrow q_1, p_4 \rightarrow q_2$.
By expanding the left-hand side of (\ref{216}) into the series in
the powers of $k_1^{-1}$ we obtain a hierarchy of equations
in the Hirota shape. The principal members of this hierarchy,
possibly after change of local parameters, are written as
\begin{equation}
(\Delta_n{\hat D}_1 - \Delta_{n-1}\Delta_1{\hat D}_1)
\theta\cdot \theta = 0.
\label{217}
\end{equation}
Its members for $n=3,5$ are of the shape (\ref{200})
and the hierarchy (\ref{216}) coincides with the hierarchy of
``holomorphic'' half-parts of the Veselov--Novikov equations or
with the DKP hierarchy.}

These propositions explain from the algebro-geometric point of
view the existence of a pair of hierarchies related to double
coverings of algebraic curves with two branch points
and close to the Kadomtsev--Petviashvili hierarchy due to their
properties.

The equivalence of a part of the Veselov--Novikov hierarchy
and the DKP hierarchy indicated above explains why there is
no commutation representation for the DKP equations by
scalar differential operators.

We do not give explicit formula for finite-zone solutions to the DKP
hierarchy. Notice only that they can be easily obtained from
(\ref{212},\ref{214},\ref{215}).

Now consider the Beauville--Debarre quadrisecant identity 7.2.2.
We apply it to unramified coverings of the class related to the
Landau--LIfschitz equation.

Let $\Gamma\rightarrow E$ be a double covering with $4g$ branch
points of a Riemann surface of genus $1$ and let the branch set be
invariant under the translation by a half-period:
$u\rightarrow u+\alpha$. Denote by $\sigma$ the involution of
$\Gamma$ generated by this translation and denote by $\omega$
the involution permuting the branches of the covering
$\Gamma\rightarrow E$.  Evidently,
$\Gamma\rightarrow\Gamma/\sigma$ is an unramified double
covering and the involutions
$\sigma$ and  $\omega$ commute:
$\sigma\omega=\omega\sigma$.  We call the covering
$\Gamma\rightarrow \Gamma/\sigma$ an {\it LL-covering}.

As it was first indicated in \cite{DJKM2} finite-zone solutions to
the LL equations are constructed from such coverings. Explicit
formulae in terms of Prym theta functions we derived in \cite{Bob}.
Choose a basis for Prym differentials  $u_k$ such that
$\sigma^*(u_k)=\omega^*(u_k)=-u_k$.  Fix a point
$p_0$, a branch point of the covering
$\Gamma\rightarrow \Gamma/\omega$, and introduce the vector
$$
\mu_k = \int^{p_0}_{\sigma(p_0)} u_k.
$$
It is evident that
$2\mu \equiv 0$ in $Pr(\Gamma,\sigma)$.
Then we have
$$
A_{Pr}(p) =
\frac{1}{2}\int^p_{\sigma(p)} u = \frac{1}{2}\left( \int^p_{p_0} u +
\mu + \int^{\sigma(p_0)}_{\sigma(p)} u \right) = \int^p_{p_0} u +
\frac{\mu}{2}.
$$

Consider a quadruple $p_1, p_2, p_3, p_4 \in \Gamma$ and denote by
$\delta$ the vector
$$
\delta_j = \int^{p_1}_{\omega(p_1)} u_j.
$$
Introduce a local parameter $\varepsilon$
in a neighbourhood of $p=p_1$ and expand the Abel--Prym mapping
into the series in the powers of $\varepsilon$
$$
A_{Pr}(\varepsilon ) = A_{Pr}(p) +
U_1\varepsilon +U_2\varepsilon^2 + \dots.
$$
Denote by $D_j$ the derivative in the direction $U_j$.

The Beauville--Debarre quadrisecant formula is of the shape
(\ref{211}) where for brevity we denote $A_{Pr}(p)$ by $p$ in arguments
of theta functions.

Consider the degeneration of the quadrisecant formula as
$p_2, p_3 \rightarrow \sigma\omega(p)$, $p_4 \rightarrow\omega(p)$.
In this event
$$
p_2, p_3 \rightarrow p_1 + \mu,  \ \ p_4
\rightarrow -p_1-\mu,
$$
$$
(p_1 + \dots + p_4), (p_1 + p_2 - p_3 -
p_4), (p_1 + p_2 - p_2 - p_4) \rightarrow \delta,
$$
$$
(p_1 + p_4
-p_2 - p_3) \rightarrow -\delta-2\mu
$$
in the Prym variety.

The simplest relations followed from this degeneration
are
\begin{equation}
(D_2 - D_1^2 + \alpha) {\vec \theta}(\delta) +
\beta {\vec \theta}(2\mu+\delta) = 0.
\label{218}
\end{equation}
Transforming (\ref{218}) into the Hirota equations
by using the binary addition theorem of Riemann we
obtain (\ref{204}) and  (\ref{205}) for
$f=\theta(z+\mu+\delta)$, $f^*\theta(z)$, $g=-i\theta(z+\delta)$, and
$g^*=-i\theta(z+\mu)$ where $\theta$ is the Prym theta function of
the covering $\Gamma\rightarrow \Gamma/\sigma$
((\ref{204}) is transformed into (\ref{205}) by the translation
$z \rightarrow z + \mu$).

For this choice of
$f,f^*,g,g^*$ (\ref{202}) is equivalent to
$$
D_1\left({\vec \theta}(\mu+\delta)+{\vec \theta}(\mu-\delta)\right)=0
$$
and  (\ref{203}) is equivalent to
$$
(D_2-D_1^2)\left({\vec
\theta}(\mu+\delta)-{\vec \theta}(\mu-\delta)\right)=0,
$$
i.e., (\ref{202}) and (\ref{203}) are trivial consequences
of the symmetry of the function ${\vec \theta}$
(${\vec \theta}(z)={\vec \theta}(-z)$).  As a result we obtain the
following.

{\bf 12.3.3.} {\sl (\cite{T93})}
{\sl The equations of the LL hierarchy
are obtained from the quadrisecant formula for LL-coverings
in the limit as $p_2, p_3 \rightarrow \sigma\omega(p_1),
p_4 \rightarrow \omega(p_1)$.}

Notice that the class of LL-coverings is rather small
comparatively with all unramified double coverings. Indeed,
LL-coverings with $g$-dimensional Prym variety are defined locally
by $2g$ parameters:  one parameter defines the conformal class of
a torus $E$ and  $2g-1$ parameters define  branch points taken up to
translation on the torus. All coverings with
$g$-dimensional Prym variety are defined by $3g$
($= \dim {\cal DU}_{g+1}$) parameters.

Note that degenerations of quadrisecant formulae were
considered in \cite{G2}.  So far , possible relation of results
obtained in \cite{G2} to soliton equations were not discussed.

{\bf 12.4. Analogue of the Riemann--Schottky problem for Prym
varieties.}

In the spirit of the Novikov conjecture the Veselov--Novikov,
BKP, and Landau--Lifschitz equations can be considered as
characterizing Abelian varieties in which theta functions
these equations are integrated. Of a special interest are the first
two hierarchies with their solutions constructed from
all Prym varieties of coverings with two branch points.
In this connection we mentioned above an analogue of
the Novikov conjecture 12.2.2.

An analogue of the Dubrovin theorem 9.2.3
was obtained in \cite{T87,T90}:

{\bf 12.4.1.}
{\sl The solvability condition on
(\ref{208}), i.e., on the effectivization equations for finite-zone
solutions to the Veselov--Novikov equation, distinguishes
a subvariety in ${\cal A}_g$ containing
$\overline{Pr({\cal DR}_g})$ as one of components.}

That is, the Veselov--Novikov equation locally solves the
Riemann-Schottky problem for the Prym varieties of coverings with two
branch points. In view of 12.1.2 the proof of 12.4.1
reduces to computing the dimension of this subvariety. In proofs of
Theorems 9.2.3, 10.1.2, and 10.2.1, solving
the Riemann--Schottky problem locally, this part was played by the
Torelli theorem.

As it turns out using of higher equations of the Veselov--Novikov
hierarchy as well as using of the BKP hierarchy for
rejecting other components met considerable difficulties.

First, potential two-dimensional Schr\"odinger operators,
finite-zone on the zero energy level, are constructed also
by singular surfaces:

1) let the spectral surface $\Gamma$
of a Schr\"odinger operator also has
double points $q_1,\dots,q_k$;

2)  all singular points $q_1, \dots, q_n$
are fixed points of the involution $\sigma$ which does not
permute the branches in their neighbourhoods
(\cite{VN1}).

\noindent
The quotient surface $\Gamma/\sigma$ also has
$k$ double points. The Jacobi varieties of the surfaces
$\Gamma$ and $\Gamma/\sigma$ are non-compact but the Prym
variety can be treated up to isogeny
as the quotient space $J(\Gamma)/J(\Gamma/\sigma)$
and it is compact and principally-polarised.

Second, at the soliton proof of Theorem 10.3.4
we told that the KP-flows span the tangent space to the Jacobi
variety of the surface $\Gamma$. For the BKP hierarchy this does
not valid and theta functional solution can be concentrated on the
non-full-dimensional winding of a principally-polarised Abelian
variety.

This explains the following statement of the theorem of Shiota
(\cite{Shi2,Shi3}).

{\bf 12.4.2.}
{\sl If for any $Z \in \C^g$ the function
$$
{\cal U}_Z = 2\partial ^2\log \theta(Ux+V_1t_1 + \dots +Z,\Omega)
$$
is a solution to the BKP hierarchy and 
$\C^g$ is the linear span of the vectors
$U, V_1, \dots$, then the principally-polarised Abelian
variety
$\C^g/\{\Z^g+\Omega\Z^g\}$ is the Prym variety of a double
covering $\Gamma\rightarrow \Gamma/\sigma$
with two non-singular branch points and finitely many double
points at which the involution $\sigma$ does not permute the branches
of $\Gamma$.}

An analogous result can be obtained for the
Veselov--Novikov hierarchy (\cite{Shi3}).

The analogue of the Riemann--Schottky problem for coverings with two
branch points was not discussed from a geometric point of view
in contrast with the case of unramified coverings.

The most close in spirit to the soliton theory is the
following result of Beauville and Debarre
(\cite{BD1}), an analogue of their Theorem 10.6.2.

{\bf 12.4.3.}
{\sl The condition of existence of quadrisecant distinguishes
a subvariety in ${\cal A}_g$
containing $\overline{Pr({\cal DU}_{g+1})}$,
a locus of Prym varieties of double unramified coverings, as one of
components.}

So far , the relation of this result to the theory of soliton
equations is not clear. It seems that the
Landau--Lifschitz equation can be used for characterizing Prym varieties
of LL-coverings.

Recalling other approaches to the Riemann--Schottky
problem we note the investigation of $Sing \Theta$ carried out
in \cite{Mm4,Sho1}.  In particular, it enables us to understand and
describe the difference between  Prym and Jacobi varieties.
Other geometric approaches were also considered
(\cite{Rad}) but we do not dwell on them.
Notice only that the Donagi conjecture on the type of
non-single-valuedness of the Prym map (see the survey \cite{Sho2}) and
the analogue of the Riemann--Schottky problem for Prym varieties are
still in question.

\vskip2.5mm

\begin{center}
{\bf Final remarks}
\end{center}

\vskip2.5mm

In conclusion we indicate some problems which are left beyond the
present survey.

1. It seems that the first application of soliton theory
to algebraic geometry was the following result of Dubrovin and
Novikov (\cite{DN}). Let $V_g$ be the moduli space of hyperelliptic
surfaces of genus $g$. Consider the bundle
$M_g \rightarrow V_g$ with a fibre $J(\Gamma), \Gamma \in
V_g$.  Construct a $2g+2$-sheeted covering over $M_g$
corresponding to fixing one of the $2g+2$ branch points of
the hyperelliptic covering
$\Gamma\rightarrow \C P^1$.  The space ${\tilde M}_g$
is called the {\it complete moduli space
of hyperellitic Jacobi varieties}.
In (\cite{DN}) it was proved that {\sl ${\tilde M}_g$ is
a rational variety}.

The finite-zone theory of equations related to
hyperelliptic surfaces consequently led to many profound
results (see their exposition in \cite{Mm2}).

2. Recently Nakayashiki, developing the ideas of Sato, has shown
that there exists a natural definition of integrable non-linear
equations expressed as commutation conditions on matrix differential
operators (\cite{Nak1,Nak2}). These equations define deformations of
(Baker--Akhieser) vector bundles on arbitrary principally-polarised
Abelian varieties and thus are not related to the geometry of Riemann
surfaces. For the moment explicit formulae for equations and their
solutions are neglected.

3. The trisecant identity is a far-reached generalization of
the addition theorem for elliptic functions
(this is discussed in \cite{Fay1}).  There exist another
possibilities of generalizing this theorem , i.e., the
addition theorems for hyperelliptic functions of large genus
(\cite{Baker}). Their relation to soliton equations has been
observed recently and now it is under investigation (\cite{BEL}).

4. The papers of Mulase \cite{Mul2} and Li and Mulase
\cite{LiMul,LiMul2} are devoted to category generalizations of the
Burchnall--Chaundy--Krichever correspondence.
In \cite{Mul2} on the language of category theory the part
of the paper of Krichever and Novikov \cite{KN2} where the KP
hierarchy was realised as deformations of
framed semi-stable vector bundles of rank $l \geq 1$ with
$c_1=lg$ over non-singular Riemann surfaces of genus
$g$ is exposed.  Theta functional identities and
formulae are not discussed in \cite{Mul2,LiMul,LiMul2}.
It seems interesting to relate with theta functional formulae the
results of \cite{LiMul} where the Prym varieties of $n$-sheeted
coverings were ineffectively characterised in terms of the $n$-component
KP hierarchies.

5. We did not address at all the relation of Riemann
surfaces with finite-dimensional integrable systems having
commutation representations (an up-to-date survey is given
in \cite{DM}) and with difference soliton equations (a survey of
these papers is contained in \cite{DKN2}).

\end{document}